\documentclass[preprint,10pt]{elsarticle}
\usepackage{hyperref}
\usepackage{amsmath,amssymb}

\usepackage{graphicx}
\usepackage{graphics}
\usepackage {graphicx,fancyhdr}
\usepackage{graphics, color}
\usepackage{flafter}
\numberwithin{equation}{section}
\usepackage{multirow}
\usepackage{mathrsfs}
\usepackage{subfigure}
\usepackage{tabularx}
\usepackage{overpic}
\usepackage{stmaryrd}
\usepackage{color}

\begin{document}
	
\begin{frontmatter}
\tnotetext[label1]{This work is partially supported by Key Program Special Fund in XJTLU (KSF-E-50, KSF-P-02) and XJTLU Research Development Funding (RDF-19-01-15). The research of Eric Chung is partially supported by the Hong Kong RGC General Research Fund (Project numbers 14304719 and 14302018) and CUHK Faculty of Science Direct Grant 2019-20.}
	
\title{Generalized multiscale approximation of a multipoint flux mixed finite element method for Darcy-Forchheimer model \tnoteref{label1}}
\author[lab1]{Zhengkang He}
\author[lab2]{Eric T. Chung}
\author[lab1,lab3]{Jie Chen\corref{cor1}}
\author[lab1,lab4]{Zhangxin Chen\corref{cor2}}
	
\cortext[cor1]{Corresponding author. Jie Chen E-mail address: jie.chen01@xjtlu.edu.cn}
\cortext[cor2]{Corresponding author. Zhangxin Chen E-mail address: zhachen@ucalgary.ca}

\cortext[lab1]{Zhengkang He E-mail address:  hzk2abc@163.com}
\cortext[lab2]{Eric T.Chung E-mail address:  tschung@math.cuhk.edu.hk.}

\address[lab1]{School of Mathematics and Statistics, Xi'an Jiaotong University, Xi'an, 710049, China}
\address[lab2]{Department of Mathematics, The Chinese University of Hong Kong (CUHK), Hong Kong SAR}
\address[lab3]{Department of Mathematical Sciences, Xi'an Jiaotong-Liverpool University, Suzhou, 215123, China}
\address[lab4]{Department of Chemical $\&$ Petroleum Engineering, Schulich School of Engineering, University of Calgary, 2500 University Drive N.W., Calgary, Alberta T2N 1N4, Canada}
	
\begin{abstract}
In this paper, we propose a multiscale method for the Darcy-Forchheimer model in highly heterogeneous porous media. The problem is solved in the framework of generalized multiscale finite element methods (GMsFEM) combined with a multipoint flux mixed finite element (MFMFE) method. 
We consider the MFMFE method that utilizes the lowest order Brezzi-Douglas-Marini ($\textrm{BDM}_1$) mixed finite element spaces for the velocity and pressure approximation. The symmetric trapezoidal quadrature rule is employed for the integration of bilinear forms relating to the velocity variables so that the local velocity elimination is allowed and leads to a cell-centered system for the pressure. 
We construct multiscale space for the pressure and solve the problem on the coarse grid following the GMsFEM framework. In the offline stage, we construct local snapshot spaces and perform spectral decompositions to get the offline space with a smaller dimension. In the online stage, we use the Newton iterative algorithm to solve the nonlinear problem and obtain the offline solution, which reduces the iteration times greatly comparing to the standard Picard iteration. Based on the offline space and offline solution, we calculate online basis functions which contain important global information to enrich the multiscale space iteratively. The online basis functions are efficient and accurate to reduce relative errors substantially. Numerical examples are provided to highlight the performance of the proposed multiscale method.
\end{abstract}
	
	
\end{frontmatter}
\section{Introduction}
In many porous-media flow applications, such as petroleum recovery and groundwater resource management, flow velocities are typically low and Darcy's law is usually used to describe the linear relationship between the velocity and the gradient of pressure. However, there are cases that flow velocities are relatively high, the relationship between
the velocity and the pressure gradient becomes nonlinear such that Darcy's law no longer holds, which is observed by Forchheimer and others. The nonlinear relationship is described by Forchheimer' law (Darcy-Forchheimer equation) which is a corrected formula of Darcy's law by supplementing a quadratic nonlinear inertial term.

From the numerical perspective, there have been many methods developed for solving the Darcy-Forchheimer model in porous media. Park \cite{park2005mixed} studied a semi-discrete mixed finite element method for generalized Forchheimer flow. In \cite{girault2008numerical}, Girault et al. employed piecewise constant elements and Crouziex-Raviart elements for the approximation of velocity and pressure, respectively. Pan et al. \cite{pan2012mixed} presented a different mixed formulation that makes use of Raviart-Thomas mixed elements or Brezzi-Douglas-Marini mixed elements to discretize the velocity and pressure variables. Rui et al. \cite{rui2012block,rui2015block,rui2017block} constructed block-centered finite difference methods. Wang and Rui \cite{wang2015stabilized} introduced a stabilized mixed finite element method using Crouzeix-Raviart elements. Xu et al. \cite{xu2017multipoint} developed an MFMFE method for the compressible Darcy-Forchheimer model. There is also a nonlinear multigrid method constructed in \cite{huang2018multigrid}, two-level methods established in \cite{rui2015two,fairag2020two}, and variational multiscale interpolating element-free Galerkin method developed in \cite{zhang2020variational} for the nonlinear Darcy-Forchheimer model.

In the real world, the geological porous media are generally governed by coefficients with high heterogeneities and complex spatial distributions. Solving these problems directly on the fine grid will result in large-scale discrete systems which are challenging to deal with. Model reduction techniques are required to reduce computational complexity. Spiridonov et al. \cite{spiridonov2019mixed} utilized the mixed generalized multiscale finite element method (mixed GMsFEM) to approximate the Darcy-Forchheimer model on the coarse grid. The mixed GMsFEM is originally developed in \cite{chung2015mixed} for Darcy's flow in heterogeneous media, the multiscale basis functions for velocity are constructed following the GMsFEM framework \cite{efendiev2013generalized,chung2014adaptive,efendiev2014generalized,chung2015residual,chan2016adaptive,chung2016adaptive} which generalizes the multiscale finite element method (MsFEM) \cite{hou1997multiscale} by enriching the coarse-grid space systematically with additional multiscale basis functions that can help to reduce the error efficiently and substantially. Recently, Chen et al. \cite{chen2020generalized} proposed another form of the mixed GMsFEM for Darcy's law, where pressure is approximated in a multiscale function space between fine-grid space and coarse-grid space, trapezoidal quadrature rule is used for local velocity elimination on rectangular meshes and velocity is solved directly in the fine-grid space.

In this paper, we develop an efficient algorithm to construct a multiscale solution on the coarse grid for the Darcy-Forchheimer model in heterogeneous porous media following the framework proposed in \cite{chen2020generalized}. The algorithm is based on the combination of the GMsFEM and MFMFE method. In the MFMFE methods \cite{wheeler2006multipoint,wheeler2012multipoint,ingram2010multipoint,klausen2006robust,hou2014numerical}, appropriate mixed finite element spaces and suitable quadrature rules are employed, which allow for local velocity elimination and lead to a cell-centered system for the pressure. Here, we consider an MFMFE method that has been studied in \cite{xu2017multipoint} for the discretization of Darcy-Forchheimer model on the fine grid meshes (reference solution) composed of simplices and perturbed parallelograms, where $\textrm{BDM}_1$ mixed finite element spaces are used for the approximation of velocity and pressure variables, and symmetric trapezoidal quadrature rule is employed for the integration of bilinear forms relating to velocity variables. The resulting mass matrix for velocity in the discrete system is block diagonal, symmetric and positive definite, which can be inverted straightly, i.e. the velocity can be solved in the fine-grid space explicitly. For the coarse-grid approximation, we follow the GMsFEM framework to calculate the multiscale basis functions for the pressure on the coarse grid. In the offline stage, we begin with the construction of the local snapshot space per coarse element by solving a series of local problems numerically, then we obtain the smaller dimensional offline space through the spectral decompositions in each local snapshot space. In the online stage, firstly, we exploit the derived offline space to solve the nonlinear problem on the coarse grid and find out the offline solution. Different from \cite{spiridonov2019mixed}, Newton iterative algorithm is used to handle the nonlinear term, which will result in much fewer iterations than the use of Picard iterative algorithm when the nonlinearity is strong. The offline solution has a good approximation of the fine-grid solution. Secondly, in order to achieve higher accuracy, based on the offline space and solution, we perform enrichments of the multiscale space with the addition of online basis functions iteratively. The online basis functions, which have been well studied in \cite{chung2015residual,chan2016adaptive,chung2016adaptive,chung2017online}, can capture important global information and are capable of realizing a substantial error reduction of the multiscale solution.

We organize the rest of the paper as follows. In section 2, we introduce the Darcy-Forchheimer model, the corresponding weak formulation and the fine-grid discretization of the problem by use of an MFMFE method. In section 3, we first construct local snapshot spaces and the offline space for approximating the pressure, then we enrich the multiscale space by adding online basis functions based on the offline space and the offline solution to improve the accuracy. In section 4, some numerical examples are presented. Finally, we give some conclusions in section 5.
\section{Darcy-Forchheimer model, weak formulation and fine-grid approximation}
\subsection{Darcy-Forchheimer model and weak formulation}
Let $\Omega$ be a bounded and simply connected porous-media domain in $\mathbb{R}^2$ with a Lipschitz continuous boundary $\partial\Omega$. Darcy-Forchheimer model that describes the single-phase flow in $\Omega$ is the coupling of Forchheimer's law and a mass conservation equation
\begin{eqnarray}
\mu{\kappa}^{-1}\mathbf{u} + \beta\rho|\mathbf{u}|\mathbf{u} + \nabla{p} = \mathbf{0} &&\textrm{in}\ \ \Omega,\label{eqn_model_1}\\
\nabla\cdot\mathbf{u} = f &&\textrm{in}\ \ \Omega,\label{eqn_model_2}
\end{eqnarray}
where $|\cdot| = (\cdot,\cdot)^{\frac{1}{2}}$ is the discrete $L_2$ norm, $\mu$ is the viscosity, $\kappa$ is the heterogeneous permeability, $\beta$ is the heterogeneous non-Darcy coefficient and $\rho$ is the density of the fluid. The boundary conditions on $\partial\Omega$ is defined as follows
\begin{eqnarray*}
\mathbf{u}\cdot\mathbf{n} = g_N &&\textrm{on}\ \ \partial\Omega_N,\\
p = g_D &&\textrm{on}\ \ \partial\Omega_D,
\end{eqnarray*}
where $\mathbf{n}$ denotes the unit outward normal vector on $\partial\Omega$, $\partial\Omega_N$ and $\partial\Omega_D$ are the Neumann and Dirichlet boundaries, respectively, and $g_N$, $g_D$ are the associated boundary data.

In order to introduce the weak formulation of the model (\ref{eqn_model_1})-(\ref{eqn_model_2}), standard notations and definitions for Sobolev spaces are used, we define the following spaces
\begin{equation*}
{V} = \{\mathbf{v} | \mathbf{v}\in({L}_3(\Omega))^2,\ \nabla\cdot\mathbf{v}\in{L}_2(\Omega)\}\quad\textrm{and}\quad W = L_2(\Omega).
\end{equation*}
The weak formulation of the model (\ref{eqn_model_1})-(\ref{eqn_model_2}) can be written as: find $(\mathbf{u},p)\in{V\times{W}}$ such that
\begin{eqnarray}
(\mu\kappa^{-1}\mathbf{u},\mathbf{v}) + (\beta\rho|\mathbf{u}|\mathbf{u},\mathbf{v}) - (p,\nabla\cdot\mathbf{v}) = (g_D,\mathbf{v}\cdot\mathbf{n})_{\partial\Omega_D} &&\forall\ \mathbf{v}\in V,\label{eqn_weak_1}\\
-(\nabla\cdot\mathbf{u},q) = -(f,q) \hspace{1.07cm}&&\forall\ q\in W.\label{eqn_weak_2}
\end{eqnarray}
The existence and uniqueness of the solution to the above weak formulation (\ref{eqn_weak_1})-(\ref{eqn_weak_2}) have been proved in \cite{pan2012mixed}.
\subsection{Fine-grid approximation by a multipoint flux mixed finite element method}
We adopt an MFMFE method to solve the Darcy-Forchheimer model on the fine grid, which has been studied in \cite{xu2017multipoint}. 
The fine grid $\mathcal{T}_h$ is a conforming shape-regular partition of $\Omega$, composed of convex quadrilaterals or triangles.
\begin{figure}
	\centering\footnotesize
	\begin{overpic}[height=4.5cm,width=10.0cm]{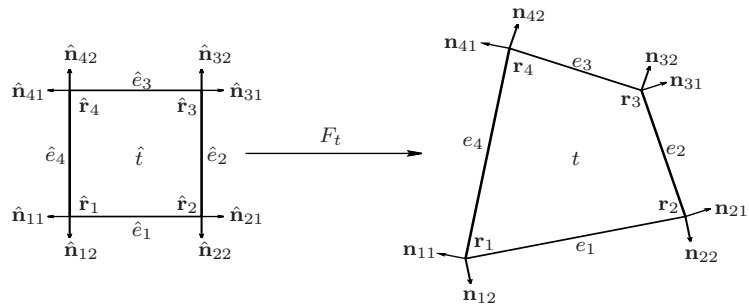}
		\put(10,15){$\hat{\mathbf r}_1$}
		\put(23,15){$\hat{\mathbf r}_2$}
		\put(23,28){$\hat{\mathbf r}_3$}
		\put(10,28){$\hat{\mathbf r}_4$}
		\put(62.5,10){$\mathbf{r}_1$}
		\put(87,15.0){$\mathbf{r}_2$}
		\put(82,29){$\mathbf{r}_3$}
		\put(67.7,33.3){$\mathbf{r}_4$}
		\put(42,24){$F_t$}
		\put(1,13.5){$\hat{\mathbf{n}}_{11}$}
		\put(8,8,8){$\hat{\mathbf{n}}_{12}$}
		\put(30,13.5){$\hat{\mathbf{n}}_{21}$}
		\put(25.8,8.8){$\hat{\mathbf{n}}_{22}$}
		\put(30,30){$\hat{\mathbf{n}}_{31}$}
		\put(25.8,35){$\hat{\mathbf{n}}_{32}$}
		\put(1,30){$\hat{\mathbf{n}}_{41}$}
		\put(8,35){$\hat{\mathbf{n}}_{42}$}	
		\put(53,9){${\mathbf{n}}_{11}$}
		\put(61,3){${\mathbf{n}}_{12}$}
		\put(94.5,14.5){${\mathbf{n}}_{21}$}
		\put(90.5,8.5){${\mathbf{n}}_{22}$}
		\put(88.5,31.5){${\mathbf{n}}_{31}$}
		\put(85.2,34.7){${\mathbf{n}}_{32}$}
		\put(58.5,36.5){${\mathbf{n}}_{41}$}
		\put(67.5,40.5){${\mathbf{n}}_{42}$}
		\put(17,11.5){$\hat{e}_{1}$}
		\put(27,21.5){$\hat{e}_{2}$}
		\put(17,31.5){$\hat{e}_{3}$}
		\put(5.5,21.5){$\hat{e}_{4}$}
		\put(76,9.5){${e}_{1}$}
		\put(88,22.5){${e}_{2}$}
		\put(75.5,34){${e}_{3}$}
		\put(61,23.5){${e}_{4}$}
		\put(17.5,21){$\hat{t}$}
		\put(75.5,21){${t}$}
	\end{overpic}
	\caption{A bijective mapping $F_t$ maps the reference square element $\hat{t}$ to a physical quadrilateral element $t$. $\hat{\mathbf{r}}_i$($\mathbf{r}_i$) and $\hat{e}_i$($e_i$), $i=1,\cdots,4$, are vertices and edges of $\hat{t}$($t$), respectively, with the unit outward normal vectors $\mathbf{\hat{n}}_{ij}$($\mathbf{n}_{ij}$), $i=1,\cdots,4,j=1,2$.}\label{Fig_mapping}
\end{figure}
Let $\hat{t}$ be the reference unit square $[0,1]^2$ with vertices $\hat{\mathbf{r}}_1=(0,0)^T$, $\hat{\mathbf{r}}_2=(1,0)^T$, $\hat{\mathbf{r}}_3=(0,1)^T$ and $\hat{\mathbf{r}}_4=(1,1)^T$, and let $t$ be any physical element in $\mathcal{T}_h$ with vertices $\mathbf{r}_i=(x_i,y_i)^T, i=1,...,4$, then there exists a bijective mapping $F_t:\hat{t}\to t$ as defined in (\ref{eqn_bi_mapping}) and shown in Figure \ref{Fig_mapping}. We denote the Jacobian matrix of $F_t$ by $DF_t$, determinant of $DF_t$ by $J_t=\textrm{det}({DF}_t)|$, inverse mapping of $F_t$ by ${F_t^{-1}}$, Jacobian matrix of ${F_t^{-1}}$ by ${DF_t^{-1}}(x,y)=({DF_t)^{-1}}(\hat{x},\hat{y})$, and determinant of ${DF_t^{-1}}$ by $J_{F_t^{-1}}=1\slash J_t(\hat{x},\hat{y})$, respectively.
\begin{equation}\label{eqn_bi_mapping}
F_t(\hat{x},\hat{y})=\mathbf{r}_1(1-\hat{x})(1-\hat{y})+\mathbf{r}_2\hat{x}(1-\hat{y})+\mathbf{r}_3\hat{x}\hat{y}+\mathbf{r}_4(1-\hat{x})\hat{y}.
\end{equation}
Using the bijective mapping $F_t$, for any scalar function $\hat{w}$ defined in $\hat{t}$, we denote the scalar transformation of $\hat{w}$ in $t$ by $w$, defined as
\begin{equation*}
w\leftrightarrow\hat{w}:w=\hat{w}\circ F_t^{-1}.
\end{equation*}
On the reference unit square $\hat{t}$, the space $\textrm{BDM}_1(\hat{t})$ is defined as
\begin{eqnarray*}
\textrm{BDM}_1(\hat{t}) = P_1(\hat{t})^2 + r\ \textrm{curl}({\hat{x}}^2\hat{y}) + s\ \textrm{curl}(\hat{x}{\hat{y}}^2)
= \left(
\begin{matrix}
\alpha_1\hat{x} + \beta_1\hat{y} + \gamma_1 + r{\hat{x}}^2 + 2s\hat{x}\hat{y}
\\
\alpha_2\hat{x} + \beta_2\hat{y} + \gamma_2 - 2r\hat{x}\hat{y} + s{\hat{y}}^2
\end{matrix}\right).
\end{eqnarray*}
where ${\alpha_i}|_{i=1,2},{\beta_i}|_{i=1,2},{\gamma_i}|_{i=1,2},r,s\in\mathbb{R}$, are arbitrary constants. We take the normal components $\hat{\mathbf{v}}\cdot\hat{\mathbf{n}}_{ij}|_{i=1,\cdots,4,j=1,2}$ at four vertices $\hat{\mathbf{r}}_i|_{i=1,\cdots,4}$ as degrees of freedom for the function $\hat{\mathbf{v}}$ in $\textrm{BDM}_1(\hat{t})$, where $\hat{\mathbf{n}}_{ij}|_{i=1,\cdots,4,j=1,2}$ shown in Figure \ref{Fig_mapping} are unit outward normal vectors of $\hat{t}$ at corners. We choose a set of basis functions such that $\textrm{BDM}_1(\hat{t}) = \textrm{span}\{\hat{\mathbf{v}}_{i,j}|_{i=1,\cdots,4,j=1,2}\}$ and these basis functions satisfy
\begin{align}
\begin{split}
\hat{\mathbf{v}}_{ij}\cdot\hat{\mathbf{n}}_{sl}|_{\hat{e}_{sl}}=\left\{
\begin{array}{ll}
1 & \text{if}\; i=s,j=l, \\
0 & \text{otherwise},
\end{array}\right.
\quad i,s=1,\cdots,4,\ j,l=1,2,
\end{split}
\end{align}
Functions in space $\textrm{BDM}_1(t)$ for any physical element $t\in \mathcal{T}_h$ are defined via the following vector transformation
\begin{equation}\label{eqn_Piola}
\mathbf{v}\leftrightarrow\hat{\mathbf{v}}:\mathbf{v}=\frac{1}{J_t}DF_t\hat{\mathbf{v}}\circ F_t^{-1},\quad \hat{\mathbf{v}}\in\textrm{BDM}_1(\hat{t}),
\end{equation}
which is known as the Piola transformation, preserving the normal components of the velocity vectors on edges, i.e.,
\begin{equation*}
\mathbf{v}\cdot\mathbf{n}_e=\frac{1}{|e|}\hat{\mathbf{v}}\cdot\hat{\mathbf{n}}_{\hat{e}},\quad \forall e\in\partial t.
\end{equation*}
We use the following $\textrm{BDM}_1$ mixed finite element spaces to approximate velocity and pressure for the discretization of the weak formulation (\ref{eqn_weak_1})-(\ref{eqn_weak_2}) on the fine grid $\mathcal{T}_h$,
\begin{eqnarray*}
&&V_h = \{\mathbf{v}\in V : \mathbf{v}|_{t}\in\textrm{BDM}_1(t) \quad\forall t\in\mathcal{T}_h \},\\
&&W_h = \{q\in W : q|_{t}\in\mathbb{P}_0(t) \quad\forall t\in\mathcal{T}_h \},
\end{eqnarray*}
where $\mathbb{P}_0$ denotes the polynomial space of degree zero.\\
We obtain the mixed form of corresponding discrete weak formulation: find a pair $(\mathbf{u}_h,p_h)\in{V_h\times{W}_h}$, such that
\begin{eqnarray}
(\mu\kappa^{-1}\mathbf{u}_h,\mathbf{v}_h) + (\beta\rho|\mathbf{u}_h|\mathbf{u}_h,\mathbf{v}_h) - (p_h,\nabla\cdot\mathbf{v}_h) = (g_D,\mathbf{v}_h\cdot\mathbf{n})_{\partial\Omega_D} &&\forall\ \mathbf{v}_h\in V_h,\label{eqn_discrete_1}\\
-(\nabla\cdot\mathbf{u}_h,q_h) = -(f,q_h) \hspace{1.08cm}&&\forall\ q_h\in W_h.\label{eqn_discrete_2}
\end{eqnarray}
It is well known that the above velocity-pressure system is in a saddle-point structure, which is computational expensive. To avoid tackling the saddle-point algebraic system, we apply the trapezoidal quadrature rule that allows for local velocity elimination and results in a symmetric and positive definite algebraic system for the pressure. 

Suppose that $\mathcal{K}$ is a symmetric tensor defined in $\Omega$, then we apply the trapezoidal quadrature to compute the integration $(\mathcal{K}\mathbf{u},\mathbf{v})$, for any $\mathbf{u},\mathbf{v}\in V_h$. By use of the bilinear mapping (\ref{eqn_bi_mapping}) and the Piola transformation (\ref{eqn_Piola}), the integration on any physical element $t\in\mathcal{T}_h$ is mapped to the reference element $\hat{t}$, that is,
\begin{equation}
(\mathcal{K}\mathbf{u},\mathbf{v})_t = (\frac{1}{J_t}DF_t^T\widehat{\mathcal{K}}DF_t\hat{\mathbf{u}},\hat{\mathbf{v}})_{\hat{t}}=(\widehat{\mathcal{M}}_t\hat{\mathbf{u}},\hat{\mathbf{v}})_{\hat{t}},
\end{equation}
where $\widehat{\mathcal{K}}=\mathcal{K}\circ F_t$, $\widehat{\mathcal{M}}_t=\frac{1}{J_t}DF_t^T\widehat{\mathcal{K}}DF_t$, and $\hat{\mathbf{u}}$, $\hat{\mathbf{v}}\in\textrm{BDM}_1(\hat{t})$ are the inverse functions of $\mathbf{u}$, $\mathbf{v}$ through the Piola transformation (\ref{eqn_Piola}), respectively. By applying the trapezoidal quadrature rule on the reference element $\hat{t}$, we get the quadrature rule on the physical element $t\in\mathcal{T}_h$ as follows
\begin{equation}
(\mathcal{K}\mathbf{u},\mathbf{v})_{Q,t}=(\widehat{\mathcal{M}}_t\hat{\mathbf{u}},\hat{\mathbf{v}})_{\hat{Q},\hat{t}}=
\frac{|\hat{t}|}{4}\sum\limits_{i=1}^{4}{\widehat{{\mathcal{M}}}_t(\hat{\mathbf{r}}_i)\hat{\mathbf{u}}(\hat{\mathbf{r}}_i)\cdot\hat{\mathbf{v}}(\hat{\mathbf{r}}_i)},
\end{equation}
and the global quadrature rule for the integration $(\mathcal{K}\mathbf{u},\mathbf{v})$ in $\Omega$ is defined as
\begin{equation}\label{eq315}
(\mathcal{K}\mathbf{u},\mathbf{v})_Q = \sum\limits_{t\in\mathcal{T}_h}{(\mathcal{K}\mathbf{u},\mathbf{v})_{Q,t}}.
\end{equation}
The above quadrature rule $(\cdot,\cdot)_Q$ only couples the two basis functions of velocity that associated with the same vertex, for example, on the reference element $\hat{t}$
\begin{eqnarray}
&&(\widehat{\mathcal{M}}_t\hat{\mathbf{v}}_{11},\hat{\mathbf{v}}_{11})_{\hat{Q},\hat{t}} = \frac{\widehat{\mathcal{M}}_{t,11}(\hat{\mathbf{r}}_1)}{4},\quad
(\widehat{\mathcal{M}}_t\hat{\mathbf{v}}_{11},\hat{\mathbf{v}}_{12})_{\hat{Q},\hat{t}} = \frac{\widehat{\mathcal{M}}_{t,12}(\hat{\mathbf{r}}_1)}{4},\label{eqn_diag_1}\\
&&(\widehat{\mathcal{M}}_t\hat{\mathbf{v}}_{11},\hat{\mathbf{v}}_{ij})_{\hat{Q},\hat{t}} = 0\quad\forall i\neq1,j=1,2.\label{eqn_diag_2}
\end{eqnarray}
where $\widehat{\mathcal{M}}_{t,ij}$ is the $ij$-th component of $\widehat{\mathcal{M}}_t$, $i=1,\cdots,4,j=1,2$, respectively.

We obtain the corresponding discrete weak formulation using the MFMFE method: find a pair $(\mathbf{u}_h,p_h)\in{V_h\times{W}_h}$, such that
\begin{eqnarray}
(\mu\kappa^{-1}\mathbf{u}_h,\mathbf{v}_h)_Q + (\beta\rho|\mathbf{u}_h|\mathbf{u}_h,\mathbf{v}_h)_Q - (p_h,\nabla\cdot\mathbf{v}_h) = (g_D,\mathbf{v}_h\cdot\mathbf{n})_{\partial\Omega_D} &&\forall\ \mathbf{v}_h\in V_h,\label{eqn_mfmfe_1}\\
-(\nabla\cdot\mathbf{u}_h,q_h) = -(f,q_h) \hspace{1.08cm}&&\forall\ q_h\in W_h.\label{eqn_mfmfe_2}
\end{eqnarray}
The Picard iterative algorithm and Newton iterative algorithm for solving the above nonlinear discrete system (\ref{eqn_mfmfe_1})-(\ref{eqn_mfmfe_2}) on the fine grid $\mathcal{T}_h$ are written as follows\\
\textit{Picard iterative algorithm:} given arbitrary $\mathbf{u}^0_h\in V_h$, find a pair $(\mathbf{u}^{n+1}_h,p^{n+1}_h)\in V_h\times W_h$, such that
\begin{eqnarray}
(\mu\kappa^{-1}\mathbf{u}^{n+1}_h,\mathbf{v}_h)_Q + (\beta\rho|\mathbf{u}^{n}_h|\mathbf{u}^{n+1}_h,\mathbf{v}_h)_Q - (p^{n+1}_h,\nabla\cdot\mathbf{v}_h)
= (g_D,\mathbf{v}_h\cdot\mathbf{n})_{\partial\Omega_D} &&\forall\ \mathbf{v}_h\in V_h,\label{eqn_Picard_1}\\
-(\nabla\cdot\mathbf{u}^{n+1}_h,q_h) = -(f,q_h) \hspace{1.06cm}&&\forall\ q_h\in W_h.\label{eqn_Picard_2}
\end{eqnarray}
\textit{Newton iterative algorithm:} given arbitrary $\mathbf{u}^0_h\in V_h$, find a pair $(\mathbf{u}^{n+1}_h,p^{n+1}_h)\in V_h\times W_h$, such that
\begin{eqnarray}
(\mu\kappa^{-1}\mathbf{u}^{n+1}_h,\mathbf{v}_h)_Q + (\beta\rho|\mathbf{u}^{n}_h|\mathbf{u}^{n+1}_h,\mathbf{v}_h)_Q + (\beta\rho\mathcal{P}^n_h\mathbf{u}^{n+1}_h,\mathbf{v}_h)_Q - (p^{n+1}_h,\nabla\cdot\mathbf{v}_h) \hspace{1.0cm}&&\hspace{1.5cm}\nonumber\\
= (g_D,\mathbf{v}_h\cdot\mathbf{n})_{\partial\Omega_D} + (\beta\rho\mathcal{P}^n_h\mathbf{u}^{n}_h,\mathbf{v}_h)_Q&&\forall\ \mathbf{v}_h\in V_h,\label{eqn_Newton_1}\\
-(\nabla\cdot\mathbf{u}^{n+1}_h,q_h) = -(f,q_h) \hspace{3.85cm}&&\forall\ q_h\in W_h,\label{eqn_Newton_2}
\end{eqnarray}
where $\mathcal{P}^n_h = \frac{\mathbf{u}^n_h\otimes\mathbf{u}^n_h}{|\mathbf{u}^{n}_h|}$, $\mathbf{u}^n_h\otimes\mathbf{u}^n_h = \mathbf{u}^n_h{\mathbf{u}^n_h}^T$.

Suppose the dimensions of $V_h$ and $W_h$ are $m_1$ and $m_2$, respectively, then the above iterative algorithms can be written into matrix forms as : given arbitrary vector ${U}^0_h\in\mathbb{R}^{m_1}$, find a pair $({U}^{n+1}_h,{P}^{n+1}_h)\in\mathbb{R}^{m_1}\times\mathbb{R}^{m_2}$, such that
\begin{equation}\label{eqn_matrix_fine}
\left(
\begin{matrix}
A^n_h\ \ &B_h\\ \\B^T_h\ &0
\end{matrix}\right)
\left(
\begin{matrix}
U^{n+1}_h\\ \\P^{n+1}_h
\end{matrix}\right)=
\left(
\begin{matrix}
G^n_h\\ \\F_h
\end{matrix}\right),
\end{equation}
where the matrix $A^n_h$ is associated with terms
\begin{equation*}
(\mu\kappa^{-1}\mathbf{u}_h^{n+1},\mathbf{v}_h)_Q + (\beta\rho|\mathbf{u}^{n}_h|\mathbf{u}_h^{n+1},\mathbf{v}_h)_Q,
\end{equation*}
in the Picard iterative algorithm (\ref{eqn_Picard_1})-(\ref{eqn_Picard_2}) or terms
\begin{equation*}
(\mu\kappa^{-1}\mathbf{u}_h^{n+1},\mathbf{v}_h)_Q + (\beta\rho|\mathbf{u}^{n}_h|\mathbf{u}_h^{n+1},\mathbf{v}_h)_Q + (\beta\rho\mathcal{P}^n_h\mathbf{u}_h^{n+1},\mathbf{v}_h)_Q,
\end{equation*}
in the Newton iterative algorithm (\ref{eqn_Newton_1})-(\ref{eqn_Newton_2}). $B_h$, $G^n_h$ and $F_h$ are associated with terms $(p_h^{n+1},\nabla\cdot\mathbf{v}_h)$, $(g_D,\mathbf{v}_h\cdot\mathbf{n})_{\partial\Omega_D} + (\beta\rho\mathcal{P}^n_h\mathbf{u}^{n}_h,\mathbf{v}_h)_Q$ and $-(f,q_h)$, respectively. From property (\ref{eqn_diag_1})-(\ref{eqn_diag_2}) and by means of the $\textrm{BDM}_1$ mixed finite element spaces and the trapezoidal quadrature rule $(\cdot,\cdot)_Q$, we know that $A^n_h$ is a block diagonal, symmetric and positive definite matrix. Thus $A^n_h$ can  be inverted easily and we can solve the system (\ref{eqn_matrix_fine}) for each iteration in the following way
\begin{equation}\label{matrixform}
-B^T_h({A^n_h})^{-1}B_h P^{n+1}_h = F_h-B^T_h({A^n_h})^{-1}G^n_h,
\end{equation}
namely, we only need to solve symmetric and positive definite systems for the pressure.
\section{Coarse-grid approximation}
\begin{figure}
	\centering\footnotesize
	\begin{overpic}[height=6.1cm,width=13cm]{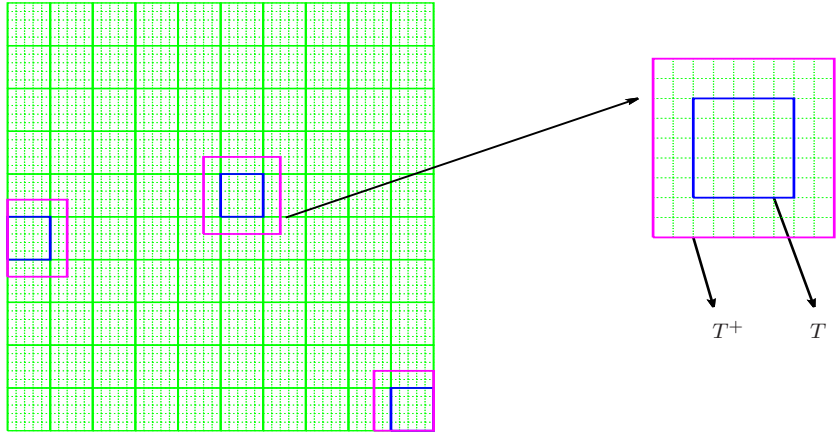}
		\put(85,11){$T$}
		\put(75,11){$T^+$}
	\end{overpic}
	\caption{The illustration of a multiscale mesh in the left and an oversampling block $T^+_i$ associated with a coarse element $T_i$ in the right.}
	\label{Fig_mesh_Oversampling}
\end{figure}

In this section, we make an approximation of the Darcy-Forchheimer model on the coarse grid by illustrating a systematic way to construct the multiscale space for pressure, which follows the GMsFEM framework. The coarse grid is denoted by $\mathcal{T}_H$. Each coarse element $T$ in the coarse grid $\mathcal{T}_H$ is a connected collection of elements in the fine grid $\mathcal{T}_h$; that is, coarse element number $i$ is formed by $T_i = \cup^{N_i}_{k=1}t_k$, where $N_i$ is the number of fine-grid elements contained in $T_i$. In the simplest case, the coarse grid turns into a uniform partition of the fine Cartesian grid so that each coarse element $T$ becomes a rectangle. See Figure \ref{Fig_mesh_Oversampling} for an example of a multiscale mesh and a coarse element $T$. We use $N_T$ to denote the total number of coarse elements included in $\mathcal{T}_H$. The coarse grid approximation is separated into two stages: offline computation and online computation.
\subsection{Offline computation}
In this subsection, we first construct the local snapshot space on each coarse element by solving a series of local problems with different boundary conditions. The snapshot space provides a solution space in each coarse element locally. Then we perform a spectral decomposition in each local snapshot space to derive the local offline space with a smaller dimension. All local offline spaces form into the offline space which is used to approximate the original problem on the coarse grid and obtain the offline solution. Both the offline space and the offline solution are of great importance and can be efficiently employed in the online computation stage to calculate the online basis functions for the enrichment of the multiscale space.
\subsubsection{Snapshot space}
Let $T_i \in \mathcal{T}_H$ be a coarse element in $\Omega$. Basis functions of the local snapshot space $W^{i}_{\textrm{snap}}$ corresponding to $T_i$ is derived by solving the following problems numerically: find $(\mathbf{\psi}_j^{i}, \phi_j^{i})$ such that
\begin{align}\label{eqn_local_problem}
\begin{split}
\mu\kappa^{-1}\mathbf{\psi}_j^{i} + \nabla \phi_j^{i} &= 0 \qquad \mathrm{in} \; T_i, \\
\nabla\cdot\mathbf{\psi}_j^{i} & = 0 \qquad \mathrm{in} \; T_i.
\end{split}
\end{align}
The boundary of coarse element can be written as a collection of fine-grid edges, given as $\partial T_i = \bigcup_{j = 1}^{J_i} e_j $, where $J_i$ is the total number of fine-grid edges on the boundary of coarse element $T_i$.
Let $\delta^{i}_j $ be a piecewise constant function defined on $\partial T_i$ with respect to the fine-grid edge such that it has value $1$ on
$e_j$ and value $0$ on the other fine-grid edges, defined as
\begin{align}\label{eqn_cellConstraint}
\begin{split}
\delta^{i}_j = \left\{
\begin{array}{ll}
1 & \text{in}\; e_j, \\
0 & \text{on\, other \, fine-grid \,edges \,on\, } \partial T_i,
\end{array}\right.
\qquad j = 1, 2, \cdots, J_i.
\end{split}
\end{align}
The boundary conditions on the boundary of coarse element $T_i$ for the local problem (\ref{eqn_local_problem}) are taken as
\begin{align}\label{eqn_cell_BC}
\phi^{i}_j = \delta^{i}_j \quad \text{on} \; \partial T_i.
\end{align}
Therefore, we obtain the local snapshot space associated with the coarse element $T_i$ as
\begin{equation*}
W^{i}_{\textrm{snap}} = \textrm{span}\{\phi^{i}_1,\phi^{i}_2,\cdots,\phi^{i}_{J_i}\}.
\end{equation*}
Additionally, we define the following local space that will be used in the local spectral decomposition to derive the local offline space
\begin{equation*}
V^{i}_{\textrm{snap}} = \textrm{span}\{\psi^{i}_1,\psi^{i}_2,\cdots,\psi^{i}_{J_i}\}.
\end{equation*}
\textit{Remark 1 :} Oversampling technique \cite{efendiev2014generalized} can be employed to get more effective snapshot spaces. Let $T^+$ be a coarse block defined by adding some fine-grid layers around $T$ such that $T\subset T^+$, as shown in the right of Figure \ref{Fig_mesh_Oversampling}. The snapshot basis functions with respect to the coarse block $T^+$ are derived form solving the problem (\ref{eqn_local_problem}) in the oversampling region $T^+$ with the similar boundary conditions with (\ref{eqn_cell_BC}) defined on $\partial T^+$.
\subsubsection{Offline space}
We construct the local offline space by performing a dimension reduction in the local snapshot space. To this end, we solve the following spectral decomposition problem to get the dominate modes in each local snapshot space $W^{i}_{\textrm{snap}}$: find a real number $\lambda_k\ge 0$ and a vector $\Phi_k$ such that
\begin{equation}\label{eqn_spectral}
A^{i}_{\textrm{off}}\Phi_k = \lambda_k S^{i}_{\textrm{off}}\Phi_k,\quad A^{i}_{\textrm{off}} = {R^{v,i}_{\textrm{off}}}^{T}A^{i}R^{v,i}_{\textrm{off}},\quad S^{i}_{\textrm{off}} = {R^{w,i}_{\textrm{off}}}^{T}S^{i}R^{w,i}_{\textrm{off}},
\end{equation}
where, $R^{v,i}_{\textrm{off}}$ and $R^{w,i}_{\textrm{off}}$ denote the coefficients matrix of snapshot basis functions in the expansion of fine-grid basis functions
\begin{equation}
R^{v,i}_{\textrm{off}} = [\psi^{i}_1,\psi^{i}_2,\cdots,\psi^{i}_{J_i}],\quad R^{w,i}_{\textrm{off}} = [\phi^{i}_1,\phi^{i}_2,\cdots,\phi^{i}_{J_i}],
\end{equation}
$A^{i}$ and $S^{i}$ are fine-grid matrices associated with the following bilinear forms,
\begin{equation}
A^{i} = [a^{i}_{rl}] = (\mu\kappa^{-1}\psi_r,\psi_l)_{Q,T_i},\quad S^{i} = [s^{i}_{rl}] = (\phi_r,\phi_l)_{T_i}.
\end{equation}
We arrange the eigenvalues of (\ref{eqn_spectral}) in increasing order,
\begin{align}
\lambda_1^{i}< \lambda_2^{i}< \cdots < \lambda_{J_i}^{i},
\end{align}
and choose the first $M^{i}_{\textrm{off}}$ eigenvalues $\lambda_k^{i}$ and the corresponding eigenvectors $\Phi_k = (\Phi_{kj})^{J_i}_{j=1}$ to form the local offline space with respect to $T_i$, where $(\Phi_{kj})$ is the $j$-th component of the vector $\Phi_k$ for $k=1,\cdots,M^{i}_{\textrm{off}}$. Hence, we define the local offline basis functions to be
\begin{equation}
\phi^{i,\textrm{off}}_k = \sum^{J_i}_{j=1}\Phi_{kj}\phi^{i}_j, \qquad k=1,\dots,M^{i}_{\textrm{off}}.
\end{equation}
We define the local offline space with respect to $T_i$ as
\begin{equation}
W^{i}_{\textrm{off}} = \textrm{span}\{\phi^{i,\textrm{off}}_1,\phi^{i,\textrm{off}}_2,\cdots,\phi^{i,\textrm{off}}_{M^{i}_{\textrm{off}}}\}.
\end{equation}
Combine all these local offline spaces $W^{i}_{\textrm{off}}$, $i=1,\cdots,N_T$, we get the global offline space $W_{\textrm{off}}$ for the pressure, and by using of the single-index notation, it can be written as
\begin{equation}
W_{\textrm{off}} = \textrm{span}\{\psi^{\textrm{off}}_{k} : 1\le k\le M_{\textrm{off}}\}.
\end{equation}
where $M_{\textrm{off}} = \sum^{N_T}_{i=1}M^i_{\textrm{off}}$ is the total dimension of the global offline space for approximating the pressure. We use the matrix $R_{\textrm{off}}$ to denote the coefficients vector of each offline basis function in the expansion of fine-grid basis functions.\\
Having gotten the offline space $W_{\textrm{off}}$ for pressure, we define the following offline spaces for velocity and pressure as
\begin{equation}
V_H=V_h \quad\textrm{and} \quad W_H=W_{\textrm{off}},\label{eqn:coarse_grid_space}
\end{equation}
and get the mixed GMsFEM system with the following iterative algorithms.\\
\textit{Picard iterative algorithm using offline space:} given arbitrary $\mathbf{u}^0_H\in V_H$, find a pair $(\mathbf{u}^{n+1}_H,p^{n+1}_H)\in V_H\times W_H$, such that
\begin{eqnarray}
(\mu\kappa^{-1}\mathbf{u}^{n+1}_H,\mathbf{v}_H)_Q + (\beta\rho|\mathbf{u}^{n}_H|\mathbf{u}^{n+1}_H,\mathbf{v}_H)_Q - (p^{n+1}_H,\nabla\cdot\mathbf{v}_H)
= (g_D,\mathbf{v}_H\cdot\mathbf{n})_{\partial\Omega_D} &&\forall\ \mathbf{v}_H\in V_H,\label{eqn_Picard_1_H}\\
-(\nabla\cdot\mathbf{u}^{n+1}_H,q_H) = -(f,q_H) \hspace{1.06cm}&&\forall\ q_H\in W_H.\label{eqn_Picard_2_H}
\end{eqnarray}
\textit{Newton iterative algorithm using offline space:} given arbitrary $\mathbf{u}^0_H\in V_H$, find a pair $(\mathbf{u}^{n+1}_H,p^{n+1}_H)\in V_H\times W_H$, such that
\begin{eqnarray}
(\mu\kappa^{-1}\mathbf{u}^{n+1}_H,\mathbf{v}_H)_Q + (\beta\rho|\mathbf{u}^{n}_H|\mathbf{u}^{n+1}_H,\mathbf{v}_H)_Q + (\beta\rho\mathcal{P}^n_H\mathbf{u}^{n+1}_H,\mathbf{v}_H)_Q - (p^{n+1}_H,\nabla\cdot\mathbf{v}_H) \hspace{0.5cm}&&\hspace{1.5cm}\nonumber\\
= (g_D,\mathbf{v}_H\cdot\mathbf{n})_{\partial\Omega_D} + (\beta\rho\mathcal{P}^n_H\mathbf{u}^{n}_H,\mathbf{v}_H)&&\forall\ \mathbf{v}_H\in V_H,\label{eqn_Newton_1_H}\\
-(\nabla\cdot\mathbf{u}^{n+1}_H,q_H) = -(f,q_H) \hspace{3.85cm}&&\forall\ q_H\in W_H,\label{eqn_Newton_2_H}
\end{eqnarray}
where $\mathcal{P}^n_H = \frac{\mathbf{u}^n_H\otimes\mathbf{u}^n_H}{|\mathbf{u}^{n}_H|}$, $\mathbf{u}^n_H\otimes\mathbf{u}^n_H = \mathbf{u}^n_H{\mathbf{u}^n_H}^{T}$.\\
The matrix form of the above coarse-grid iterative algorithms can be written as : given arbitrary vector ${U}^0_H\in\mathbb{R}^{m_1}$, find a pair $({U}^{n+1}_H,{P}^{n+1}_H)\in\mathbb{R}^{m_1}\times\mathbb{R}^{M_{\textrm{off}}}$, such that
\begin{equation}\label{eqn:systemGMs}
\left(
\begin{matrix}
A^n_H\ \ &B_hR_{\textrm{off}}\\ \\R_{\textrm{off}}^TB^T_h\ &0
\end{matrix}\right)
\left(
\begin{matrix}
U^{n+1}_H\\ \\P^{n+1}_H
\end{matrix}\right)=
\left(
\begin{matrix}
G^n_H\\ \\R_{\textrm{off}}^TF_h
\end{matrix}\right).
\end{equation}
Since the matrix $A^n_H$ is defined in the same way as the matrix $A^n_h$ in ({\ref{eqn_matrix_fine}) whcih is block diagonal, symmetric and positive definite and can be inverted easily, we solve the system (\ref{eqn:systemGMs}) for each iteration in the following way to get the offline solution
\begin{equation}\label{matrixformGMs}
-R_{\textrm{off}}^TB^T_h({A^n_H})^{-1}B_hR_{\textrm{off}} P^{n+1}_H = R_{\textrm{off}}^TF_h-R_{\textrm{off}}^TB^T_H({A^n_H})^{-1}G^n_H.
\end{equation}
\subsection{Online computation}
In this subsection, we compute online basis functions adaptively in selected regions based on the offline space and residual indicators to enrich the multiscale space. Online basis functions contain useful global information that offline basis functions can't capture. Let $(\mathbf{u}_{\textrm{off}}, p_{\textrm{off}})$ be the solution at the end of the Newton iterations (\ref{eqn_Newton_1_H})-(\ref{eqn_Newton_2_H}), namely the offline solution. We denote $W^0_{\textrm{ms}} = W_{\textrm{off}}$ as the initial online multiscale space and $(\mathbf{u}^0_{\textrm{ms}},p^0_{\textrm{ms}}) = (\mathbf{u}_{\textrm{off}},p_{\textrm{off}})$ as the initial solution for the online multiscale space enrichment procedure. We make use of the index $m\ge 1$ to represent the enrichment level of the online multiscale space. In the enrichment level $m$, we use $W^m_{\textrm{ms}}$ and $(\mathbf{u}^m_{\textrm{ms}},p^m_{\textrm{ms}})$ to denote the corresponding online multiscale space and multiscale solution, respectively. The computataion of online basis functions is illustrated in the following.

Suppose the $m$-th level online multiscale space $W^m_{\textrm{ms}}$ and the corresponding multiscale solution $(\mathbf{u}^m_{\textrm{ms}},p^m_{\textrm{ms}})$ are already known, and we need to construct online basis function $\phi$ on the coarse element $T$ to enrich the multiscale space such that $W^{m+1}_{\textrm{ms}} = W^m_{\textrm{ms}} + \textrm{span}\{\phi\}$. Let $T^+$ be the coarse block inclusive of $T$, $T\subset T^+$, defined by adding one fine-grid layer around $T$. We solve the following problem: find $(\psi^+,\phi^+)\in V_h(T^+)\times W_h(T^+)$, such that
\begin{eqnarray}
	(\mu\kappa^{-1}\psi^+,v)_Q + (\beta\rho|\mathbf{u}^m_{\textrm{ms}}|\psi^+,v)_Q - (\phi^+,\nabla\cdot v)
	= 0 &&\forall\ v\in V_h(T^+),\label{eqn_online_1}\\
	(\nabla\cdot \psi^+,q) = (f-\nabla\cdot\mathbf{u}^m_{\textrm{ms}},q) &&\forall\ q\in W_h(T^+),\label{eqn_online_2}
\end{eqnarray}
and satisfies the homogeneous Neumann boundary condition $\psi^+\cdot n = 0$ on $\partial T^+$, where $\phi^+$ is uniquely decided with the condition that $\phi^+=0$ on the boundary element of $T^+$. Restrict $\phi^+$ on $T$, we get the online basis function on the coarse element $T$, i.e., $\phi = \phi^+|_T$.

With the online multiscale space $W^{m+1}_{\textrm{ms}}$ already known, we get the corresponding multiscale solution by solving the following problem: find $(\mathbf{u}^{m+1}_{\textrm{ms}},p^{m+1}_{\textrm{ms}})\in V_h\times W^{m+1}_{\textrm{ms}}$, such that
\begin{eqnarray}
	(\mu\kappa^{-1}\mathbf{u}^{m+1}_{\textrm{ms}},\mathbf{v}_h)_Q + (\beta\rho|\mathbf{u}^{m}_{\textrm{ms}}|\mathbf{u}^{m+1}_{\textrm{ms}},\mathbf{v}_h)_Q - (p^{n+1}_h,\nabla\cdot\mathbf{v}_h)
	= (g_D,\mathbf{v}_h\cdot\mathbf{n})_{\partial\Omega_D} &&\forall\ \mathbf{v}_h\in V_h,\label{eqn_OnlinePicard_1}\\
	-(\nabla\cdot\mathbf{u}^{m+1}_{\textrm{ms}},q_h) = -(f,q_h) \hspace{1.06cm}&&\forall\ q_h\in W^{m+1}_{ms}.\label{eqn_OnlinePicard_2}
\end{eqnarray}
\textit{Remark 2 :} In the matrix form of the system (\ref{eqn_OnlinePicard_1})-(\ref{eqn_OnlinePicard_2}), we need to update the mass matrix for velocity $A^m_H$ associated with terms $(\mu\kappa^{-1}\mathbf{u}^{m+1}_{\textrm{ms}},\mathbf{v}_h)_Q + (\beta\rho|\mathbf{u}^{m}_{\textrm{ms}}|\mathbf{u}^{m+1}_{\textrm{ms}},\mathbf{v}_h)_Q$ in each level of online multiscale space enrichment. So we also test the online multiscale space enrichment in the case that we use the term $(\beta\rho|\mathbf{u}_f|\cdot,\cdot)_Q$ to replace the term $(\beta\rho|\mathbf{u}_{\textrm{off}}|\cdot,\cdot)_Q$ in (\ref{eqn_online_1_remark2}) and (\ref{eqn_OnlinePicard_1_remark2}), respectively, i.e., the online basis function and the multiscale solution are solved in the following way
\begin{eqnarray}
(\mu\kappa^{-1}\psi^+,v)_Q + (\beta\rho|\mathbf{u}_{\textrm{off}}|\psi^+,v)_Q - (\phi^+,\nabla\cdot v)
= 0 &&\forall\ v\in V_h(T^+),\label{eqn_online_1_remark2}\\
(\nabla\cdot \psi    ^+,q) = (f-\nabla\cdot\mathbf{u}^m_{\textrm{ms}},q) &&\forall\ q\in W_h(T^+),\label{eqn_online_2_remark2}
\end{eqnarray}
and
\begin{eqnarray}
(\mu\kappa^{-1}\mathbf{u}^{m+1}_{\textrm{ms}},\mathbf{v}_h)_Q + (\beta\rho|\mathbf{u}_{\textrm{off}}|\mathbf{u}^{m+1}_{\textrm{ms}},\mathbf{v}_h)_Q - (p^{n+1}_h,\nabla\cdot\mathbf{v}_h)
= (g_D,\mathbf{v}_h\cdot\mathbf{n})_{\partial\Omega_D} &&\forall\ \mathbf{v}_h\in V_h,\label{eqn_OnlinePicard_1_remark2}\\
-(\nabla\cdot\mathbf{u}^{m+1}_{\textrm{ms}},q_h) = -(f,q_h) \hspace{1.06cm}&&\forall\ q_h\in W^{m+1}_{ms}.\label{eqn_OnlinePicard_2_remark2}
\end{eqnarray}
\section{Numerical tests}
In this section, we present some numerical examples to demonstrate the performance of the proposed multiscale method for Darcy-Forchheimer model in heterogeneous porous media. In the following examples, we set $\mu=1$ and $\rho=1$. The Darcy-Forchheimer coefficient are taken to be $\beta = \beta_0{\kappa}^{-1}$ \cite{spiridonov2019mixed,li2001literature,muljadi2016impact}, where the parameter $\beta_0$ control the influence of the nonlinear term and we will test cases with $\beta_0=1$, $\beta_0=10$, $\beta_0=100$, $\beta_0=1000$ and $\beta_0=10000$, respectively. Denote the fine-grid solution by $(p_f,\mathbf{u}_f)$, suppose the multiscale solution is denoted by $(p_{\textrm{ms}},\mathbf{u}_{\textrm{ms}})$, then the relative $L_2$ errors for pressure and velocity are denoted as follows
\begin{equation}
\textrm{Erp}(p_{\textrm{ms}}):=\|p_{\textrm{ms}} - p_f\|\slash\|p_f\| \quad\textrm{and}\quad \textrm{Eru}(\mathbf{u}_{\textrm{ms}}):=\|\mathbf{u}_{\textrm{ms}} - \mathbf{u}_f\|\slash\|\mathbf{u}_f\|.\nonumber
\end{equation}

\subsection{Offline solution}
We first investigate the performance of the offline solution. The offline solution is derived by solving the Newton iterative algorithm (\ref{eqn_Newton_1_H})-(\ref{eqn_Newton_2_H}), denoted by $(\mathbf{u}_{\textrm{off}},p_{\textrm{off}})$. When we get the offline solution for velocity $\mathbf{u}_{\textrm{off}}$, we can update the local snapshot spaces with respect to the coarse elements where the residuals are large. The updated local snapshot spaces may capture the fine-scale information more accurately than the original local snapshot spaces. The residual $R_i$ for the coarse element $T_i$ is computed as
\begin{equation}\label{eqn_residual_offline}
R_i = \int_{T_i}{|f-\nabla\cdot\mathbf{u}_{\textrm{off}}|^2}\textrm{d}x.
\end{equation}
We arrange the above $N_T$ residuals in decreasing order, that is $R_1 \ge R_2 \ge \cdots \ge R_{N_T}$, and we choose the coarse elements where the local snapshot spaces need to be updated by choosing the smallest integer $N_{\textrm{update}}$, such that
\begin{equation}\label{eqn_criterion_offline}
\sum^{N_{\textrm{update}}}_{i=1}R_i \ge \theta\sum^{N_T}_{i=1}R_i.
\end{equation}
where $1>\theta>0$ is a real number to be chosen. In the following examples of this subsection, we take $\theta=3/4$. The local snapshot spaces in regard to the selected coarse elements are updated by solving the following local problems : find $(\mathbf{\psi}_j^{(i)}, \phi_j^{(i)})$ such that
\begin{eqnarray}
(\mu\kappa^{-1}+\beta\rho|\mathbf{u}_{\textrm{off}}|)\mathbf{\psi}_j^{(i)} + \nabla \phi_j^{(i)} = 0 &\quad\mathrm{in} \; T_i,\label{eqn_cell_update_1}\\
\nabla\cdot\mathbf{\psi}_j^{(i)} = 0 &\quad\mathrm{in} \; T_i,\label{eqn_cell_update_2}
\end{eqnarray}
with the same boundary conditions as (\ref{eqn_cell_BC}), for $i=1,\cdots,N_{\textrm{update}}$. Then we implement the same spectral decomposition as (\ref{eqn_spectral}) in the updated snapshot spaces to update the offline space. We denote the partially updated offline space according to the criterion (\ref{eqn_criterion_offline}) by $\widehat{W}_{\textrm{off}}$, the corresponding offline solution by $(\hat{\mathbf{u}}_{\textrm{off}},\hat{p}_{\textrm{off}})$, and the associated $L_2$ relative errors by $\textrm{Errp}(\hat{p}_{\textrm{off}})$, $\textrm{Erru}(\hat{\mathbf{u}}_{\textrm{off}})$, respectively. For comparison, we also denote the totally updated offline space for all coarse elements ($N_{\textrm{update}}=N_T$) by $\widetilde{W}_{\textrm{off}}$, the corresponding offline solution by $(\tilde{\mathbf{u}}_{\textrm{off}},\tilde{p}_{\textrm{off}})$, and the associated $L_2$ relative errors by $\textrm{Erp}(\tilde{p}_{\textrm{off}})$, $\textrm{Erru}(\tilde{\mathbf{u}}_{\textrm{off}})$, respectively.\\

{\bf Example 1 :}
\begin{figure}[h!]
	\centering
	\includegraphics[width=0.45\textwidth]{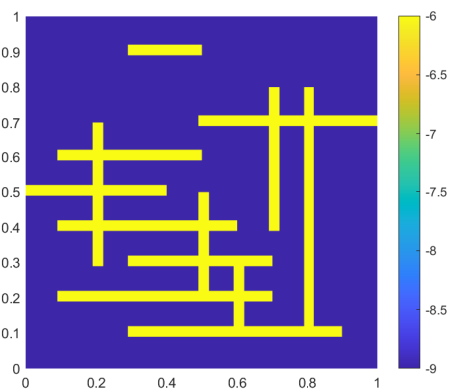}
	\caption{The distribution of permeability field $\kappa$ in logarithmic scale for example 1.}
	\label{fig_perm_1}
\end{figure}
The computational domain is set to be $\Omega=[0, 1]^2$, the fine grid is a $100\times100$ uniform mesh, and the coarse grid is a $10\times10$ uniform mesh. The permeability field $\kappa$ is shown in Figure \ref{fig_perm_1}.

Table \ref{tab_iteration_number_1} shows respectively the iteration number of Picard iterative algorithm (\ref{eqn_Picard_1_H})-(\ref{eqn_Picard_2_H}) and Newton iterative algorithm (\ref{eqn_Newton_1_H})-(\ref{eqn_Newton_2_H}) with the parameter $\beta_0$ taken different values, it is obvious that the number of Newton iterations is much less than Picard iterations, which is especially noticeable when the parameter $\beta_0$ takes large values.
\begin{table}[h!]
	\caption{(Example 1) The number of iterations.}\label{tab_iteration_number_1}
	\centering
	\begin{tabular*}
		{0.5\textwidth}{@{\extracolsep{\fill}}|c|c|c|}
		\hline
		$\beta_0$  & $\#$Newton iterations  & $\#$Picard iterations\\
		\hline
		1                &7                       &48   \\
		10               &9                       &153  \\
		100              &11                      &491  \\
		1000             &12                      &1376  \\
		10000            &14                      &3057  \\
		\hline
	\end{tabular*}
\end{table}

Tables \ref{tab_err_ex1_4}, \ref{tab_err_ex1_6}, and \ref{tab_err_ex1_8} show the $L_2$ relative errors $\textrm{Erp}(p_{\textrm{off}})$, $\textrm{Eru}(\mathbf{u}_{\textrm{off}})$, $\textrm{Erp}(\hat{p}_{\textrm{off}})$, $\textrm{Eru}(\hat{\mathbf{u}}_{\textrm{off}})$, $\textrm{Erp}(\tilde{p}_{\textrm{off}})$ and $\textrm{Eru}(\tilde{\mathbf{u}}_{\textrm{off}})$ with parameter $\beta_0$ taken different values with respect to $4$, $6$ and $8$ offline basis functions per coarse element, respectively, where `Dof per T' denotes the number of offline basis functions per coarse element. We find that the offline basis functions provides good approximations of the problem on the coarse grid and the accuracy of the offline solution is improved by applying more offline basis functions per coarse element. The updated offline basis functions are capable of reducing the relative errors and it becomes evident when the parameter $\beta_0$ becomes large. Moreover, comparisons between the results of partially and totally updated offline space demonstrate that the criterion (\ref{eqn_criterion_offline}) is valid to save the computation for the offline space update, that is, we can realize the comparable error reduction by only updating a part of local snapshot and offline spaces.

Figure \ref{fig_pressurevelocity_1} plots the fine-grid solution $(p_f,\mathbf{u}_f)$, offline solution $(p_\textrm{off},\mathbf{u}_{\textrm{off}})$ and the updated offline solution $(\hat{p}_{\textrm{off}},\hat{\mathbf{u}}_{\textrm{off}})$, respectively, with $4$ offline basis functions per coarse element and the parameter $\beta_0=100$. We can see that the offline solutions are able to achieve good accuracies though a small number of offline basis functions are used.
\begin{table}[h!]
	\caption{(Example 1) Relative errors $\textrm{Erp}(p_{\textrm{off}})$, $\textrm{Eru}(\mathbf{u}_{\textrm{off}})$, $\textrm{Erp}(\hat{p}_{\textrm{off}})$, $\textrm{Eru}(\hat{\mathbf{u}}_{\textrm{off}})$, $\textrm{Erp}(\tilde{p}_{\textrm{off}})$ and $\textrm{Eru}(\tilde{\mathbf{u}}_{\textrm{off}})$ with $4$ offline basis functions per coarse element, $\theta=3/4$.}\label{tab_err_ex1_4}
	\centering
	\begin{tabular*}
		{0.92\textwidth}{@{\extracolsep{\fill}}|c|c|c|c|c|c|c|c|c|}
		\hline
		&\multicolumn{2}{c|}{Dof per $T=4$}  &\multicolumn{3}{c|}{Dof per $T=4$} &\multicolumn{3}{c|}{Dof per $T=4$}\\
		\hline
		$\beta_0$ &$\textrm{Erp}(p_{\textrm{off}})$ &$\textrm{Eru}(\mathbf{u}_{\textrm{off}})$ &$\textrm{Erp}(\hat{p}_{\textrm{off}})$ &$\textrm{Eru}(\hat{\mathbf{u}}_{\textrm{off}})$ &$N_{\textrm{update}}$ &$\textrm{Erp}(\tilde{p}_{\textrm{off}})$ &$\textrm{Eru}(\tilde{\mathbf{u}}_{\textrm{off}})$ &$N_{\textrm{update}}$\\
		\hline
		0         &0.0053      &0.0596               &-           &-       &-                         &-                 &-       &-  \\
		1         &0.0050      &0.0778               &0.0051      &0.0747  &24                        &0.0051            &0.0753  &100\\
		10        &0.0051      &0.1032               &0.0050      &0.0901  &29                        &0.0050            &0.0858  &100\\
		100       &0.0063      &0.1341               &0.0055      &0.0975  &33                        &0.0052            &0.0950  &100\\
		1000      &0.0075      &0.1507               &0.0062      &0.1036  &34                        &0.0055            &0.1004  &100\\
		10000     &0.0080      &0.1566               &0.0065      &0.1065  &35                        &0.0056            &0.1025  &100\\
		\hline
	\end{tabular*}
\end{table}
\begin{table}[h!]
	\caption{(Example 1) Relative errors $\textrm{Erp}(p_{\textrm{off}})$, $\textrm{Eru}(\mathbf{u}_{\textrm{off}})$, $\textrm{Erp}(\hat{p}_{\textrm{off}})$, $\textrm{Eru}(\hat{\mathbf{u}}_{\textrm{off}})$, $\textrm{Erp}(\tilde{p}_{\textrm{off}})$ and $\textrm{Eru}(\tilde{\mathbf{u}}_{\textrm{off}})$ with $6$ offline basis functions per coarse element, $\theta=3/4$.}\label{tab_err_ex1_6}
	\centering
	\begin{tabular*}
		{0.92\textwidth}{@{\extracolsep{\fill}}|c|c|c|c|c|c|c|c|c|}
		\hline
		&\multicolumn{2}{c|}{Dof per $T=6$} &\multicolumn{3}{c|}{Dof per $T=6$} &\multicolumn{3}{c|}{Dof per $T=6$}\\
		\hline
		$\beta_0$ &$\textrm{Erp}(p_{\textrm{off}})$ &$\textrm{Eru}(\mathbf{u}_{\textrm{off}})$ &$\textrm{Erp}(\hat{p}_{\textrm{off}})$ &$\textrm{Eru}(\hat{\mathbf{u}}_{\textrm{off}})$ &$N_{\textrm{update}}$ &$\textrm{Erp}(\tilde{p}_{\textrm{off}})$ &$\textrm{Eru}(\tilde{\mathbf{u}}_{\textrm{off}})$ &$N_{\textrm{update}}$\\
		\hline
		0         &0.0024      &0.0058               &-           &-       &-                         &-                 &-       &-\\
		1         &0.0024      &0.0268               &0.0023      &0.0148  &24                        &0.0024            &0.0143  &100\\
		10        &0.0028      &0.0485               &0.0024      &0.0212  &31                        &0.0025            &0.0217  &100\\
		100       &0.0041      &0.0753               &0.0033      &0.0317  &36                        &0.0027            &0.0286  &100\\
		1000      &0.0053      &0.0910               &0.0039      &0.0351  &38                        &0.0030            &0.0323  &100\\
		10000     &0.0057      &0.0970               &0.0042      &0.0370  &38                        &0.0031            &0.0337  &100\\
		\hline
	\end{tabular*}
\end{table}
\begin{table}[h!]
	\caption{(Example 1) Relative errors $\textrm{Erp}(p_{\textrm{off}})$, $\textrm{Eru}(\mathbf{u}_{\textrm{off}})$, $\textrm{Erp}(\hat{p}_{\textrm{off}})$, $\textrm{Eru}(\hat{\mathbf{u}}_{\textrm{off}})$, $\textrm{Erp}(\tilde{p}_{\textrm{off}})$ and $\textrm{Eru}(\tilde{\mathbf{u}}_{\textrm{off}})$ with $8$ offline basis functions per coarse element, $\theta=3/4$.}\label{tab_err_ex1_8}
	\centering
	\begin{tabular*}
		{0.92\textwidth}{@{\extracolsep{\fill}}|c|c|c|c|c|c|c|c|c|}
		\hline
		&\multicolumn{2}{c|}{Dof per $T=8$} &\multicolumn{3}{c|}{Dof per $T=8$} &\multicolumn{3}{c|}{Dof per $T=8$}\\
		\hline
		$\beta_0$ &$\textrm{Erp}(p_{\textrm{off}})$ &$\textrm{Eru}(\mathbf{u}_{\textrm{off}})$ &$\textrm{Erp}(\hat{p}_{\textrm{off}})$ &$\textrm{Eru}(\hat{\mathbf{u}}_{\textrm{off}})$ &$N_{\textrm{update}}$ &$\textrm{Erp}(\tilde{p}_{\textrm{off}})$ &$\textrm{Eru}(\tilde{\mathbf{u}}_{\textrm{off}})$ &$N_{\textrm{update}}$\\
		\hline
		0         &0.0013      &0.0011               &-           &-       &-           &-           &-       &-\\
		1         &0.0014      &0.0233               &0.0014      &0.0111  &17          &0.0014      &0.0100  &100\\
		10        &0.0018      &0.0364               &0.0016      &0.0158  &25          &0.0015      &0.0135  &100\\
		100       &0.0028      &0.0466               &0.0023      &0.0199  &31          &0.0017      &0.0164  &100\\
		1000      &0.0035      &0.0522               &0.0028      &0.0227  &34          &0.0019      &0.0181  &100\\
		10000     &0.0039      &0.0545               &0.0031      &0.0241  &34          &0.0020      &0.0187  &100\\
		\hline
	\end{tabular*}
\end{table}
\begin{figure}[h!]
	\mbox{\hspace{-0.65cm}}
	\begin{minipage}[b]{0.34\textwidth}
		\centering
		\includegraphics[scale=0.45]{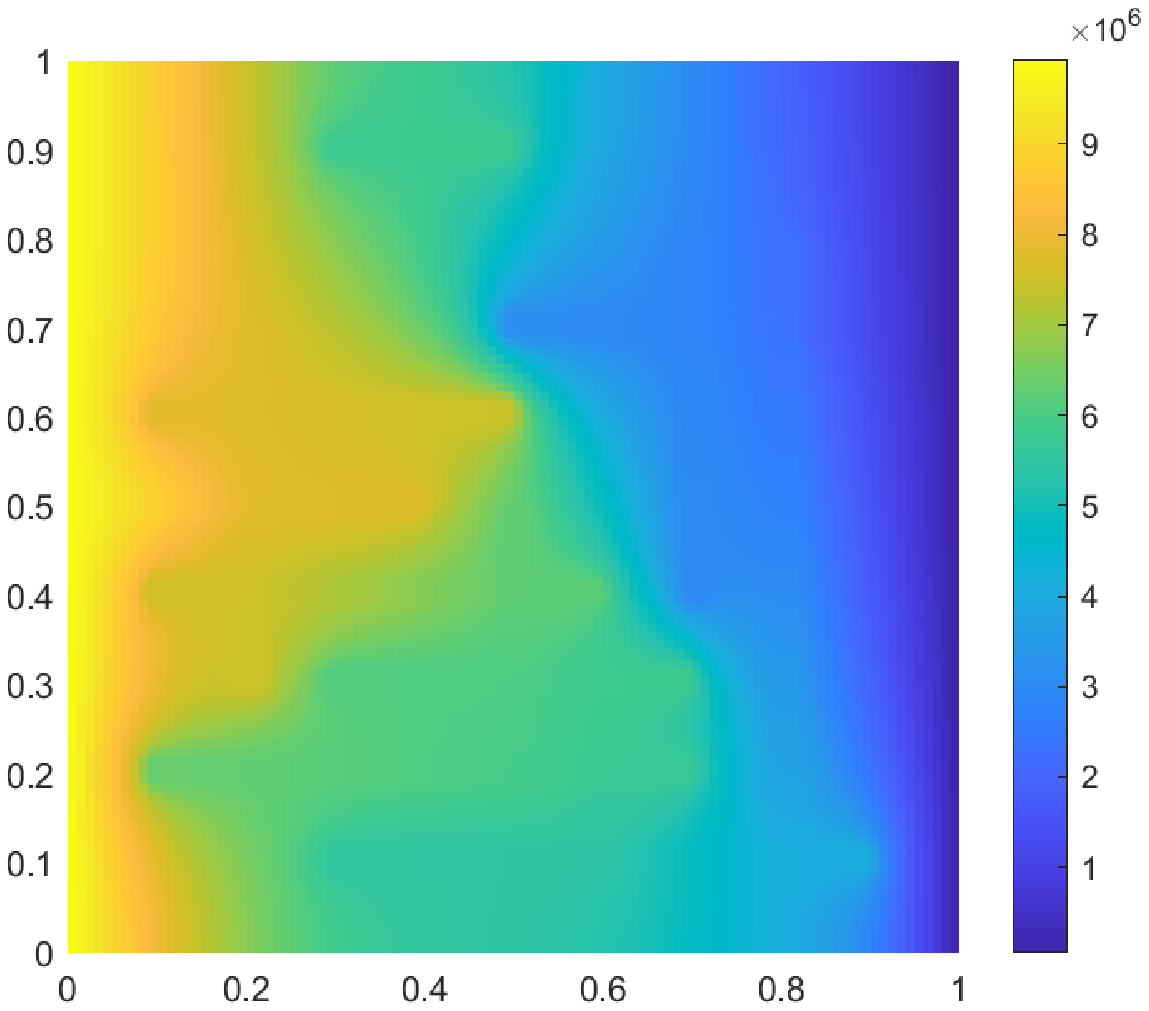}
	\end{minipage}
	\mbox{\hspace{0.00cm}}
	\begin{minipage}[b]{0.34\textwidth}
		\centering
		\includegraphics[scale=0.45]{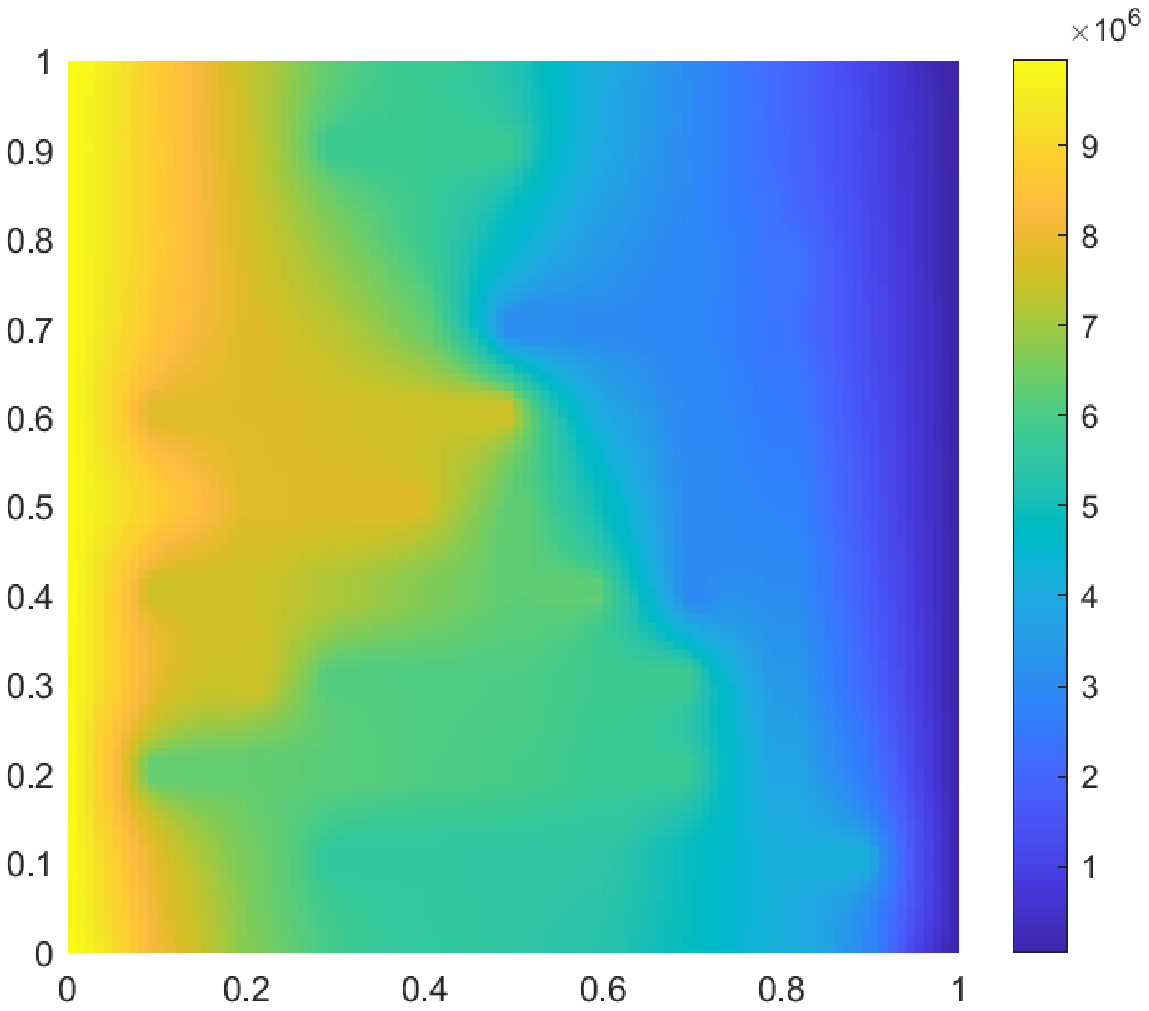}
	\end{minipage}
	\mbox{\hspace{0.00cm}}
	\begin{minipage}[b]{0.34\textwidth}
		\centering
		\includegraphics[scale=0.45]{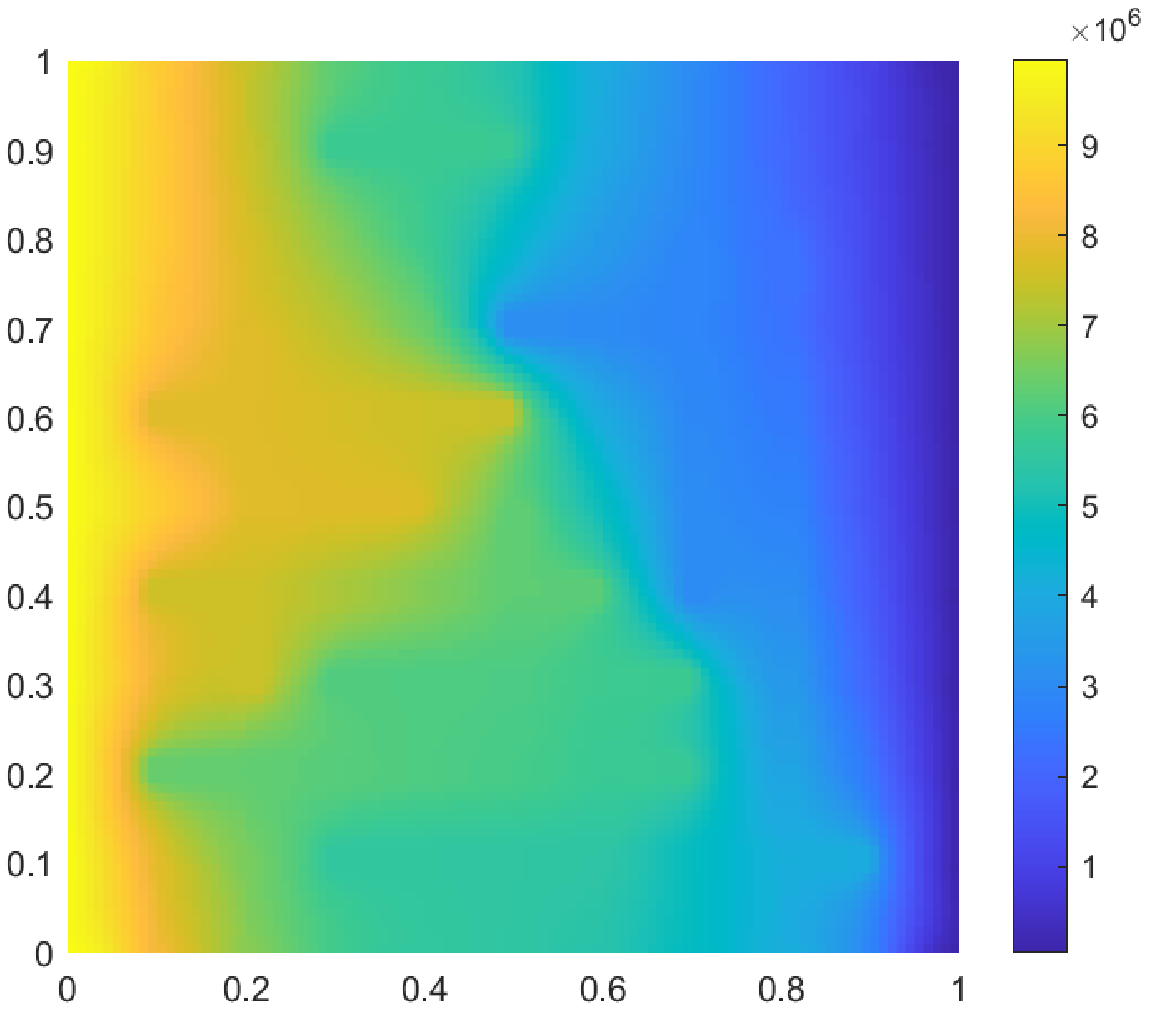}
	\end{minipage}
	\\
	\mbox{\hspace{-0.65cm}}
	\begin{minipage}[b]{0.34\textwidth}
		\centering
		\includegraphics[scale=0.45]{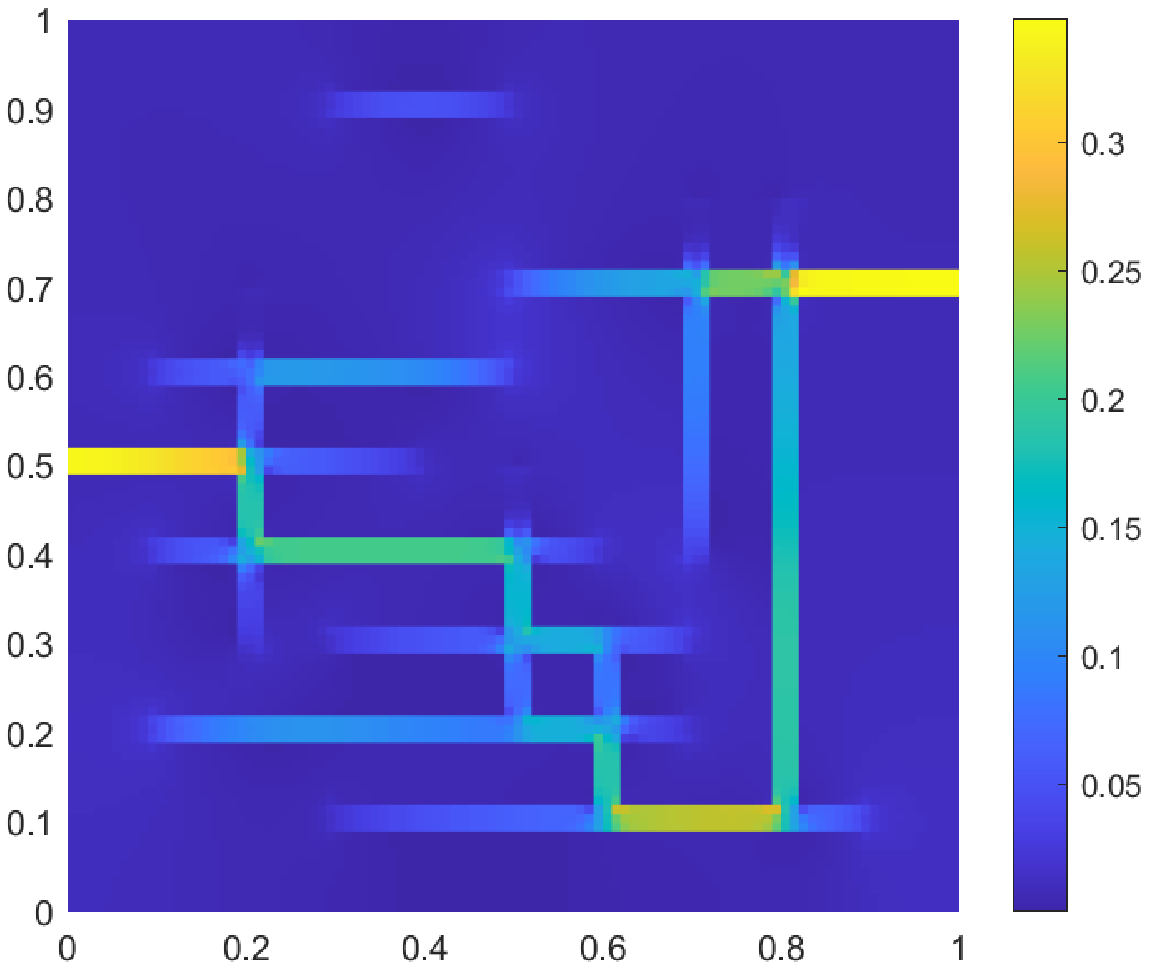}
	\end{minipage}
	\mbox{\hspace{0.00cm}}
	\begin{minipage}[b]{0.34\textwidth}
		\centering
		\includegraphics[scale=0.45]{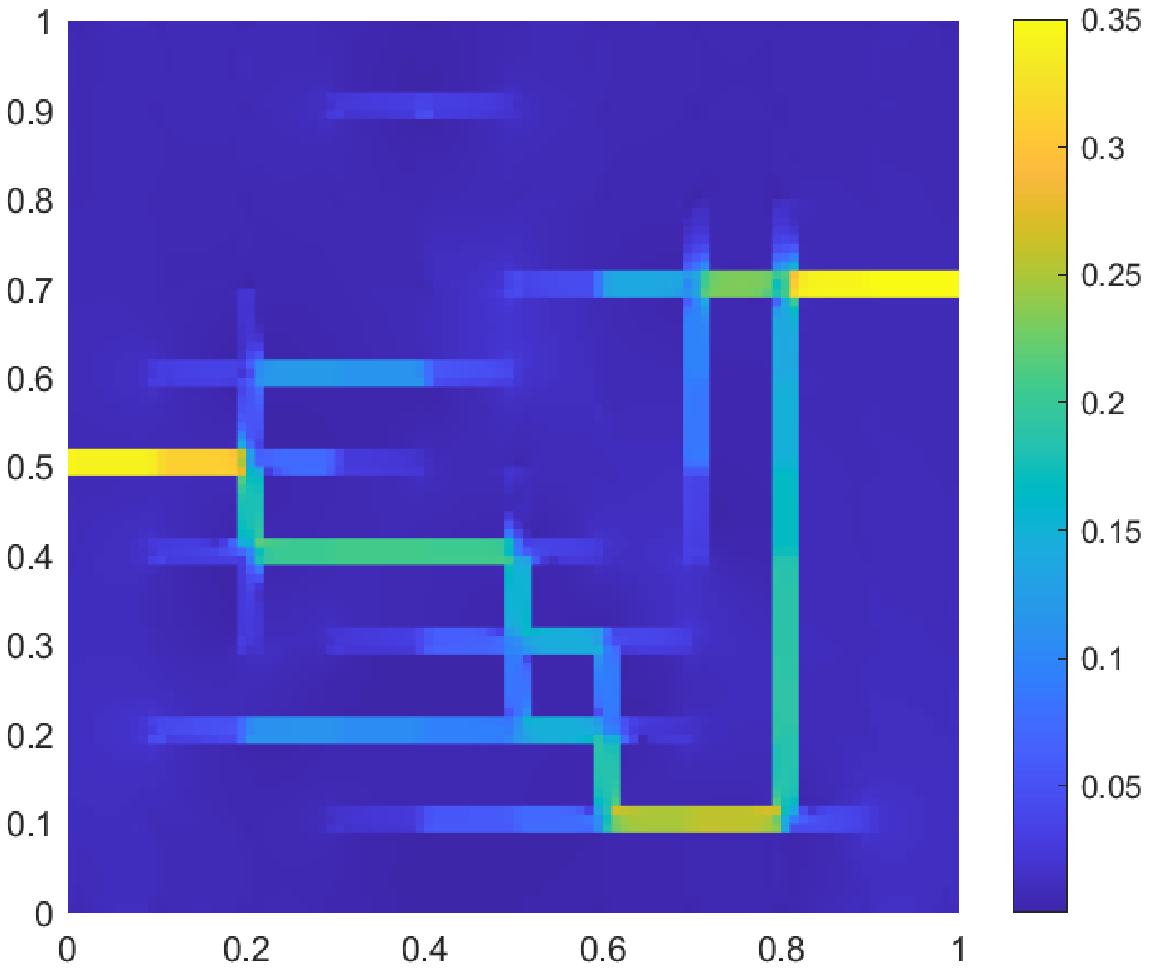}
	\end{minipage}
	\mbox{\hspace{0.00cm}}
	\begin{minipage}[b]{0.34\textwidth}
		\centering
		\includegraphics[scale=0.45]{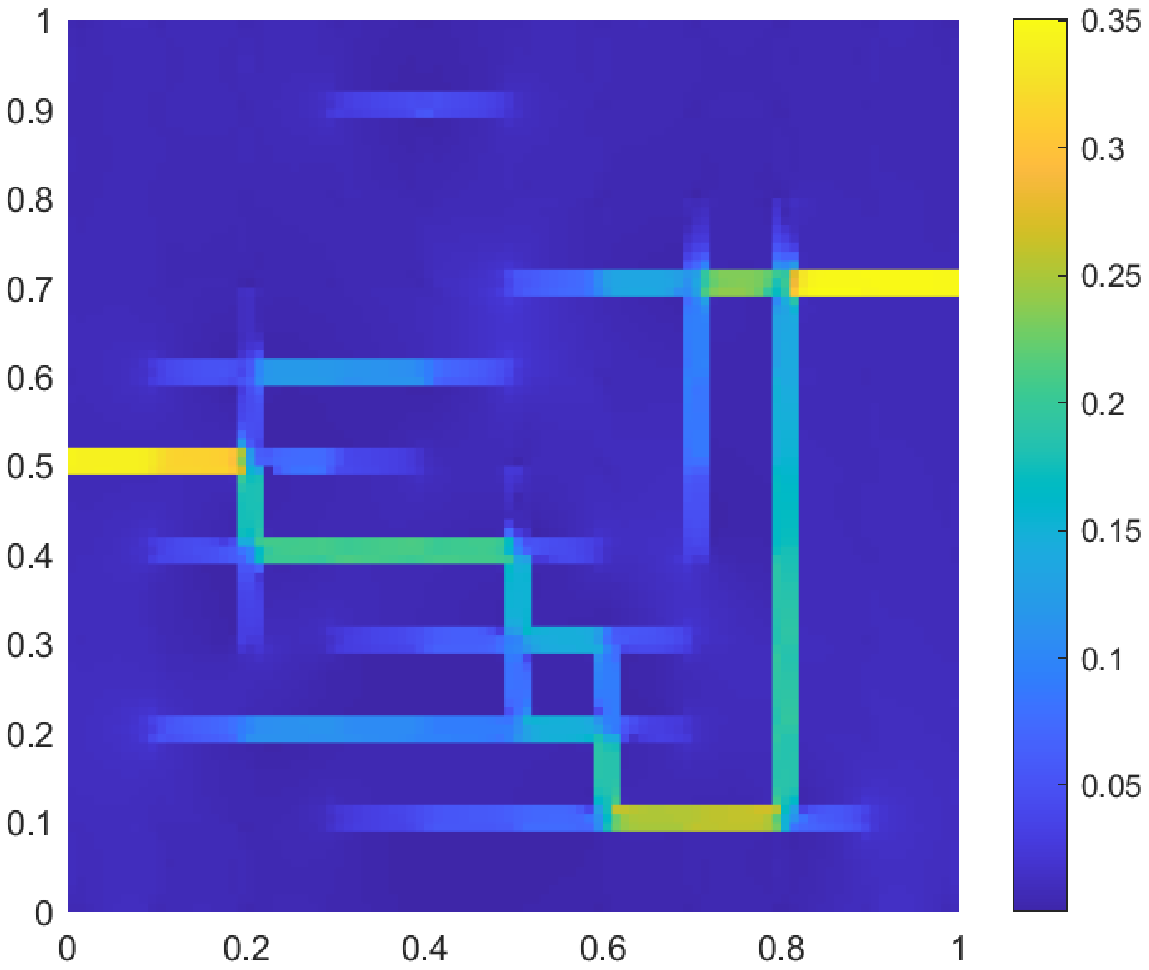}
	\end{minipage}
	\caption{(Example 1) The pressure and the velocity using $4$ offline basis functions per coarse element with $\beta_0=100$ and $\theta=3/4$. Left: fine-grid solution $(p_f,\mathbf{u}_f)$. Middle: offline solution $(p_{\textrm{off}},\mathbf{u}_{\textrm{off}})$. Right: updated offline solution $(\hat{p}_{\textrm{off}},\hat{\mathbf{u}}_{\textrm{off}})$.}\label{fig_pressurevelocity_1}
\end{figure}

{\bf Example 2 :}
\begin{figure}[h!]
	\centering
	\includegraphics[width=0.78 \textwidth]{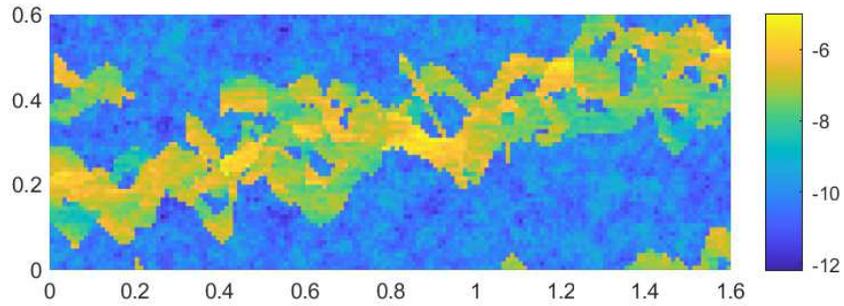}
	\caption{The distribution of permeability field $\kappa$ in logarithmic scale for example 2.}
	\label{fig_perm_2}
\end{figure}
The computational domain is set to be $\Omega=[0, 1.6]\times[0, 0.6]$, the fine grid is a $160\times60$ uniform mesh, and the coarse grid is a $16\times6$ uniform mesh. The permeability field $\kappa$ is a part of the horizontal permeability from the SPE10 data set, as shown in Figure \ref{fig_perm_2}.

Table \ref{tab_iteration_number_2} shows the iteration number of Picard iterative algorithm (\ref{eqn_Picard_1_H})-(\ref{eqn_Picard_2_H}) and Newton iterative algorithm (\ref{eqn_Newton_1_H})-(\ref{eqn_Newton_2_H}) with parameter $\beta_0$ taken different values, as the result in example 1, the number of Newton iterations is much less than the Picard iterations.
\begin{table}[h!]
	\caption{(Example 2) The number of iterations.}\label{tab_iteration_number_2}
	\centering
	\begin{tabular*}
		{0.5\textwidth}{@{\extracolsep{\fill}}|c|c|c|}
		\hline
		$\beta_0$  & $\#$Newton iterations  & $\#$Picard iterations\\
		\hline
		1                &7                       &56   \\
		10               &9                       &182  \\
		100              &10                      &522  \\
		1000             &12                      &1501  \\
		10000            &14                      &4257  \\
		\hline
	\end{tabular*}
\end{table}

Tables \ref{tab_err_ex2_4}, \ref{tab_err_ex2_6}, and \ref{tab_err_ex2_8} show the $L_2$ relative errors $\textrm{Erp}(p_{\textrm{off}})$, $\textrm{Eru}(\mathbf{u}_{\textrm{off}})$, $\textrm{Erp}(\hat{p}_{\textrm{off}})$, $\textrm{Eru}(\hat{\mathbf{u}}_{\textrm{off}})$, $\textrm{Erp}(\tilde{p}_{\textrm{off}})$ and $\textrm{Eru}(\tilde{\mathbf{u}}_{\textrm{off}})$ with the parameter $\beta_0$ taken different values with respect to $4$, $6$ and $8$ offline basis functions per coarse element, respectively. The results are similar to that of Example 1, although the relative errors in Tables \ref{tab_err_ex2_4}, \ref{tab_err_ex2_6}, and \ref{tab_err_ex2_8} are a little bit bigger than the relative errors in Tables \ref{tab_err_ex1_4}, \ref{tab_err_ex1_6}, and \ref{tab_err_ex1_8} of Example 1, it seems that the offline basis functions remain effective for the permeability filed $\kappa$ shown in Figure \ref{fig_perm_2}.
\begin{table}[h!]
	\caption{(Example 2) Relative errors $\textrm{Erp}(p_{\textrm{off}})$, $\textrm{Eru}(\mathbf{u}_{\textrm{off}})$, $\textrm{Erp}(\hat{p}_{\textrm{off}})$, $\textrm{Eru}(\hat{\mathbf{u}}_{\textrm{off}})$, $\textrm{Erp}(\tilde{p}_{\textrm{off}})$ and $\textrm{Eru}(\tilde{\mathbf{u}}_{\textrm{off}})$ with $4$ offline basis functions per coarse element, $\theta=3/4$.}\label{tab_err_ex2_4}
	\centering
	\begin{tabular*}
		{0.92\textwidth}{@{\extracolsep{\fill}}|c|c|c|c|c|c|c|c|c|}
		\hline
		&\multicolumn{2}{c|}{Dof per $T=4$}  &\multicolumn{3}{c|}{Dof per $T=4$} &\multicolumn{3}{c|}{Dof per $T=4$}\\
		\hline
		$\beta_0$ &$\textrm{Erp}(p_{\textrm{off}})$ &$\textrm{Eru}(\mathbf{u}_{\textrm{off}})$ &$\textrm{Erp}(\hat{p}_{\textrm{off}})$ &$\textrm{Eru}(\hat{\mathbf{u}}_{\textrm{off}})$ &$N_{\textrm{update}}$ &$\textrm{Erp}(\tilde{p}_{\textrm{off}})$ &$\textrm{Eru}(\tilde{\mathbf{u}}_{\textrm{off}})$ &$N_{\textrm{update}}$\\
		\hline
		0         &0.0091      &0.0891               &-           &-       &-                         &-                 &-       &- \\
		1         &0.0092      &0.0975               &0.0091      &0.0923  &26                        &0.0091            &0.0921  &96\\
		10        &0.0095      &0.1270               &0.0090      &0.1076  &27                        &0.0088            &0.1045  &96\\
		100       &0.0091      &0.1594               &0.0083      &0.1290  &29                        &0.0078            &0.1225  &96\\
		1000      &0.0084      &0.1773               &0.0075      &0.1445  &30                        &0.0069            &0.1370  &96\\
		10000     &0.0081      &0.1846               &0.0078      &0.1542  &31                        &0.0074            &0.1463  &96\\
		\hline
	\end{tabular*}
\end{table}
\begin{table}[h!]
	\caption{(Example 2) Relative errors $\textrm{Erp}(p_{\textrm{off}})$, $\textrm{Eru}(\mathbf{u}_{\textrm{off}})$, $\textrm{Erp}(\hat{p}_{\textrm{off}})$, $\textrm{Eru}(\hat{\mathbf{u}}_{\textrm{off}})$, $\textrm{Erp}(\tilde{p}_{\textrm{off}})$ and $\textrm{Eru}(\tilde{\mathbf{u}}_{\textrm{off}})$ with $4$ offline basis functions per coarse element, $\theta=3/4$.}\label{tab_err_ex2_6}
	\centering
	\begin{tabular*}
		{0.92\textwidth}{@{\extracolsep{\fill}}|c|c|c|c|c|c|c|c|c|}
		\hline
		&\multicolumn{2}{c|}{Dof per $T=6$}  &\multicolumn{3}{c|}{Dof per $T=6$} &\multicolumn{3}{c|}{Dof per $T=6$}\\
		\hline
		$\beta_0$ &$\textrm{Erp}(p_{\textrm{off}})$ &$\textrm{Eru}(\mathbf{u}_{\textrm{off}})$ &$\textrm{Erp}(\hat{p}_{\textrm{off}})$ &$\textrm{Eru}(\hat{\mathbf{u}}_{\textrm{off}})$ &$N_{\textrm{update}}$ &$\textrm{Erp}(\tilde{p}_{\textrm{off}})$ &$\textrm{Eru}(\tilde{\mathbf{u}}_{\textrm{off}})$ &$N_{\textrm{update}}$\\
		\hline
		0         &0.0018      &0.0354               &-           &-       &-                         &-                 &-       &- \\
		1         &0.0018      &0.0449               &0.0018      &0.0363  &19                        &0.0018            &0.0351  &96\\
		10        &0.0022      &0.0820               &0.0018      &0.0489  &22                        &0.0017            &0.0409  &96\\
		100       &0.0026      &0.1157               &0.0027      &0.0607  &24                        &0.0017            &0.0509  &96\\
		1000      &0.0029      &0.1324               &0.0029      &0.0694  &26                        &0.0018            &0.0589  &96\\
		10000     &0.0030      &0.1382               &0.0034      &0.0730  &27                        &0.0019            &0.0627  &96\\
		\hline
	\end{tabular*}
\end{table}
\begin{table}[h!]
	\caption{(Example 2) Relative errors $\textrm{Erp}(p_{\textrm{off}})$, $\textrm{Eru}(\mathbf{u}_{\textrm{off}})$, $\textrm{Erp}(\hat{p}_{\textrm{off}})$, $\textrm{Eru}(\hat{\mathbf{u}}_{\textrm{off}})$, $\textrm{Erp}(\tilde{p}_{\textrm{off}})$ and $\textrm{Eru}(\tilde{\mathbf{u}}_{\textrm{off}})$ with $4$ offline basis functions per coarse element, $\theta=3/4$.}\label{tab_err_ex2_8}
	\centering
	\begin{tabular*}
		{0.92\textwidth}{@{\extracolsep{\fill}}|c|c|c|c|c|c|c|c|c|}
		\hline
		&\multicolumn{2}{c|}{Dof per $T=8$}  &\multicolumn{3}{c|}{Dof per $T=8$} &\multicolumn{3}{c|}{Dof per $T=8$}\\
		\hline
		$\beta_0$ &$\textrm{Erp}(p_{\textrm{off}})$ &$\textrm{Eru}(\mathbf{u}_{\textrm{off}})$ &$\textrm{Erp}(\hat{p}_{\textrm{off}})$ &$\textrm{Eru}(\hat{\mathbf{u}}_{\textrm{off}})$ &$N_{\textrm{update}}$ &$\textrm{Erp}(\tilde{p}_{\textrm{off}})$ &$\textrm{Eru}(\tilde{\mathbf{u}}_{\textrm{off}})$ &$N_{\textrm{update}}$\\
		\hline
		0         &0.0007      &0.0201               &-           &-       &-                         &-                 &-       &- \\
		1         &0.0009      &0.0345               &0.0007      &0.0219  &18                        &0.0007            &0.0203  &96\\
		10        &0.0016      &0.0734               &0.0009      &0.0364  &20                        &0.0007            &0.0284  &96\\
		100       &0.0024      &0.1051               &0.0018      &0.0502  &23                        &0.0008            &0.0400  &96\\
		1000      &0.0027      &0.1208               &0.0020      &0.0599  &24                        &0.0009            &0.0471  &96\\
		10000     &0.0029      &0.1264               &0.0026      &0.0610  &26                        &0.0010            &0.0497  &96\\
		\hline
	\end{tabular*}
\end{table}

In addition, Figure \ref{fig_pressurevelocity_2} plots the fine-grid solutions $(p_f,\mathbf{u}_f)$, offline solution $(p_\textrm{off},\mathbf{u}_{\textrm{off}})$ and updated offline solution $(\hat{p}_{\textrm{off}},\hat{\mathbf{u}}_{\textrm{off}})$, respectively, with $4$ offline basis functions per coarse element and $\beta_0=100$. Obviously, the offline solutions still have competitive performances.
\begin{figure}[h!]
	\mbox{\hspace{-0.90cm}}
	\begin{minipage}[b]{0.52\textwidth}
		\centering
		\includegraphics[scale=0.66]{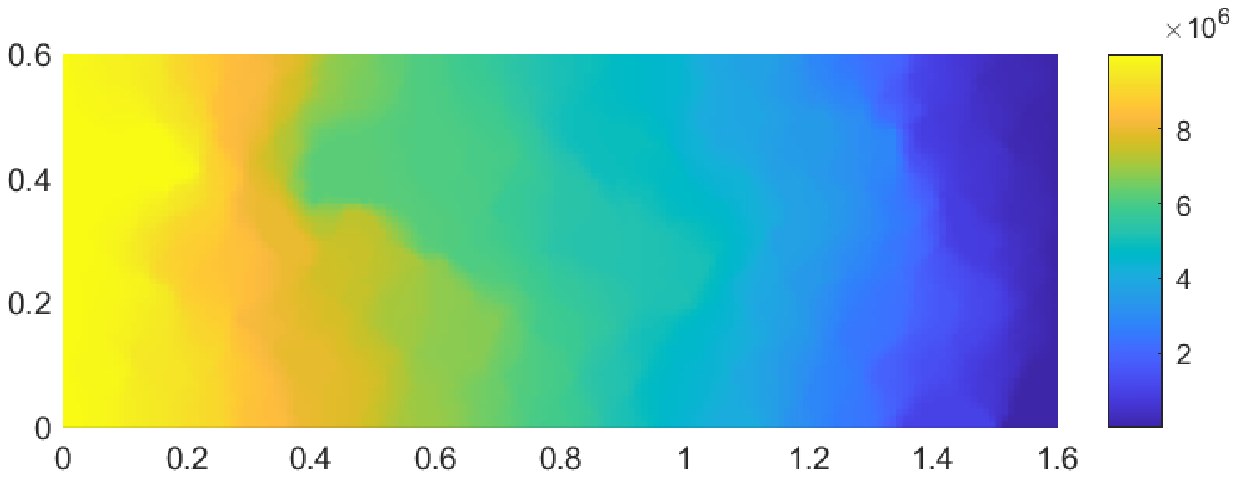}
	\end{minipage}
	\mbox{\hspace{0.00cm}}
	\begin{minipage}[b]{0.52\textwidth}
		\centering
		\includegraphics[scale=0.66]{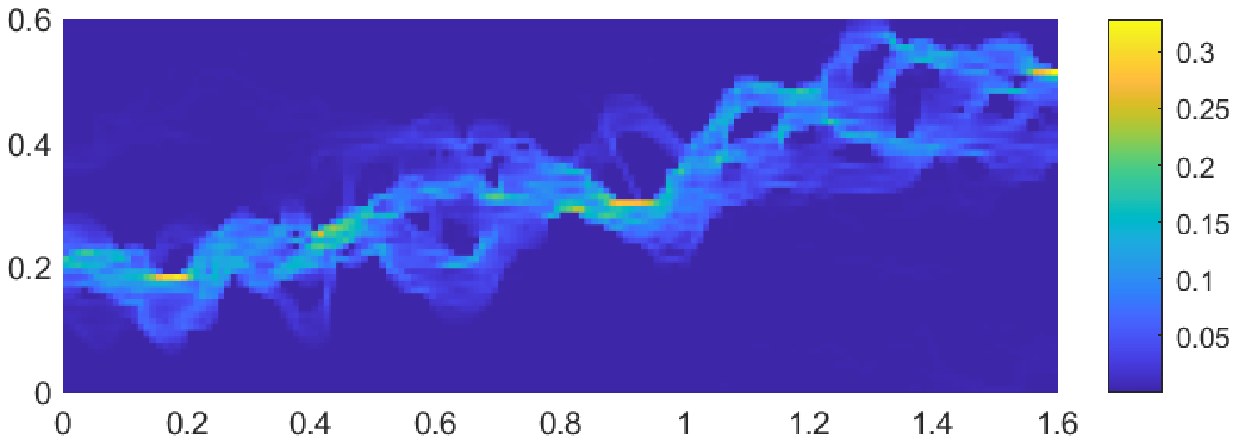}
	\end{minipage}
	\mbox{\hspace{-0.90cm}}
	\begin{minipage}[b]{0.52\textwidth}
		\centering
		\includegraphics[scale=0.66]{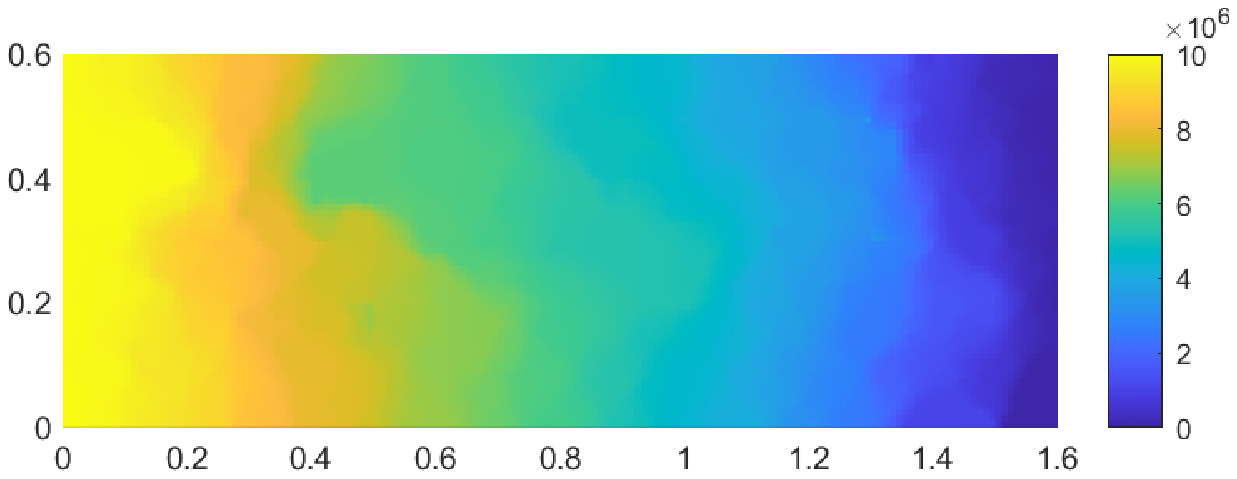}
	\end{minipage}
	\mbox{\hspace{0.00cm}}
	\begin{minipage}[b]{0.52\textwidth}
		\centering
		\includegraphics[scale=0.66]{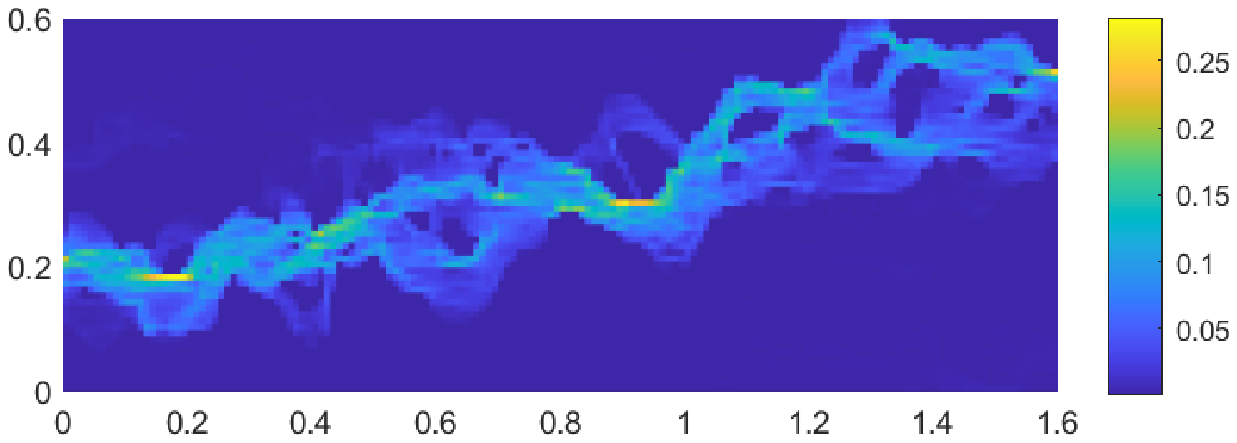}
	\end{minipage}
	\mbox{\hspace{-0.90cm}}
	\begin{minipage}[b]{0.52\textwidth}
		\centering
		\includegraphics[scale=0.66]{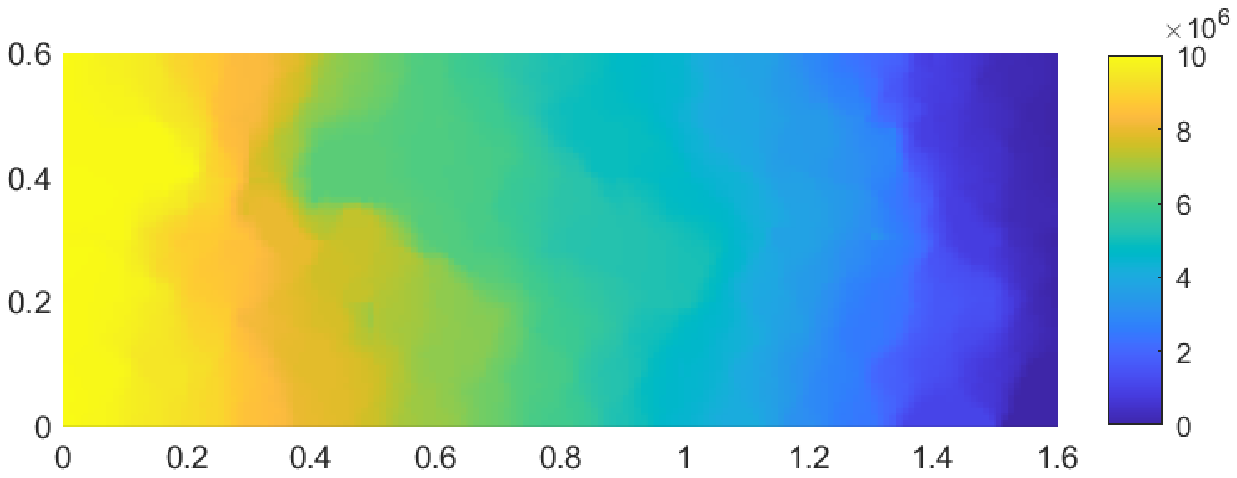}
	\end{minipage}
	\mbox{\hspace{0.00cm}}
	\begin{minipage}[b]{0.52\textwidth}
		\centering
		\includegraphics[scale=0.66]{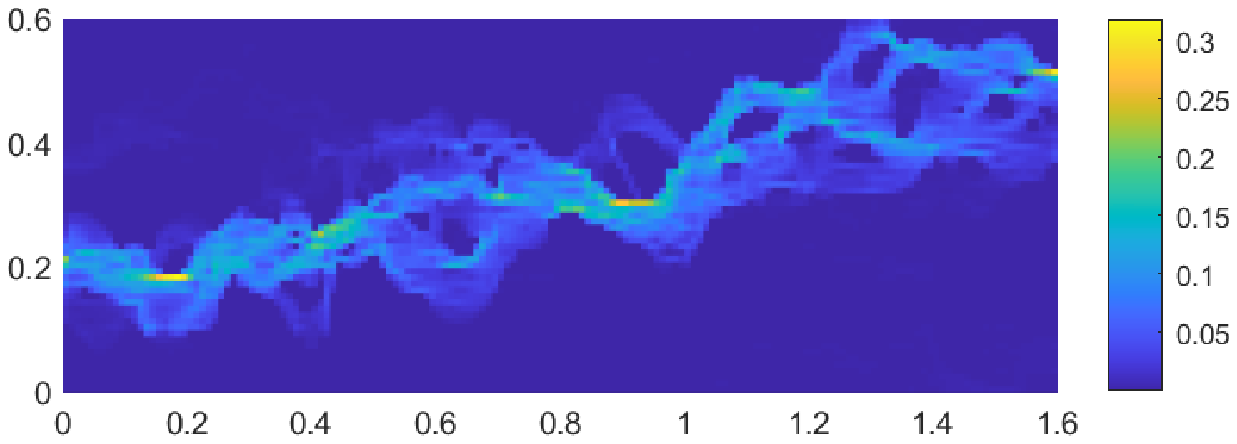}
	\end{minipage}
	\caption{(Example 2) The pressure and the velocity using $4$ offline basis functions per coarse element with $\beta_0=100$ and $\theta=3/4$. Top: fine-grid solution $(p_f,\mathbf{u}_f)$. Middle: offline solution $(p_{\textrm{off}},\mathbf{u}_{\textrm{off}})$. Bottom: updated offline solution $(\hat{p}_{\textrm{off}},\hat{\mathbf{u}}_{\textrm{off}})$.}\label{fig_pressurevelocity_2}
\end{figure}

\subsection{Online computation}
In this subsection, we investigate the performance of the online computation. We will test the multiscale space enrichment uniformly for all coarse elements and adaptively based on residuals, respectively. The online basis functions are calculated and added into the multiscale space only in disjoint regions at a time, for convenience, we use a two-index notation to enumerate all coarse elements, i.e., the coarse elements are indexed by $T_{ij}$, with $i=1,2,\cdots,N_x$ and $j=1,2,\cdots,N_y$, where $N_x$ and $N_y$ are the number of partitions of the coarse grid $\mathcal{T}_H$ along the $x$ and $y$ directions, respectively. Let $I_x=\{1,2,\cdots,N_x\}$ and $I_y=\{1,2,\cdots,N_y\}$. we denote $I_{x,1}$, $I_{x,2}$ be the subsets contain the odd, even indices of $I_x$ respectively, and $I_{y,1}$, $I_{y,2}$ be the subsets contain the odd, even indices of $I_y$ respectively. We can separate the set that contains all coarse elements into four disjoint subsets $S_1$, $S_2$, $S_3$ and $S_4$ with $S_1=I_{x,1}\times I_{y,1}$, $S_2=I_{x,1}\times I_{y,2}$, $S_3=I_{x,2}\times I_{y,1}$ and $S_4=I_{x,2}\times I_{y,2}$, respectively. Each iteration of the online multiscale space enrichment contains four subiterations, in particular, these four subiterations are defined by adding online basis functions with respect to the coarse elements $T_{ij}\in S_1$, $T_{ij}\in S_2$, $T_{ij}\in S_3$ and $T_{ij}\in S_4$, respectively. We use Example 2 in the previous subsection for the numerical test of the online computation.
\subsubsection{Uniform online enrichment}
We first enrich the multiscale space uniformly by adding one online basis function per coarse element in each enrichment iteration. We compare the performance of the uniform online enrichment with the parameter $\beta_0$ taken different values about $3$, $4$, $5$ and $6$ initial basis functions (belong to offline space $W_{\textrm{off}}$) per coarse element, respectively. The results are shown in Figure \ref{fig_velocityerr_1}, where we plot the logarithm of the relative error of velocity against the dimension of the multiscale space $W_{\textrm{ms}}$ in each enrichment iteration. We find that the accuracy of the multiscale solution is improved a lot by several level of online enrichment. The convergence rate becomes slow when $\beta_0$ becomes large. We also compare the performance of the uniform online computation with respect to different number of initial basis functions per coarse element for $\beta_0$ taken different values, the results are shown in Figure \ref{fig_velocityerr_2}, where we plot the relative errors of velocity against dimensions of the multiscale space $W_{\textrm{ms}}$ for $\beta_0$ taken four different values : $\beta_0=10, 100, 10000$ and $10000$, respectively, It can be observed that the performance of the online multiscale space is better than the offline multiscale space with the same dimension.

\begin{figure}
	\centering
	\mbox{\hspace{0.00cm}}
	\begin{minipage}[b]{0.40\textwidth}
		\centering
		\includegraphics[scale=0.48]{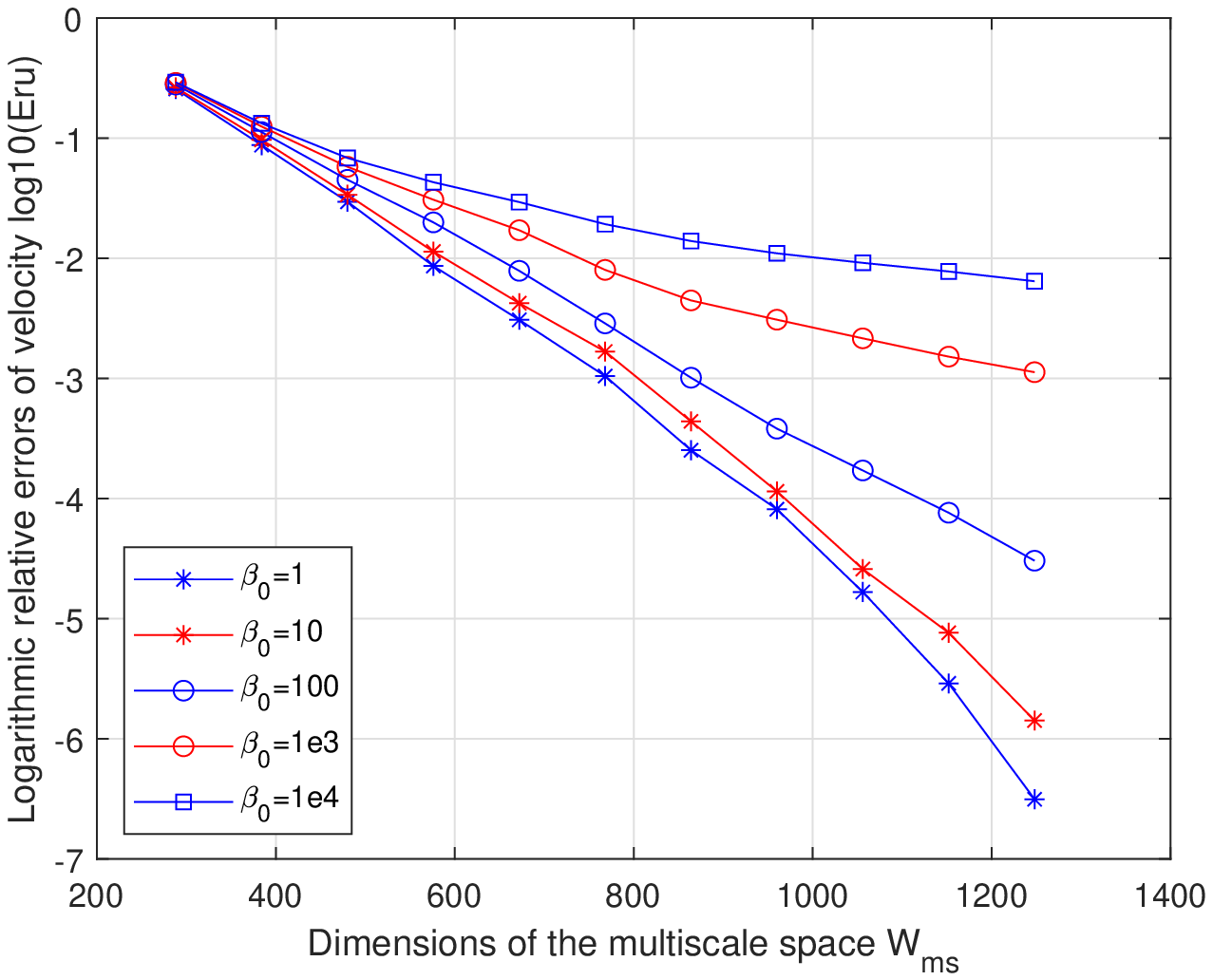}
	\end{minipage}
	\mbox{\hspace{0.00cm}}
	\begin{minipage}[b]{0.40\textwidth}
		\centering
		\includegraphics[scale=0.48]{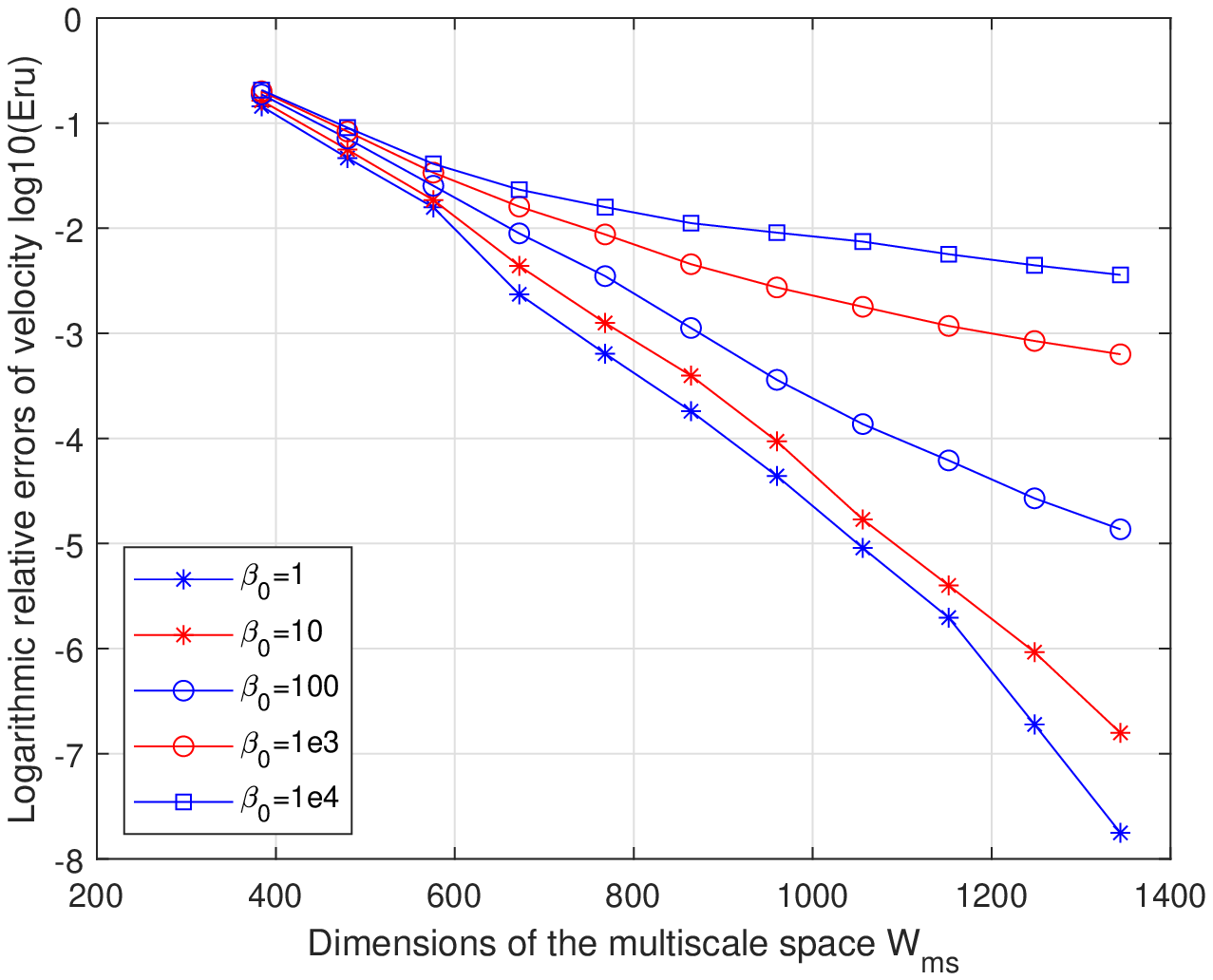}
	\end{minipage}
	\\
	\mbox{\hspace{0.00cm}}
	\begin{minipage}[b]{0.40\textwidth}
		\centering
		\includegraphics[scale=0.48]{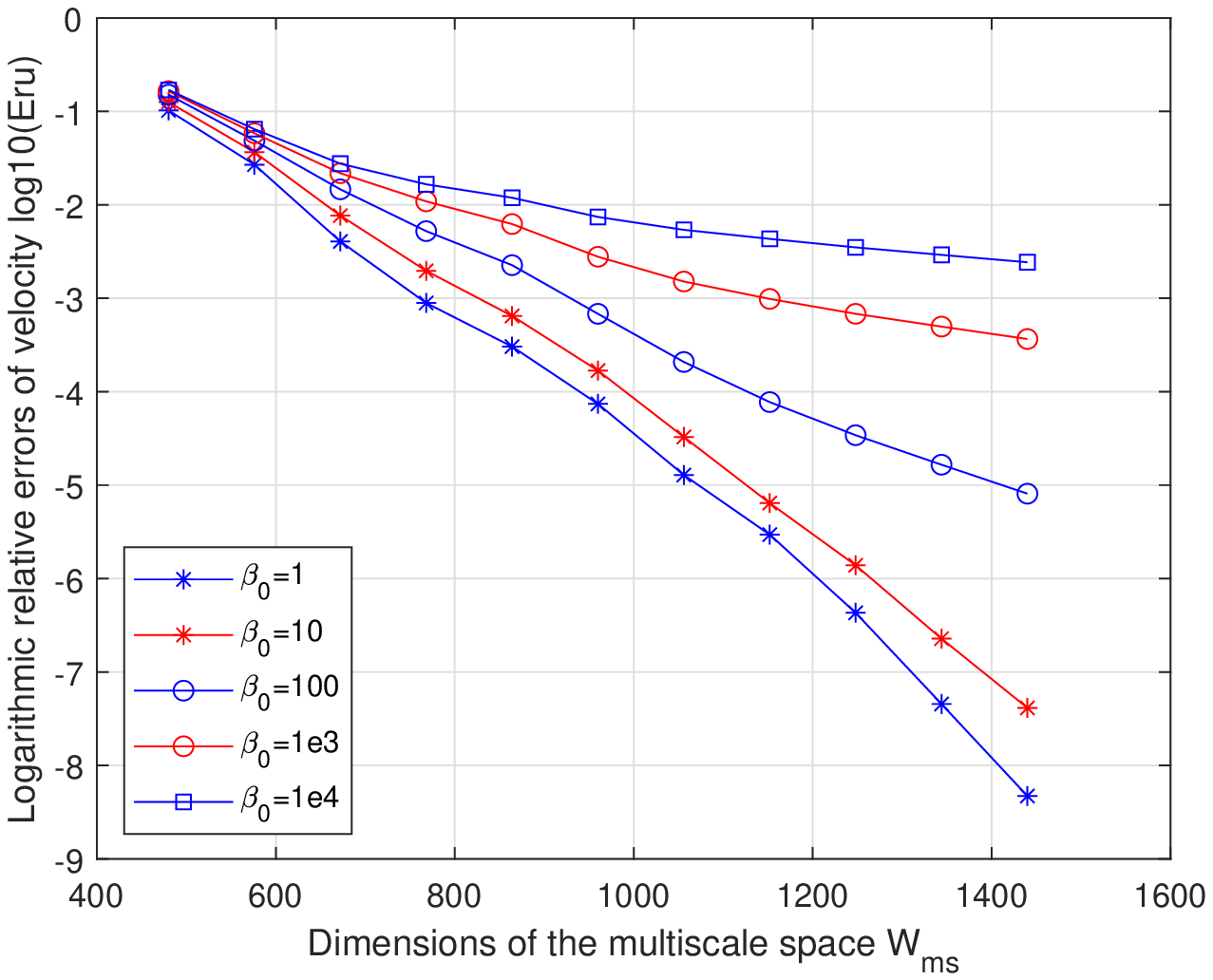}
	\end{minipage}
	\mbox{\hspace{0.00cm}}
	\begin{minipage}[b]{0.40\textwidth}
		\centering
		\includegraphics[scale=0.48]{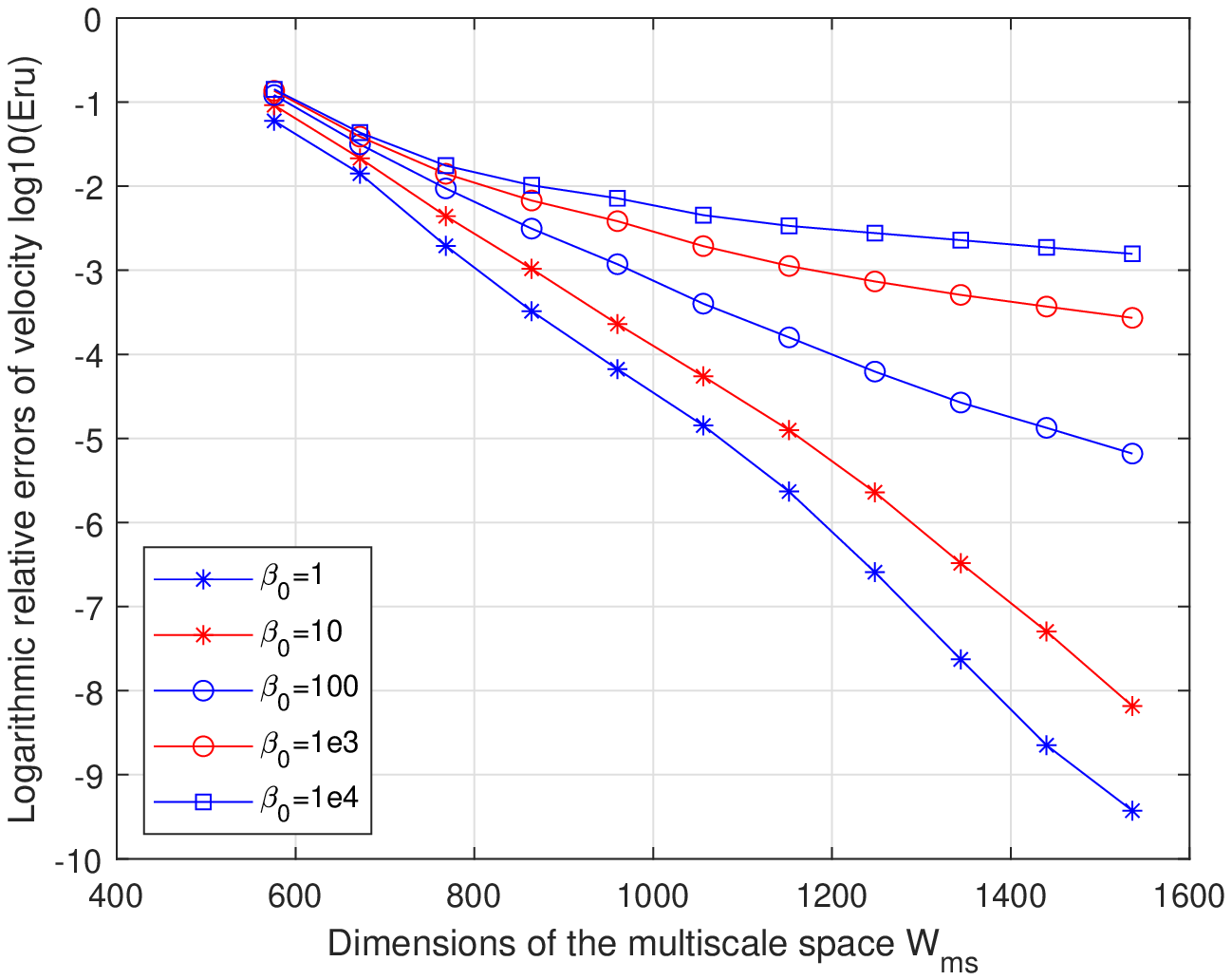}
	\end{minipage}
	\caption{(Uniform online enrichment) Convergence comparison with $\beta_0$ taken different values. Top left: 3 initial basis functions. Top right: 4 initial basis functions. Bottom left: 5 initial basis functions. Bottom right: 6 initial basis functions.}\label{fig_velocityerr_1}
\end{figure}
\begin{figure}
	\centering
	\mbox{\hspace{0.00cm}}
	\begin{minipage}[b]{0.40\textwidth}
		\centering
		\includegraphics[scale=0.48]{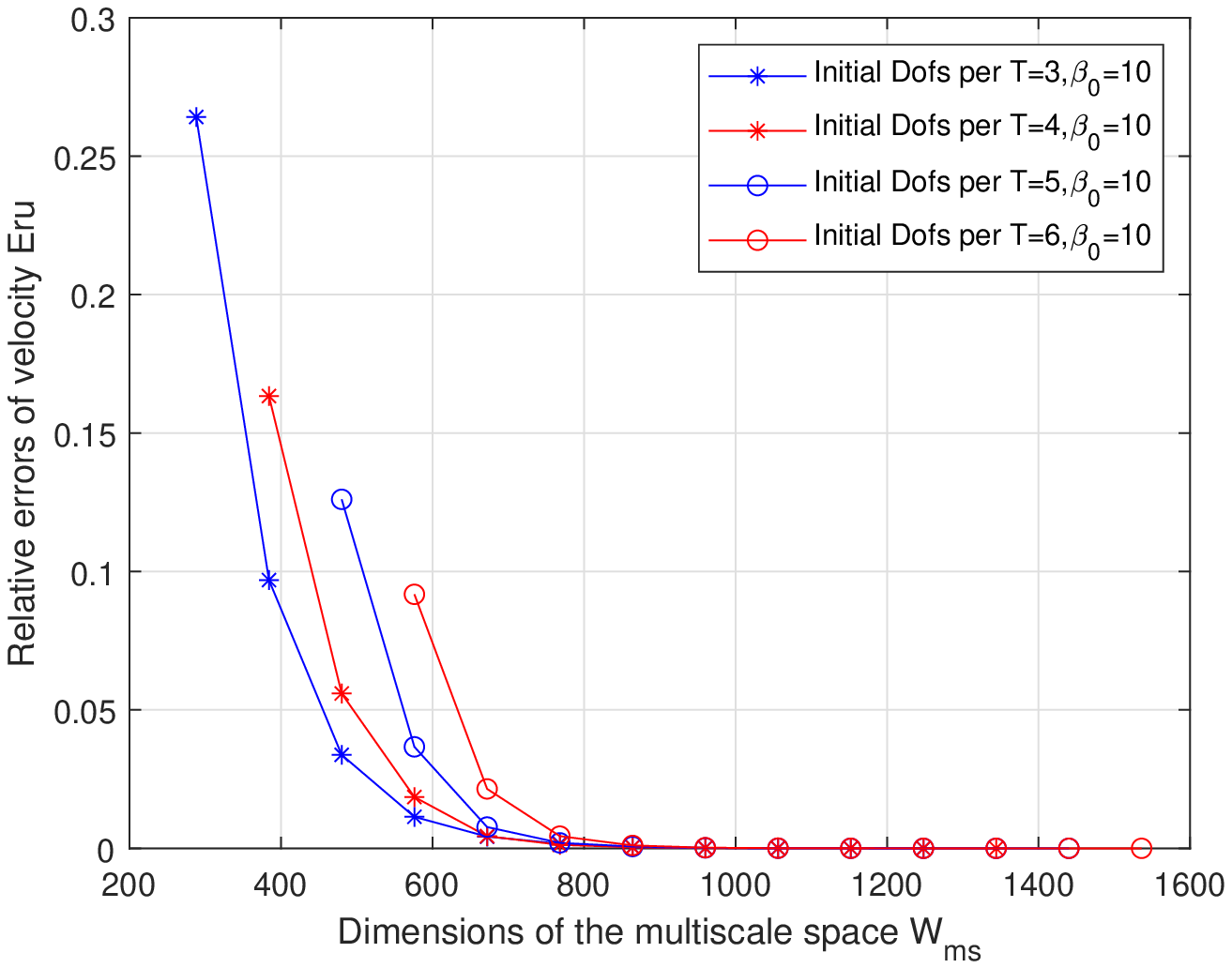}
	\end{minipage}
	\mbox{\hspace{0.00cm}}
	\begin{minipage}[b]{0.40\textwidth}
		\centering
		\includegraphics[scale=0.48]{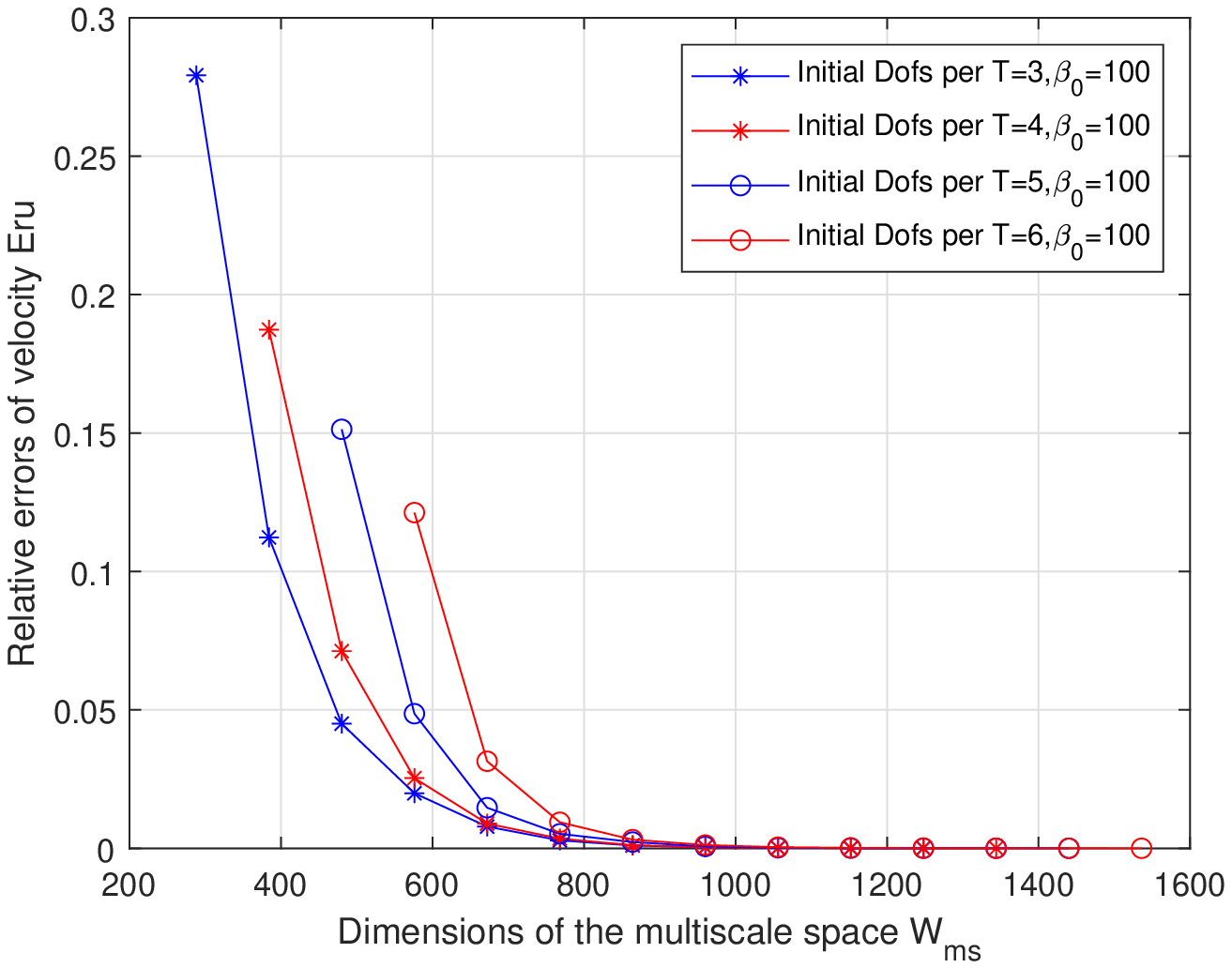}
	\end{minipage}
	\\
	\mbox{\hspace{0.00cm}}
	\begin{minipage}[b]{0.40\textwidth}
		\centering
		\includegraphics[scale=0.48]{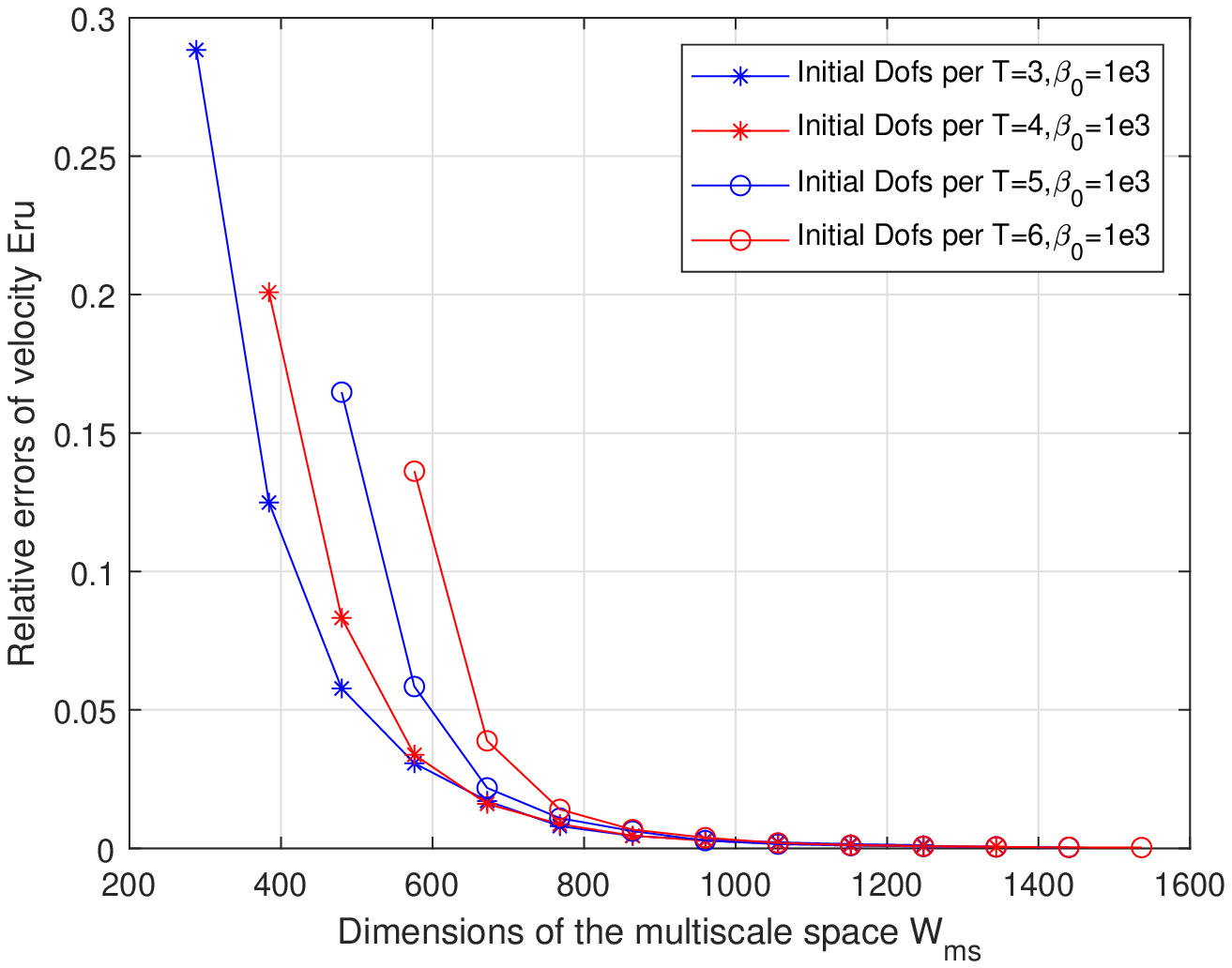}
	\end{minipage}
	\mbox{\hspace{0.00cm}}
	\begin{minipage}[b]{0.40\textwidth}
		\centering
		\includegraphics[scale=0.48]{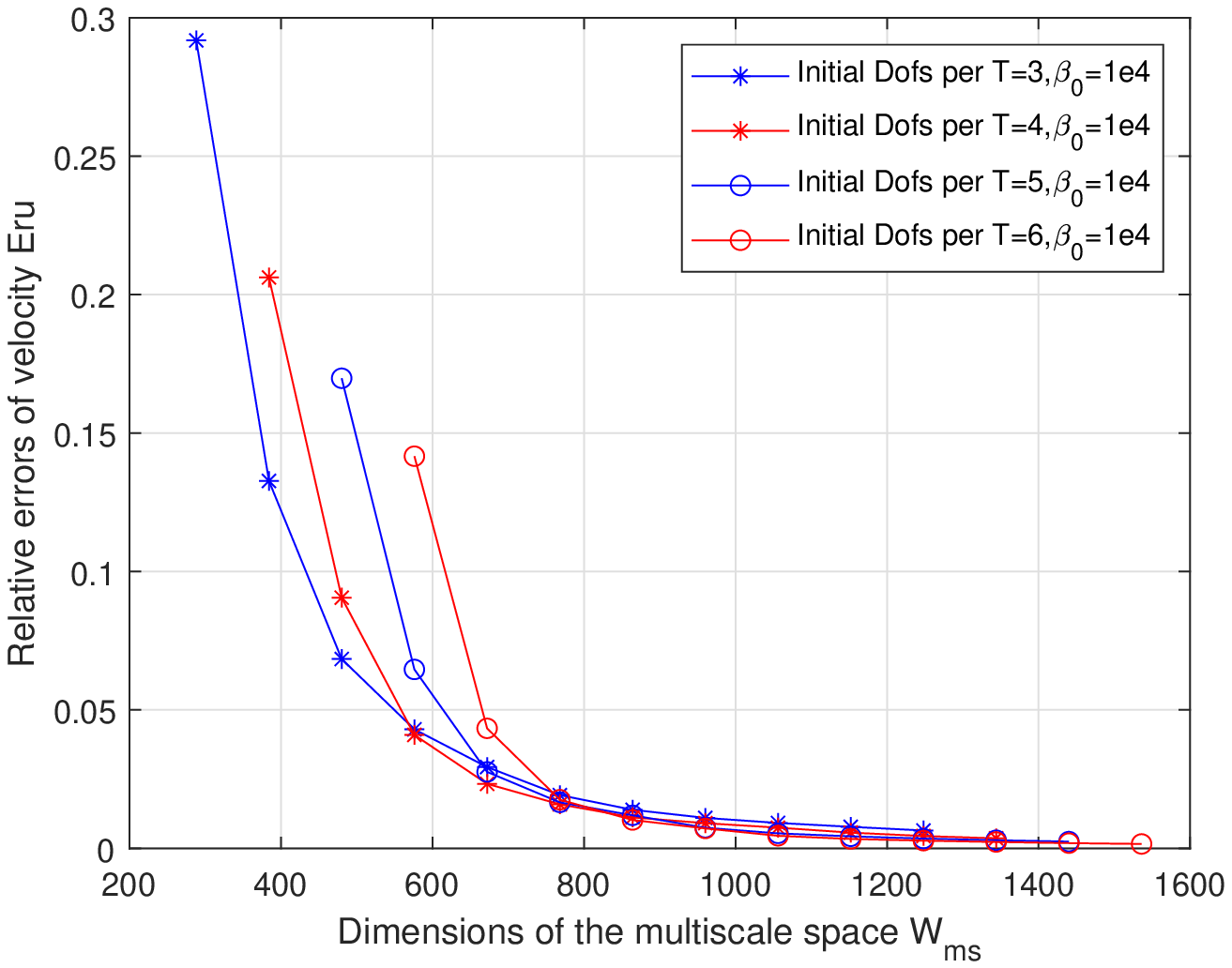}
	\end{minipage}
	\caption{(Uniform online enrichment) Convergence comparison with different choices of the number of initial basis functions. Top left: $\beta_0=10$. Top right: $\beta_0=100$. Bottom left: $\beta_0=1e3$. Bottom right: $\beta_0=1e4$.}\label{fig_velocityerr_2}
\end{figure}

Under the circumstance in Remark 2 where the mass matrix for velocity $A^m_H$ is not updated in the calculation of online basis functions using (\ref{eqn_online_1_remark2})-(\ref{eqn_online_2_remark2}) and multiscale solutions using (\ref{eqn_OnlinePicard_1_remark2})-(\ref{eqn_OnlinePicard_2_remark2}) for each enrichment iteration, in the same way, we compare the performance of the uniform online computation (\ref{eqn_online_1_remark2})-(\ref{eqn_OnlinePicard_2_remark2}) with the parameter $\beta_0$ taken different values for $3$, $4$, $5$ and $6$ initial basis functions per coarse element, respectively, the results are shown in Figure \ref{fig_velocityerr_1_rem2}. The comparisons with respect to different initial basis functions for $\beta_0$ taken different values : $\beta_0=10, 100, 10000$ and $10000$, are shown in Figure \ref{fig_velocityerr_2_rem2}. We find that in each case, the relative error decreases to a constant and no longer reduce after several levels of online enrichment, which is caused by the difference between $|\mathbf{u}_f|$ and $|\mathbf{u}_{\textrm{off}}|$ relating to the second term $(\beta\rho|\mathbf{u}_{\textrm{off}}|\cdot,\cdot)_Q$ in the left-hand side of (\ref{eqn_online_1_remark2}) and (\ref{eqn_OnlinePicard_1_remark2}). In other words, if the we know the fine-grid solution $\mathbf{u}_f$ and use the term $(\beta\rho|\mathbf{u}_f|\cdot,\cdot)_Q$ to replace the term $(\beta\rho|\mathbf{u}_{\textrm{off}}|\cdot,\cdot)_Q$ in (\ref{eqn_online_1_remark2}) and (\ref{eqn_OnlinePicard_1_remark2}), respectively, then the online multiscale solution will converge to the fine-grid solution as we increase the number of iterations of online enrichment. Anyway, totally speaking, these constant relative errors are acceptable. Furthermore, it is observed form Figure \ref{fig_velocityerr_1_rem2} that with the increase of the parameter $\beta_0$, the relative errors also increase, but the increase rates of the relative errors slow down when $\beta_0$ becomes large; from Figure \ref{fig_velocityerr_2_rem2}, it indicates that the more number of initial basis functions per coarse element are used, the smaller of the constant relative errors can be derived and the online basis functions behave better than the offline basis functions.

\begin{figure}
	\centering
	\mbox{\hspace{0.00cm}}
	\begin{minipage}[b]{0.40\textwidth}
		\centering
		\includegraphics[scale=0.48]{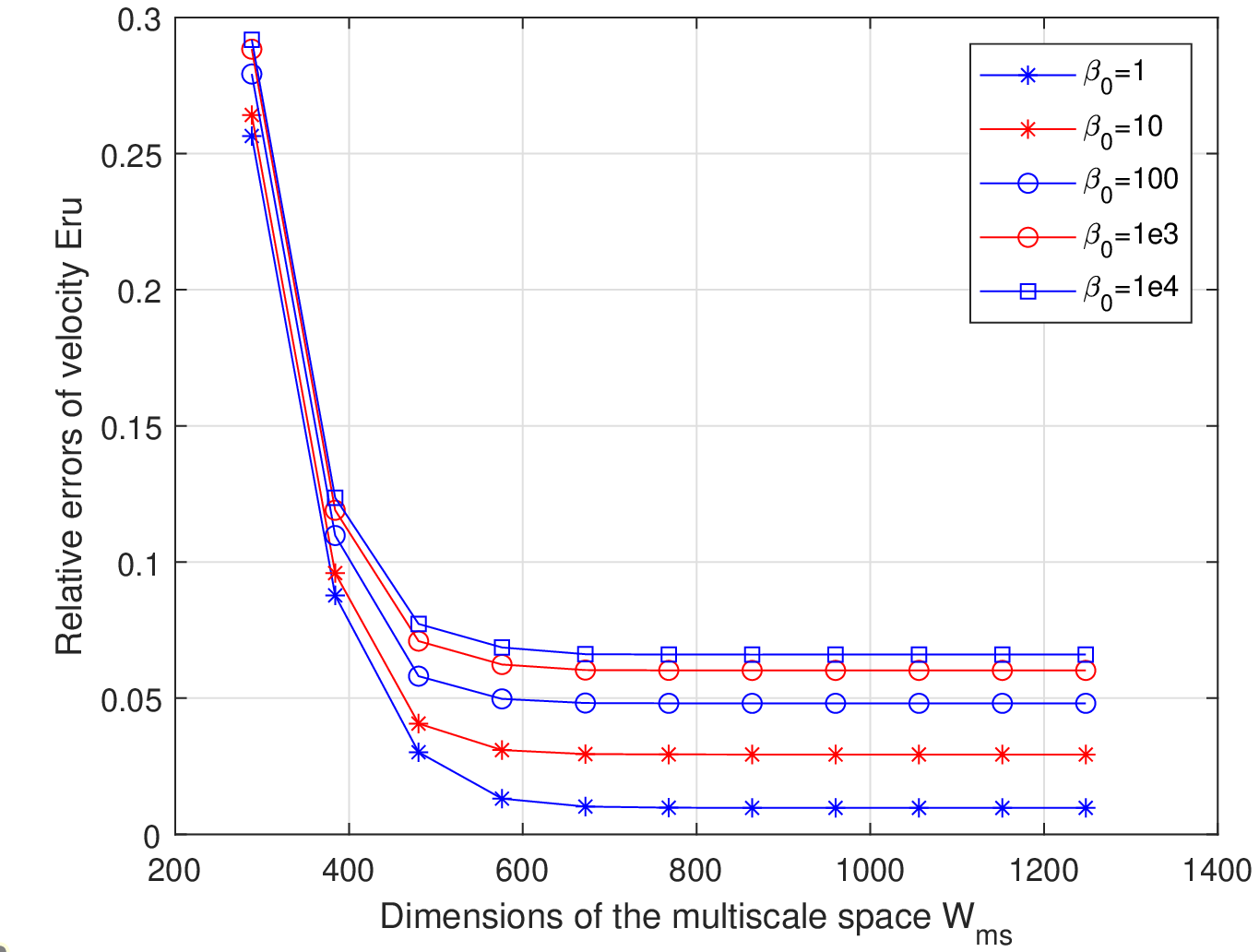}
	\end{minipage}
	\mbox{\hspace{0.00cm}}
	\begin{minipage}[b]{0.40\textwidth}
		\centering
		\includegraphics[scale=0.48]{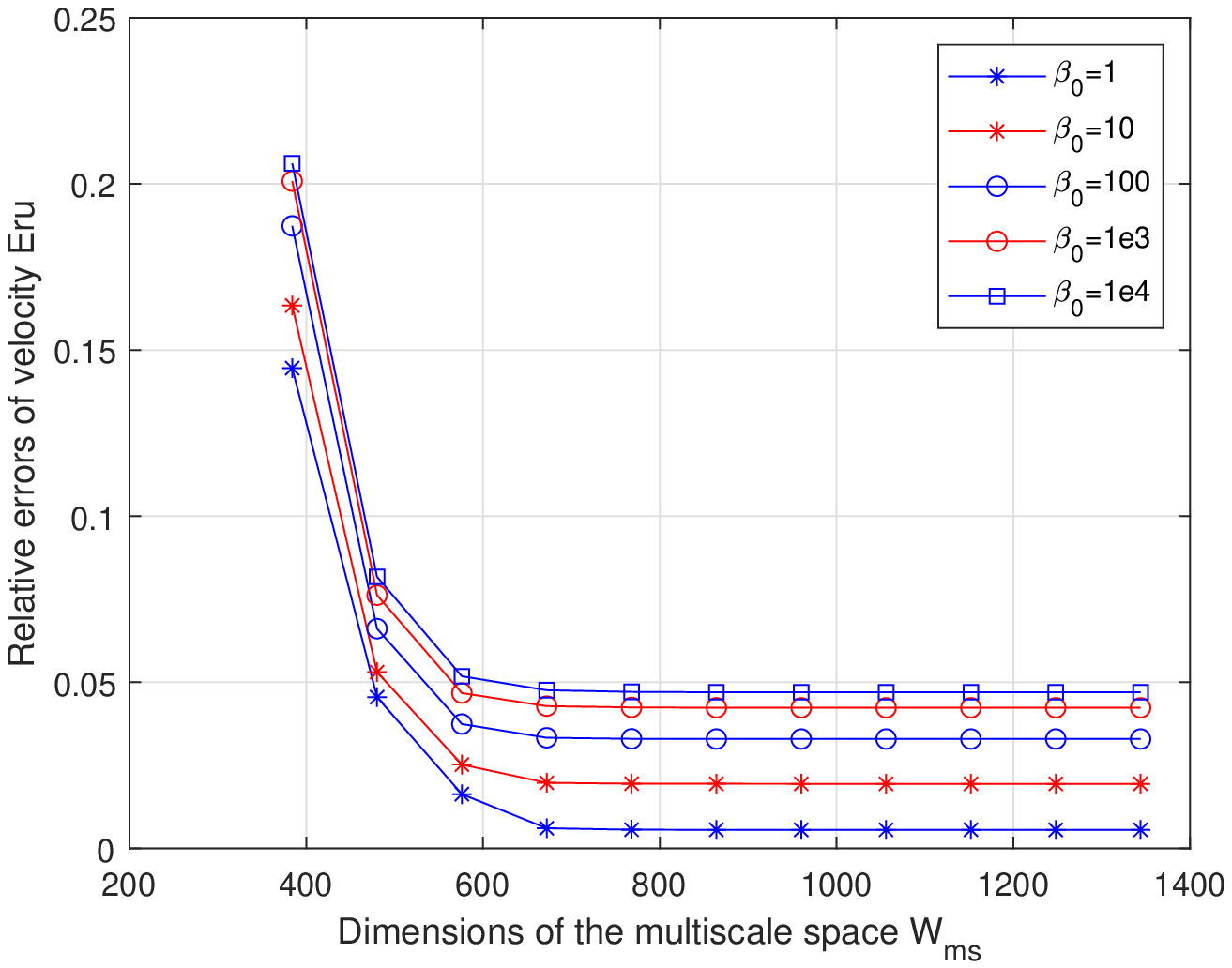}
	\end{minipage}
	\\
	\mbox{\hspace{0.00cm}}
	\begin{minipage}[b]{0.40\textwidth}
		\centering
		\includegraphics[scale=0.48]{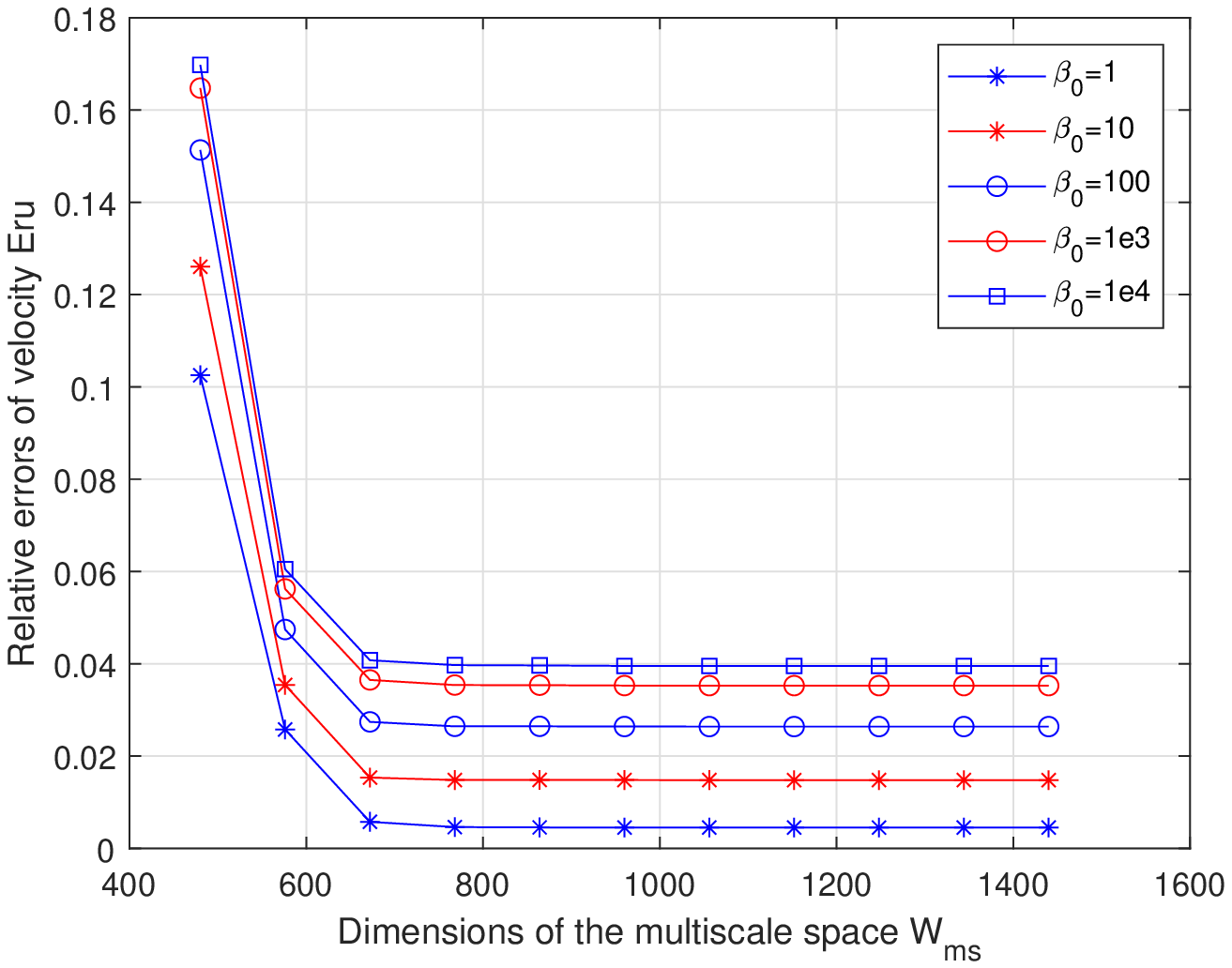}
	\end{minipage}
	\mbox{\hspace{0.00cm}}
	\begin{minipage}[b]{0.40\textwidth}
		\centering
		\includegraphics[scale=0.48]{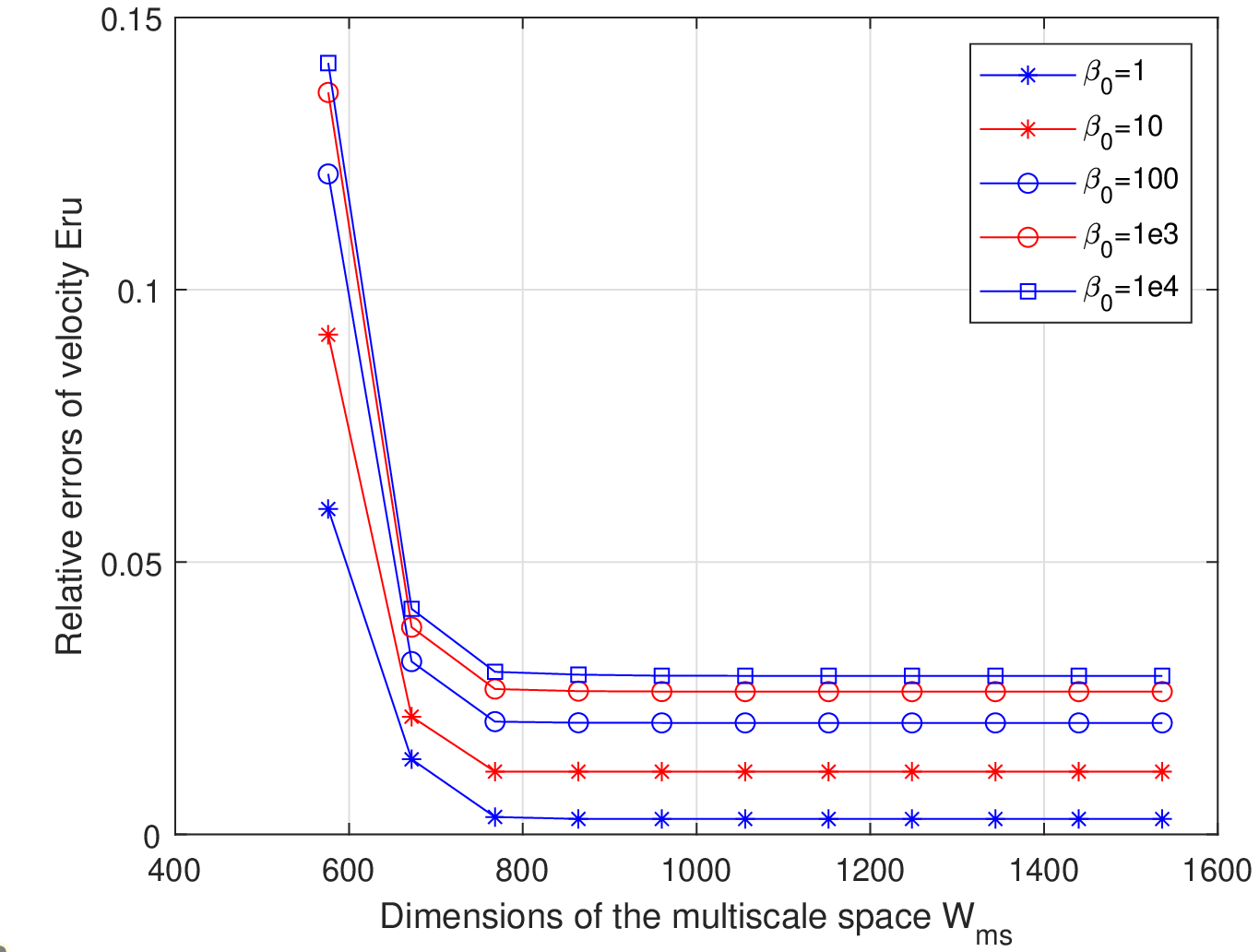}
	\end{minipage}
	\caption{(Uniform online enrichment, Remark 2) Convergence comparison with $\beta_0$ taking different values. Top left: 3 initial basis functions. Top right: 4 initial basis functions. Bottom left: 5 initial basis functions. Bottom right: 6 initial basis functions.}\label{fig_velocityerr_1_rem2}
\end{figure}
\begin{figure}
	\centering
	\mbox{\hspace{0.00cm}}
	\begin{minipage}[b]{0.40\textwidth}
		\centering
		\includegraphics[scale=0.48]{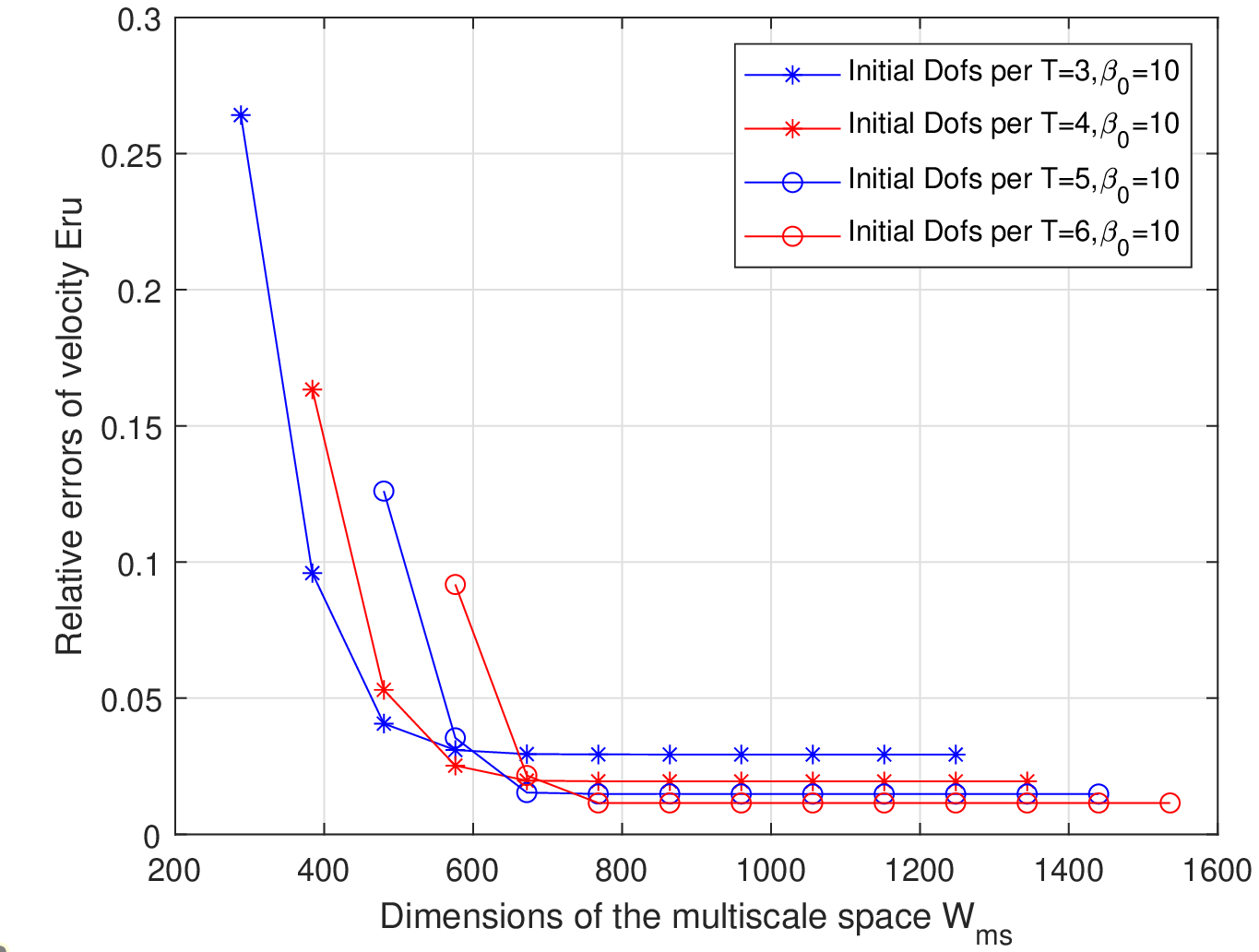}
	\end{minipage}
	\mbox{\hspace{0.00cm}}
	\begin{minipage}[b]{0.40\textwidth}
		\centering
		\includegraphics[scale=0.48]{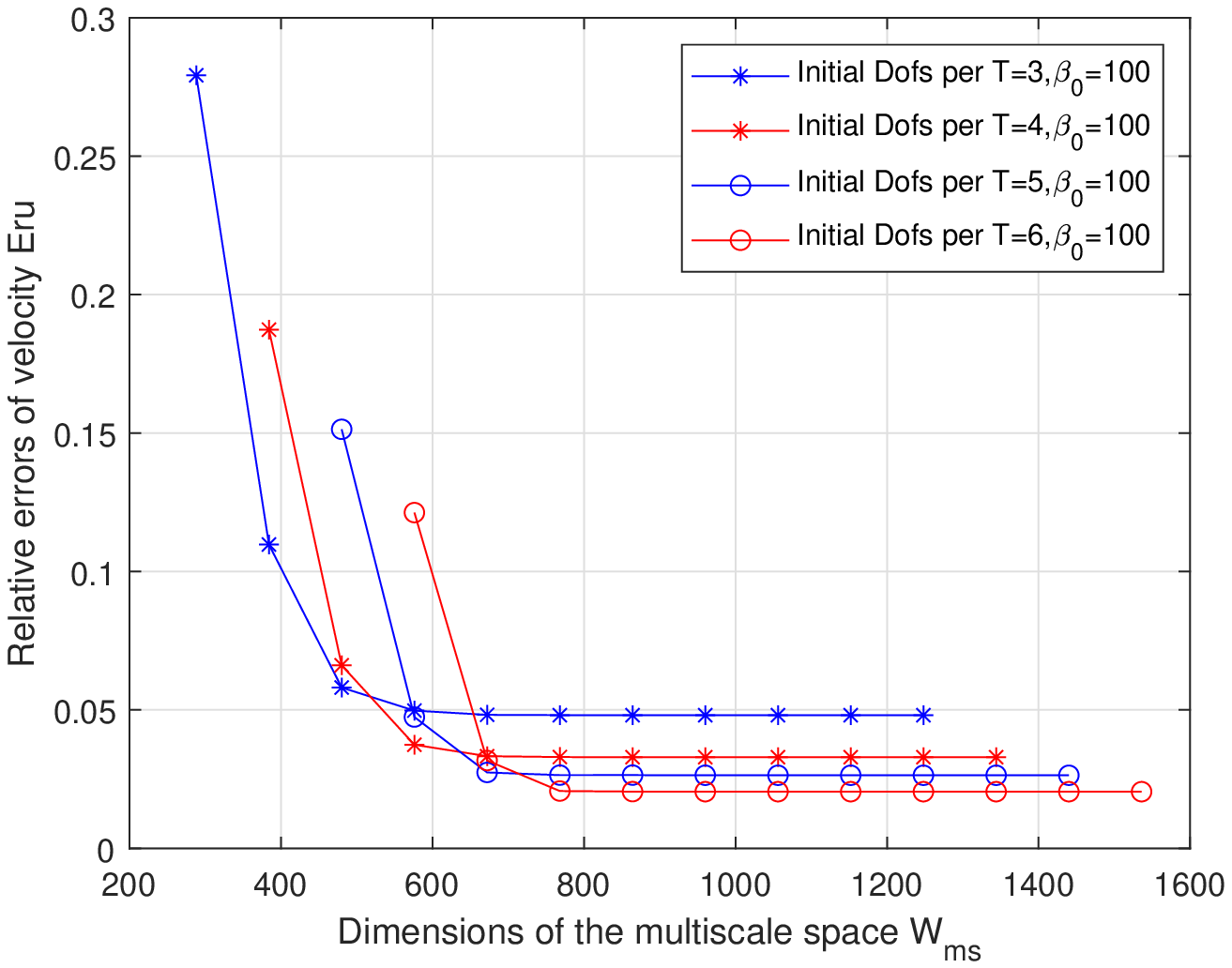}
	\end{minipage}
	\\
	\mbox{\hspace{0.00cm}}
	\begin{minipage}[b]{0.40\textwidth}
		\centering
		\includegraphics[scale=0.48]{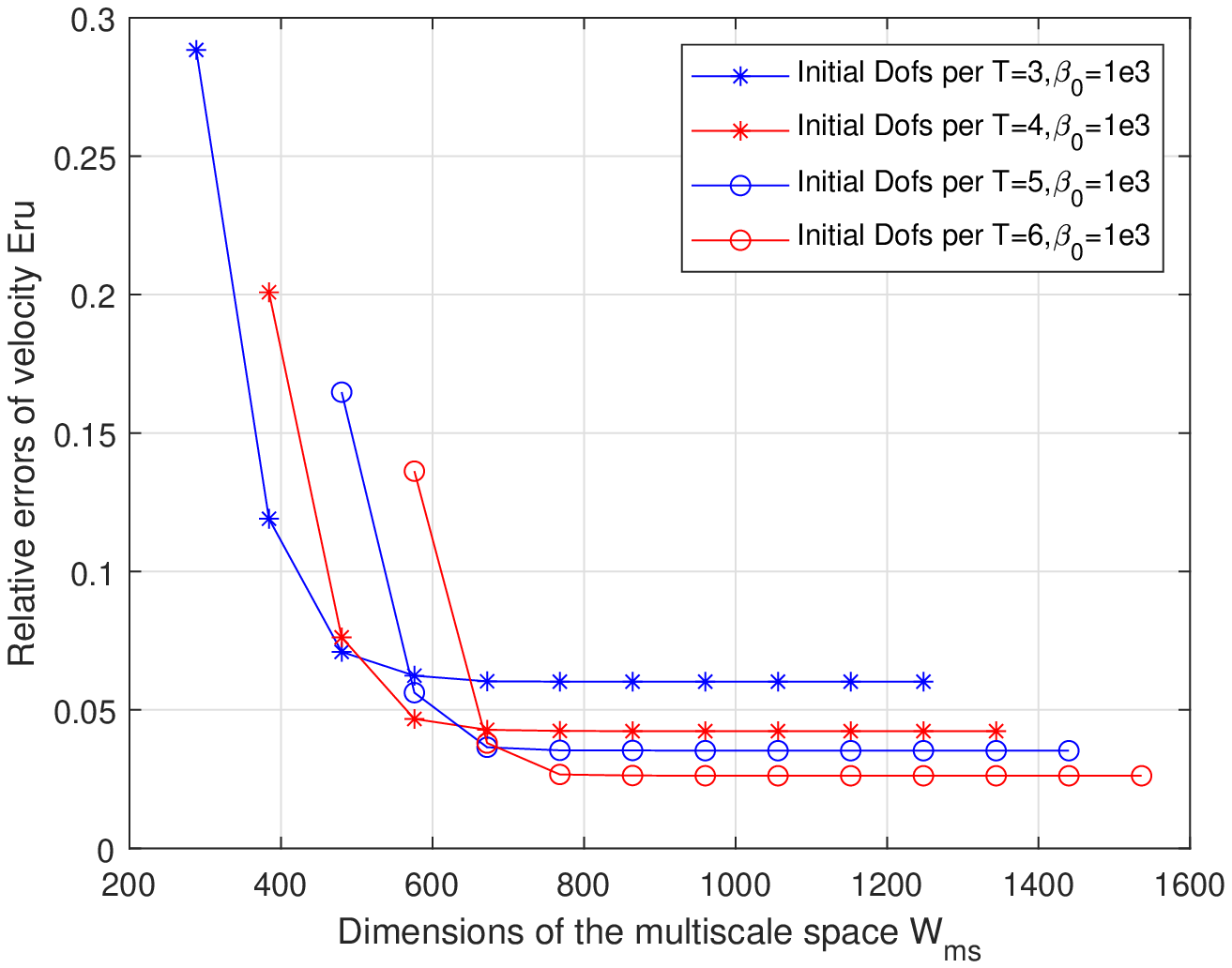}
	\end{minipage}
	\mbox{\hspace{0.00cm}}
	\begin{minipage}[b]{0.40\textwidth}
		\centering
		\includegraphics[scale=0.48]{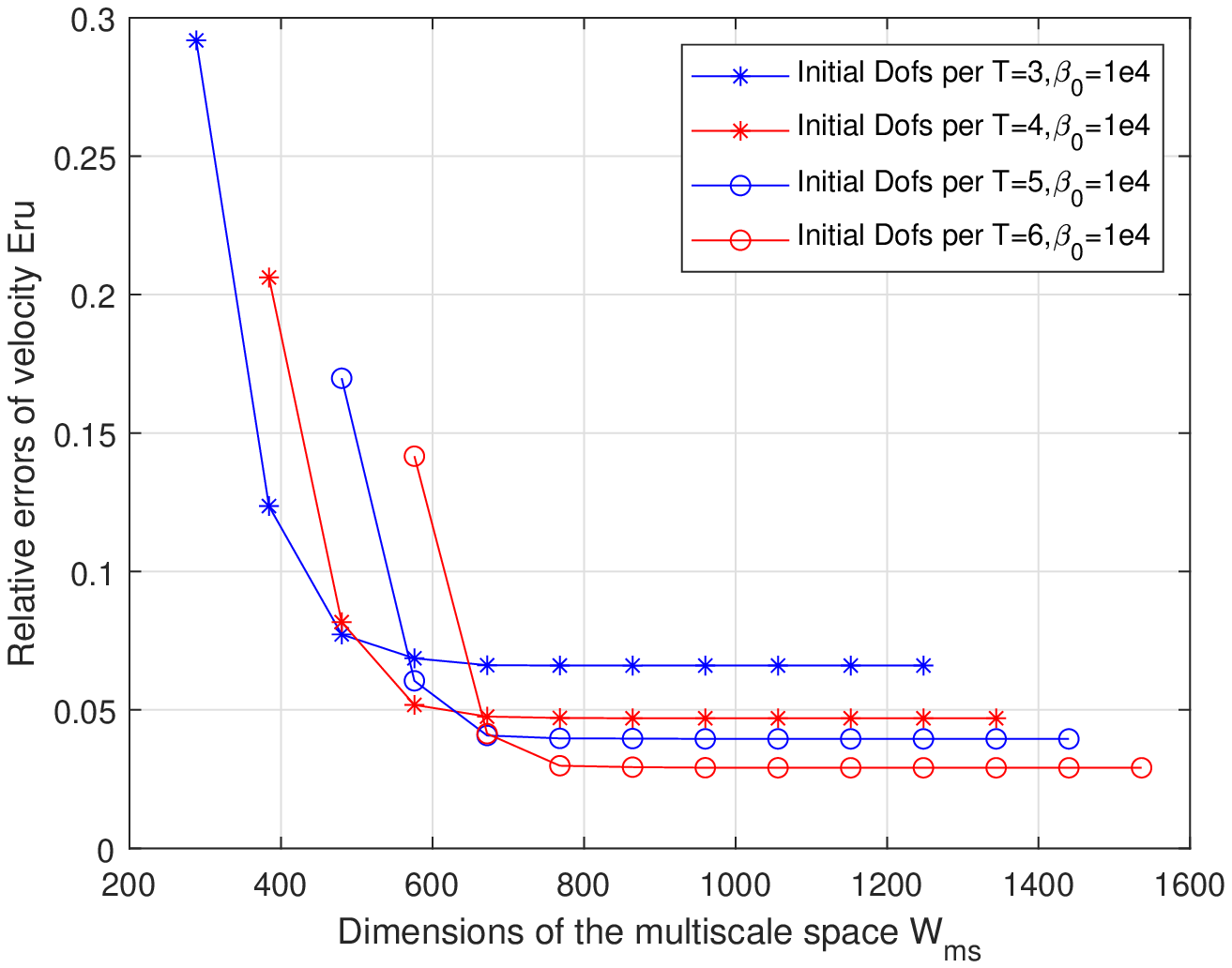}
	\end{minipage}
	\caption{(Uniform online enrichment, Remark 2) Convergence comparison for different choices of the number of initial basis. Top left: $\beta_0=1$. Top right: $\beta_0=10$. Middle left: $\beta_0=100$. Middle right: $\beta_0=1e3$. Bottom: $\beta_0 = 1e4$.}\label{fig_velocityerr_2_rem2}
\end{figure}
\subsubsection{Adaptive online enrichment}
In the following, we carry out the online multiscale space enrichment adaptively by adding online basis functions on coarse elements where the residuals are large. For the online enrichment in level $m+1$, similar with (\ref{eqn_residual_offline}), we define the online residual $R^m_i$ for each coarse element $T_i$ as
\begin{equation}\label{eqn_residual_online}
	R^m_i = \int_{T_i}{|f-\nabla\cdot\mathbf{u}^m_{\textrm{ms}}|^2}\textrm{d}x,
\end{equation}
we arrange these online residuals in decreasing order, $R^m_1 \ge R^m_2 \ge \cdots \ge R^m_{N_T}$, then we choose the coarse elements where online basis functions need to be added by choosing the smallest integer $N_{\textrm{add}}$, such that the cumulative online residuals of the selected coarse elements is the $\xi$ fraction of the total online residuals of the multiscale solution in level $m$, that is
\begin{equation}\label{eqn_criterion_online}
	\sum^{N_{\textrm{add}}}_{i=1}R^m_i \ge \xi\sum^{N_T}_{i=1}R^m_i,
\end{equation}
where $1>\xi>0$ is a real number to be chosen. In the following tests, we take $\xi=3/4$. As done in numerical tests of the uniform online enrichment, in Figure \ref{fig_velocityerr_3}, we plot the convergence history of the adaptive online computation with the parameter $\beta_0$ taken different values regarding $3$, $4$, $5$ and $6$ initial basis functions per coarse element, respectively; and in Figure \ref{fig_velocityerr_4}, we show the convergence history of the adaptive online computation with different number of initial basis functions for $\beta_0$ taken four different values : $\beta_0=10, 100, 10000$ and $10000$, respectively. Compare with the results in Figure \ref{fig_velocityerr_1} and \ref{fig_velocityerr_2} of the uniform online enrichment, we can observe that the criterion (\ref{eqn_criterion_online}) is effective, the number of basis functions are reduced to obtain the same relative error as the uniform online enrichment.

\begin{figure}
	\centering
	\mbox{\hspace{0.00cm}}
	\begin{minipage}[b]{0.40\textwidth}
		\centering
		\includegraphics[scale=0.48]{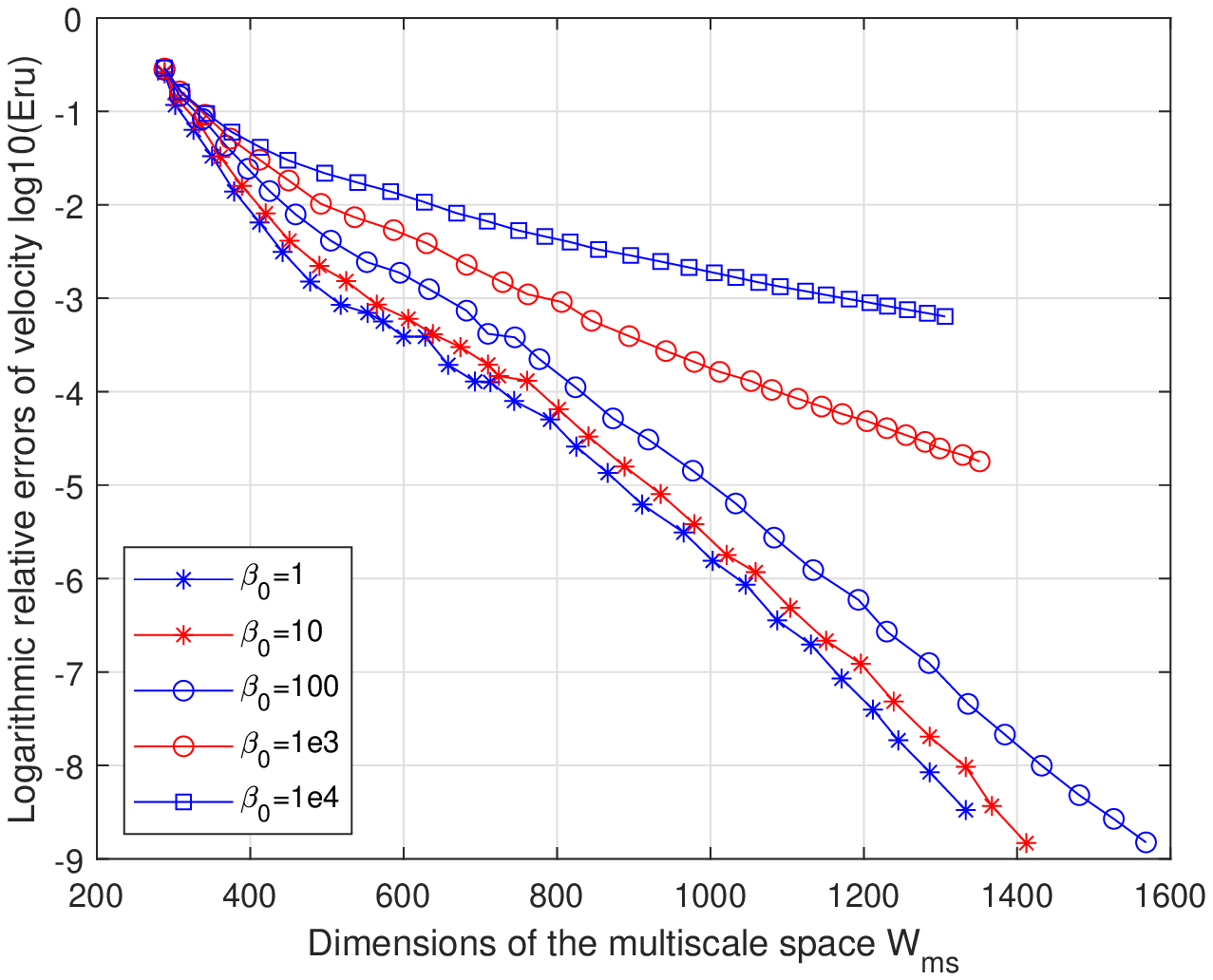}
	\end{minipage}
	\mbox{\hspace{0.00cm}}
	\begin{minipage}[b]{0.40\textwidth}
		\centering
		\includegraphics[scale=0.48]{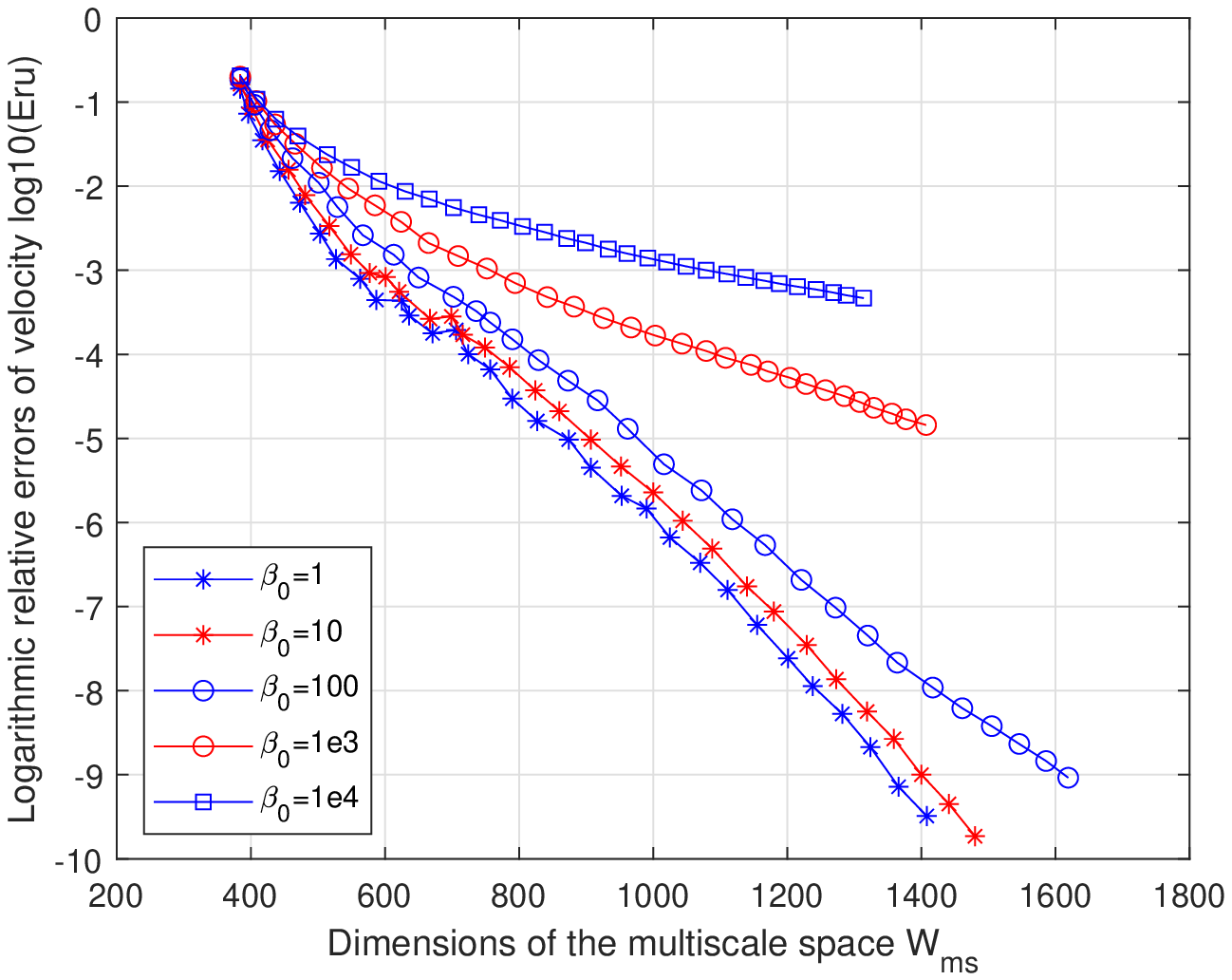}
	\end{minipage}
	\\
	\mbox{\hspace{0.00cm}}
	\begin{minipage}[b]{0.40\textwidth}
		\centering
		\includegraphics[scale=0.48]{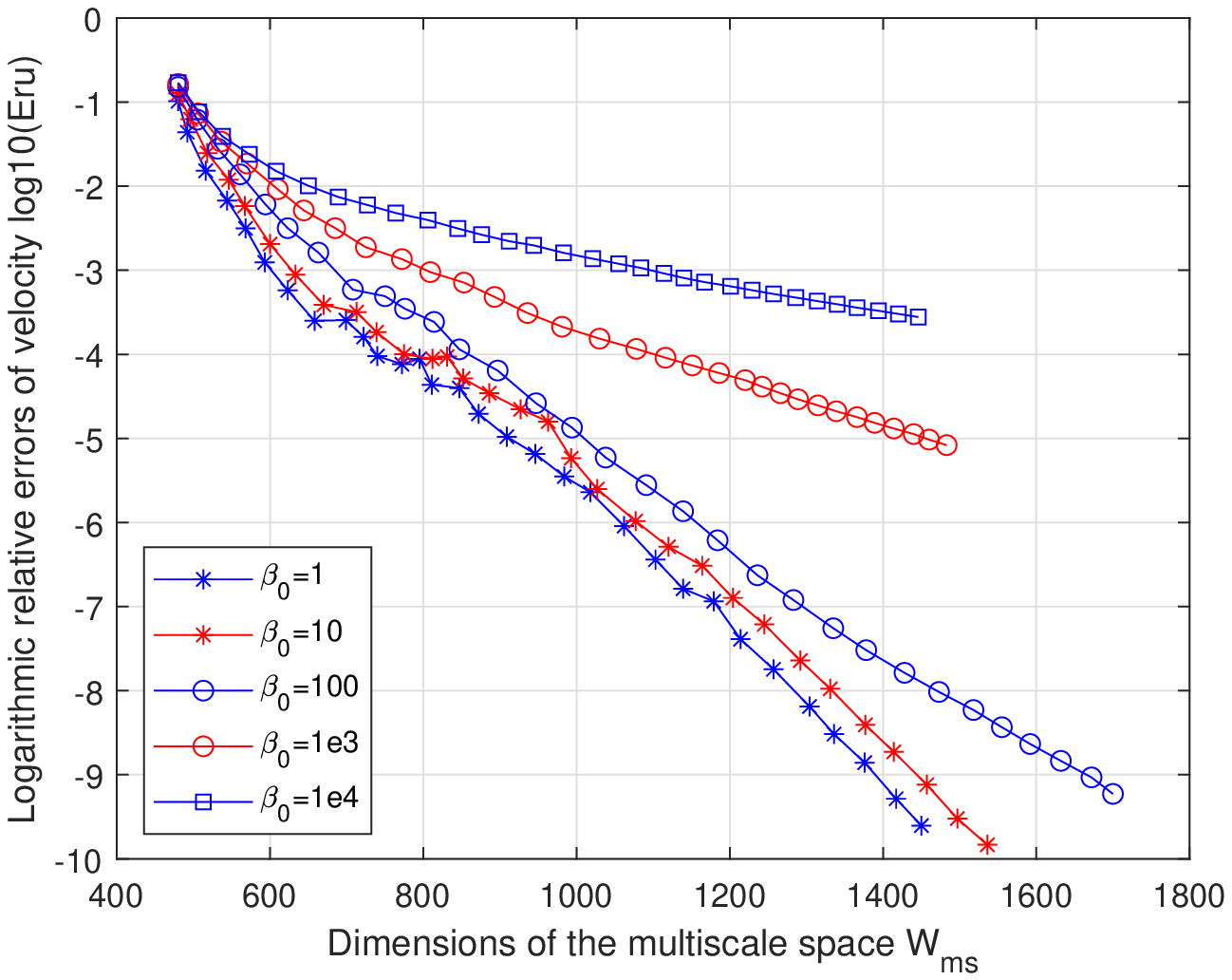}
	\end{minipage}
	\mbox{\hspace{0.00cm}}
	\begin{minipage}[b]{0.40\textwidth}
		\centering
		\includegraphics[scale=0.48]{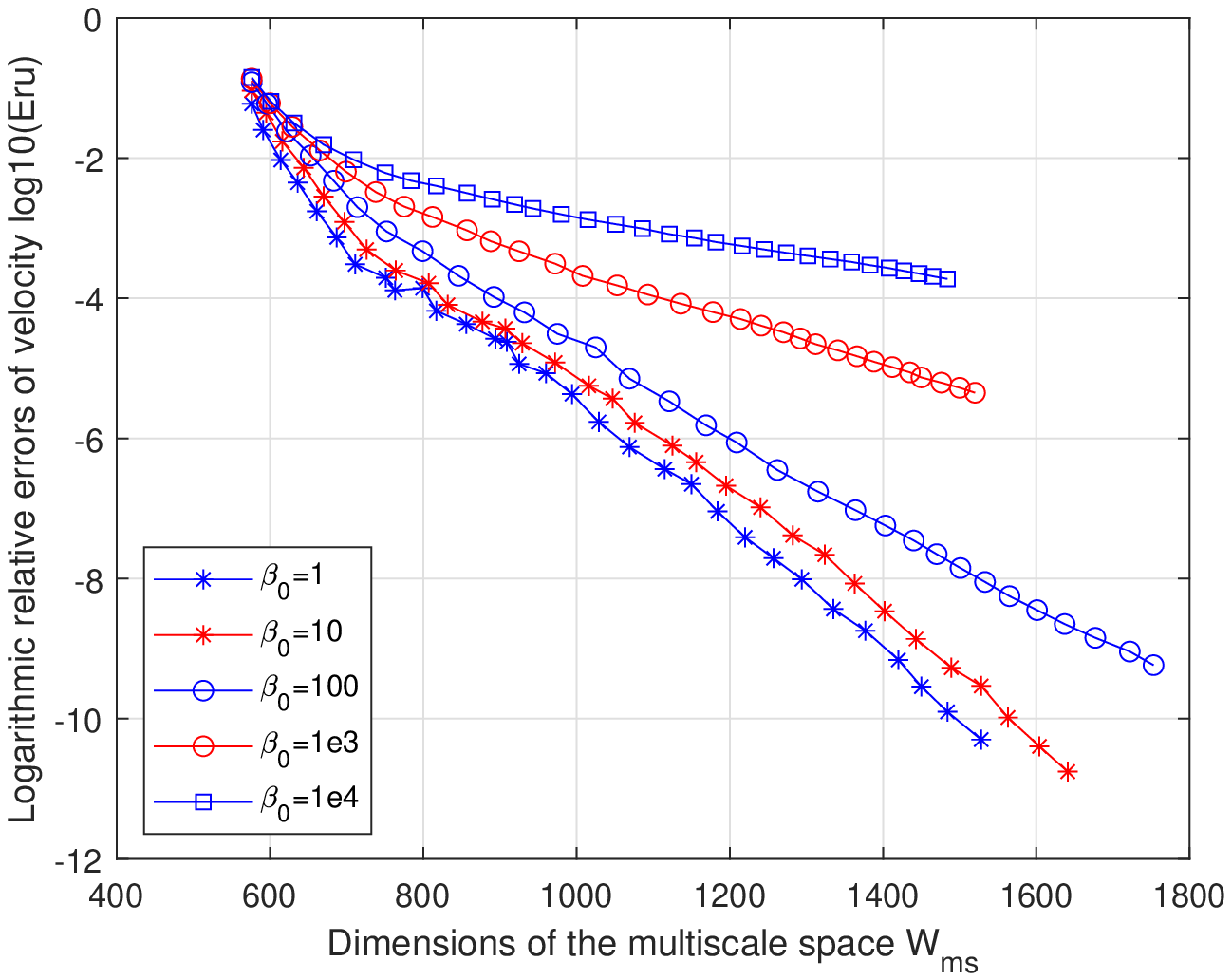}
	\end{minipage}
	\caption{(Adaptive online enrichment) Convergence comparison with $\beta_0$ taking different values. Top left: 3 initial basis functions. Top right: 4 initial basis functions. Bottom left: 5 initial basis functions. Bottom right: 6 initial basis functions.}\label{fig_velocityerr_3}
\end{figure}
\begin{figure}
	\centering
	\mbox{\hspace{0.00cm}}
	\begin{minipage}[b]{0.40\textwidth}
		\centering
		\includegraphics[scale=0.48]{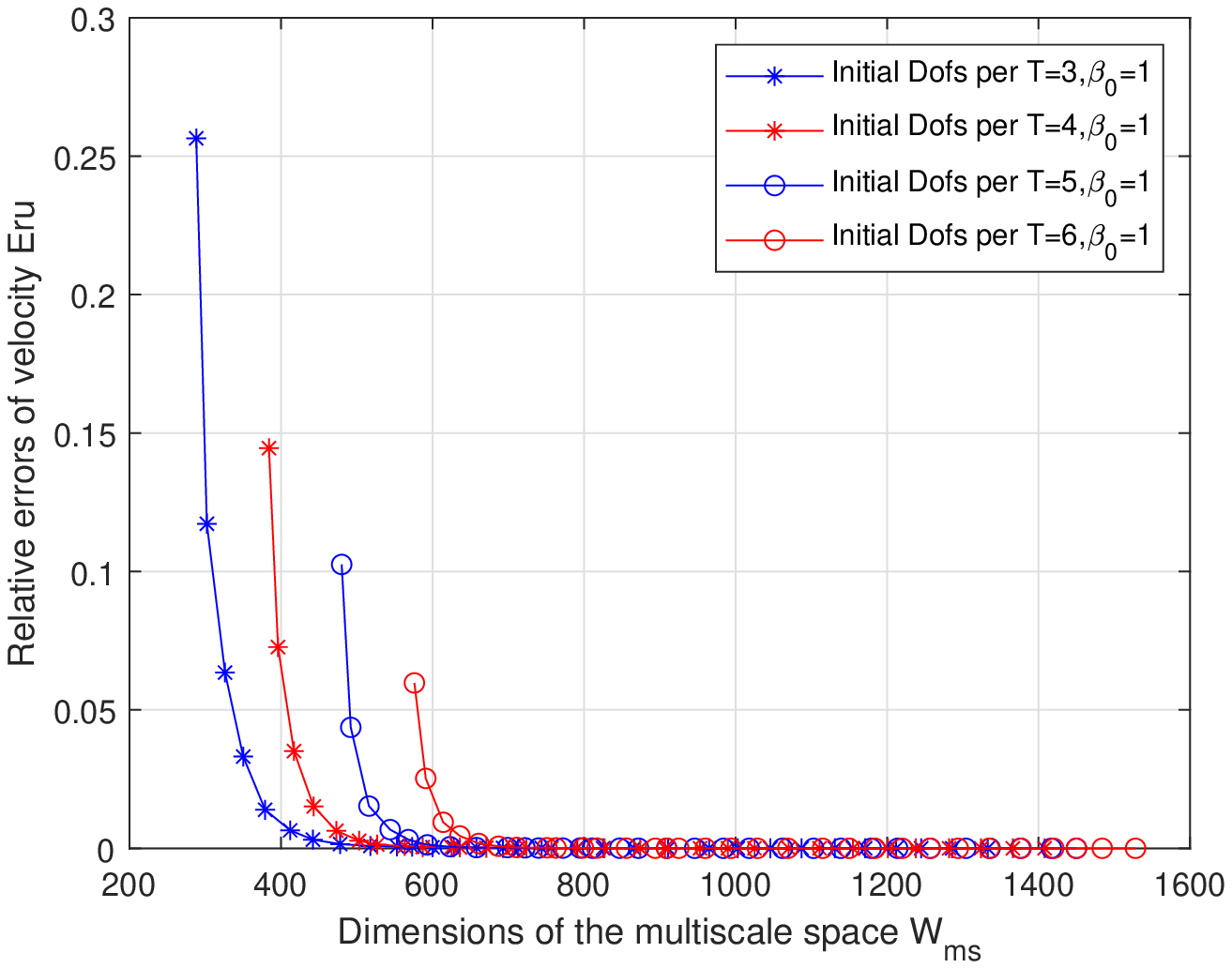}
	\end{minipage}
	\mbox{\hspace{0.00cm}}
	\begin{minipage}[b]{0.40\textwidth}
		\centering
		\includegraphics[scale=0.48]{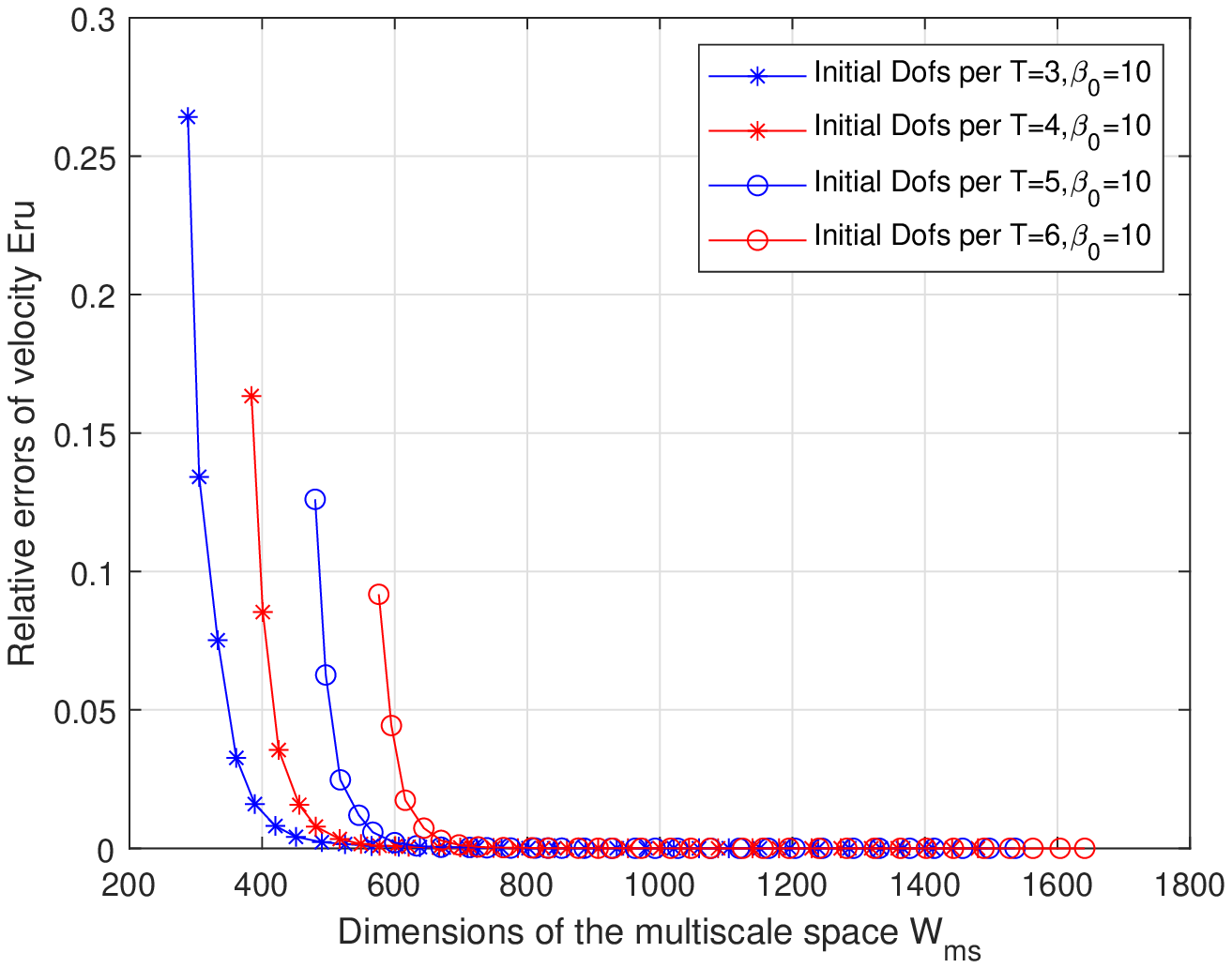}
	\end{minipage}
	\\
	\mbox{\hspace{0.00cm}}
	\begin{minipage}[b]{0.40\textwidth}
		\centering
		\includegraphics[scale=0.48]{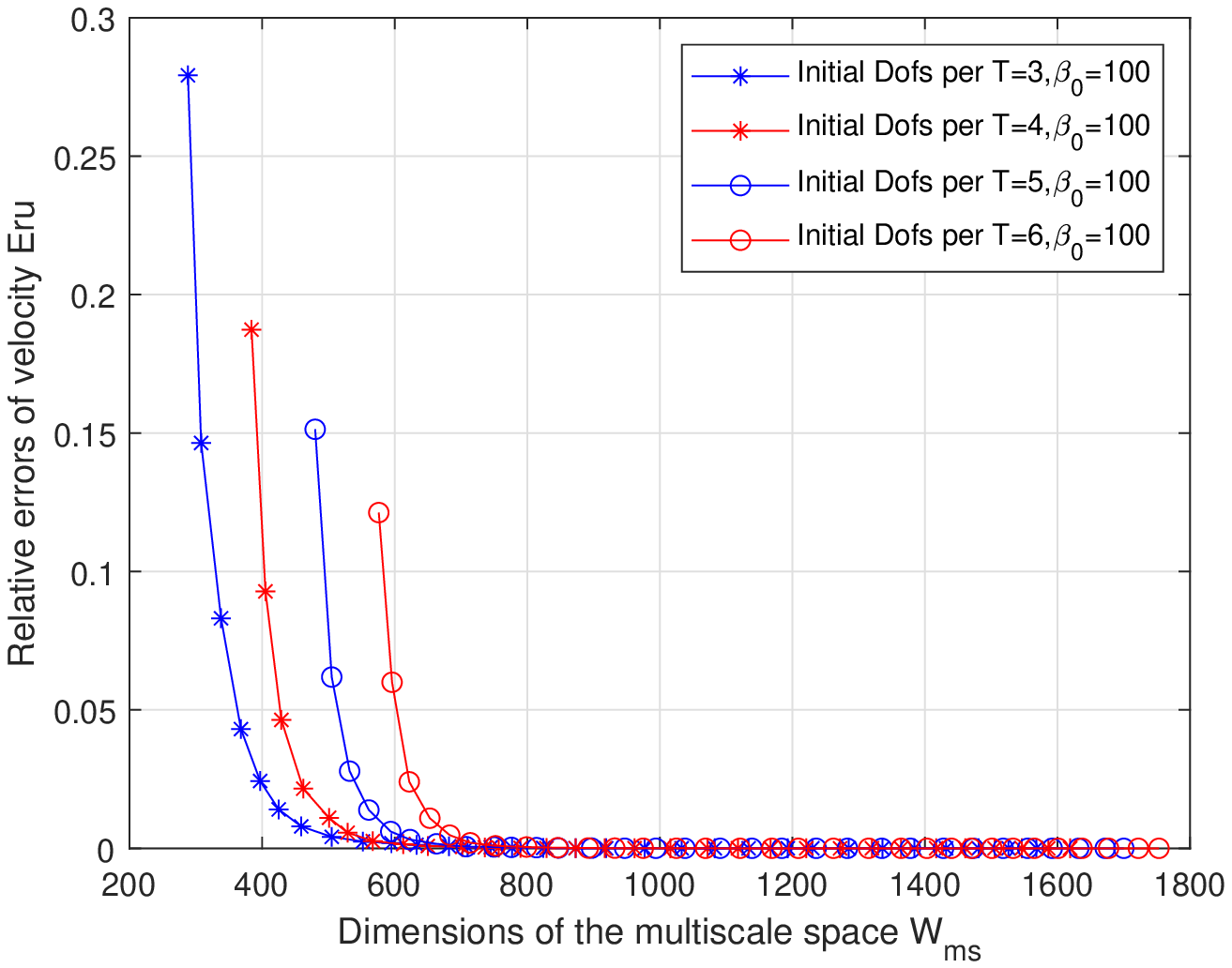}
	\end{minipage}
	\mbox{\hspace{0.00cm}}
	\begin{minipage}[b]{0.40\textwidth}
		\centering
		\includegraphics[scale=0.48]{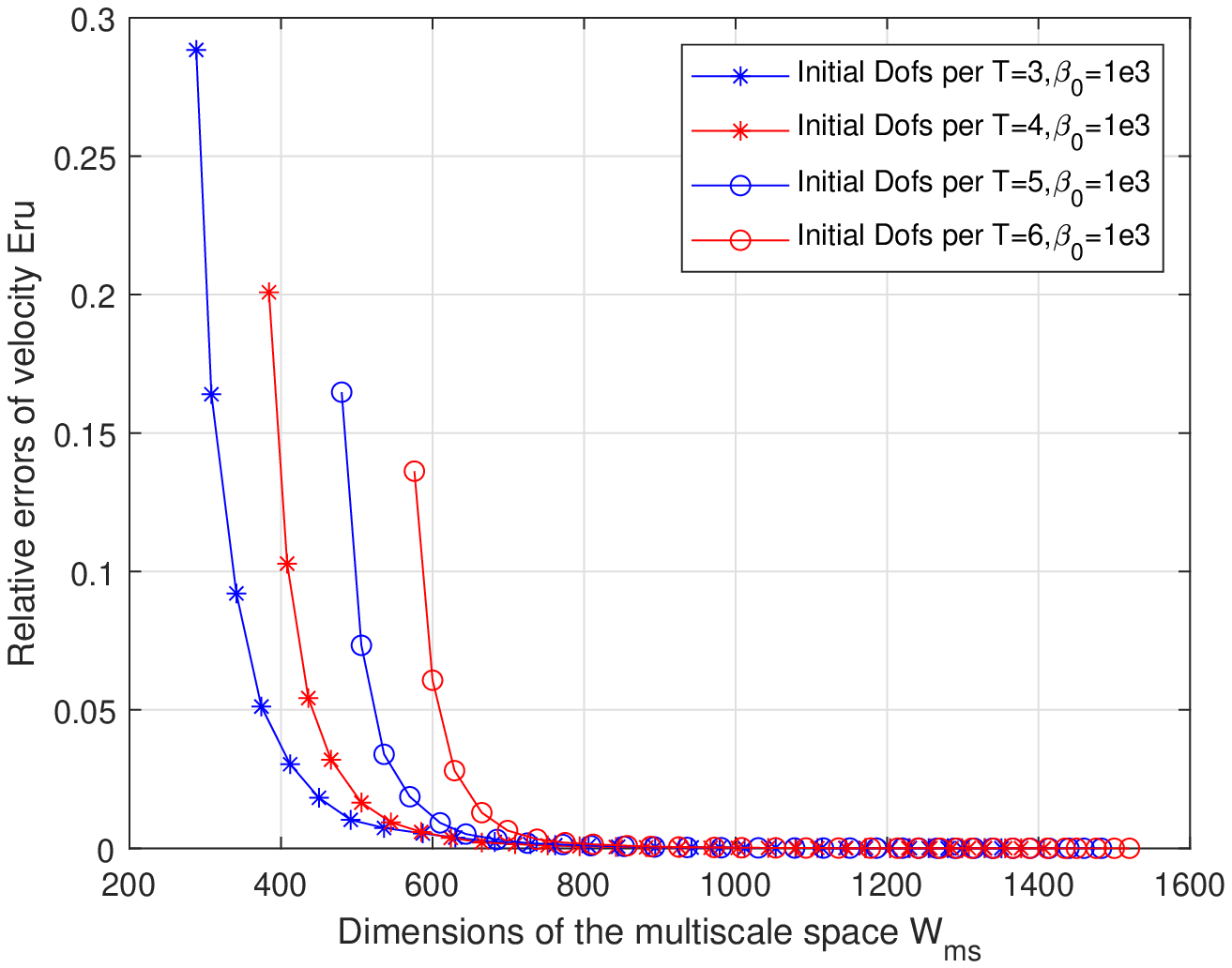}
	\end{minipage}
	\caption{(Adaptive online enrichment) Convergence comparison for different choices of the number of initial basis. Top left: $\beta_0=1$. Top right: $\beta_0=10$. Middle left: $\beta_0=100$. Middle right: $\beta_0=1e3$. Bottom: $\beta_0 = 1e4$.}\label{fig_velocityerr_4}
\end{figure}

We also test the situation in Remark 2 where the mass matrix for velocity $A^m_H$ is keep fixed in each iteration of the multiscale space enrichment.
The results of the comparison with the parameter $\beta_0$ taken different values for $3$, $4$, $5$ and $6$ initial basis functions per coarse element, respectively, are shown in Figure \ref{fig_velocityerr_3_rem2}; and the results of the comparison with respect to different initial basis functions per coarse element for $\beta_0$ taken four different values : $\beta_0=10, 100, 10000$ and $10000$, respectively, are shown in Figure \ref{fig_velocityerr_4_rem2}. As the results in Figure \ref{fig_velocityerr_1_rem2} and \ref{fig_velocityerr_2_rem2} of the uniform online enrichment, the relative errors also convergent to constant errors, and these constant errors in Figure \ref{fig_velocityerr_1_rem2} and \ref{fig_velocityerr_2_rem2} are practically coincide with the constant errors in \ref{fig_velocityerr_3_rem2} and \ref{fig_velocityerr_4_rem2}. In addition, compare with the results in Figure \ref{fig_velocityerr_1_rem2} and \ref{fig_velocityerr_2_rem2}, we observe once again that the adaptive online enrichment requires smaller number of basis functions than the uniform online enrichment to achieve the same accuracy of the multiscale solution.

\begin{figure}
	\centering
	\mbox{\hspace{0.00cm}}
	\begin{minipage}[b]{0.40\textwidth}
		\centering
		\includegraphics[scale=0.48]{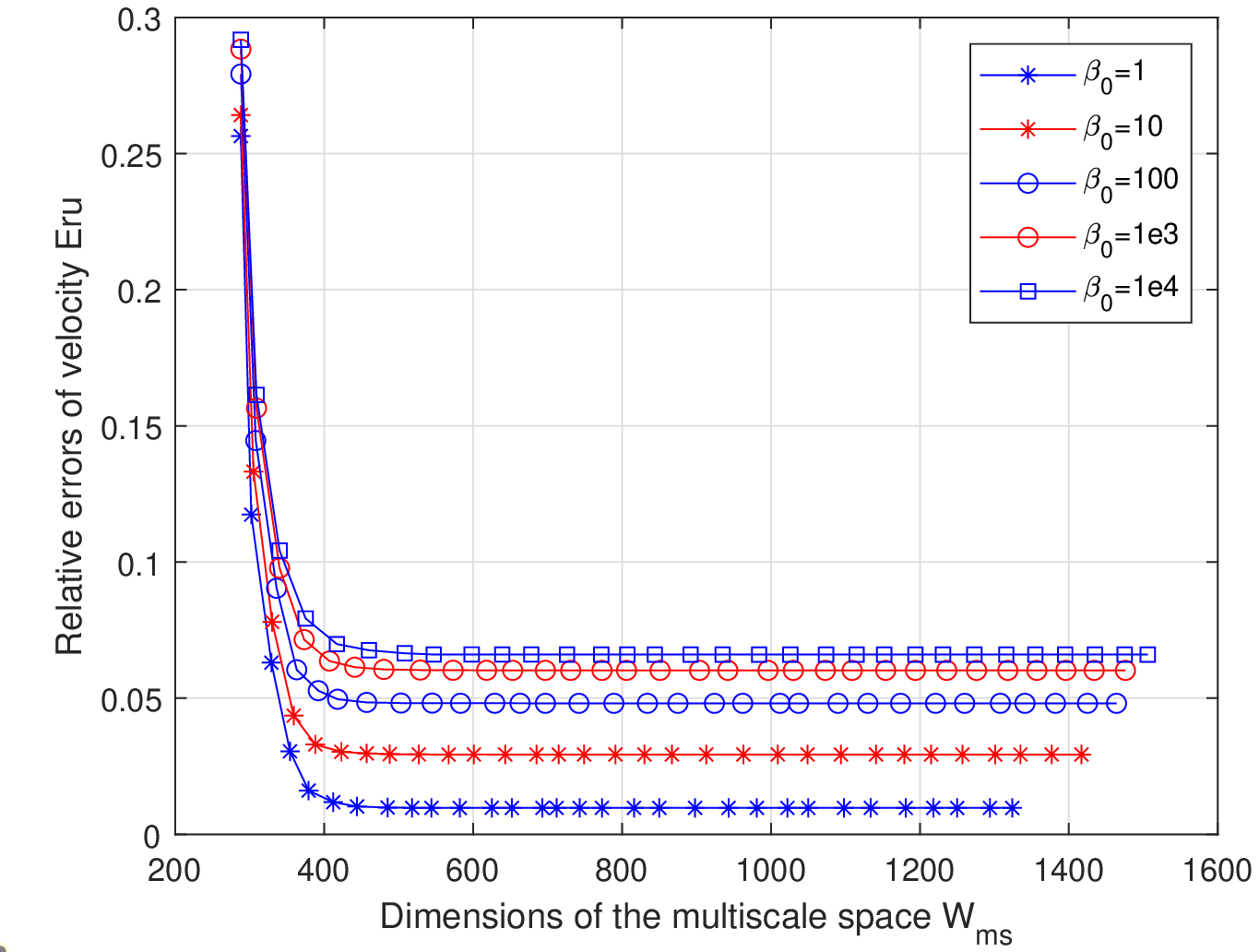}
	\end{minipage}
	\mbox{\hspace{0.00cm}}
	\begin{minipage}[b]{0.40\textwidth}
		\centering
		\includegraphics[scale=0.48]{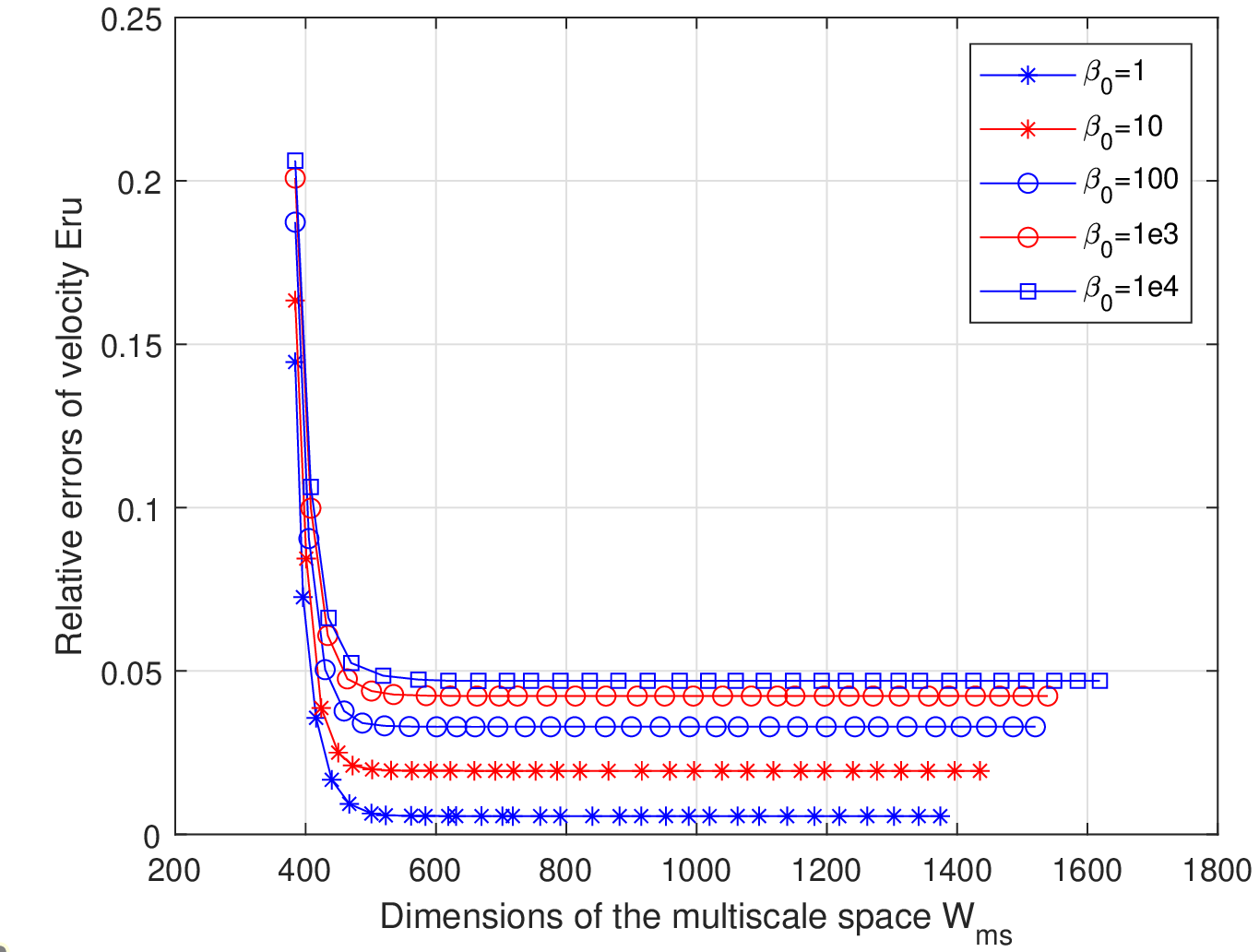}
	\end{minipage}
	\\
	\mbox{\hspace{0.00cm}}
	\begin{minipage}[b]{0.40\textwidth}
		\centering
		\includegraphics[scale=0.48]{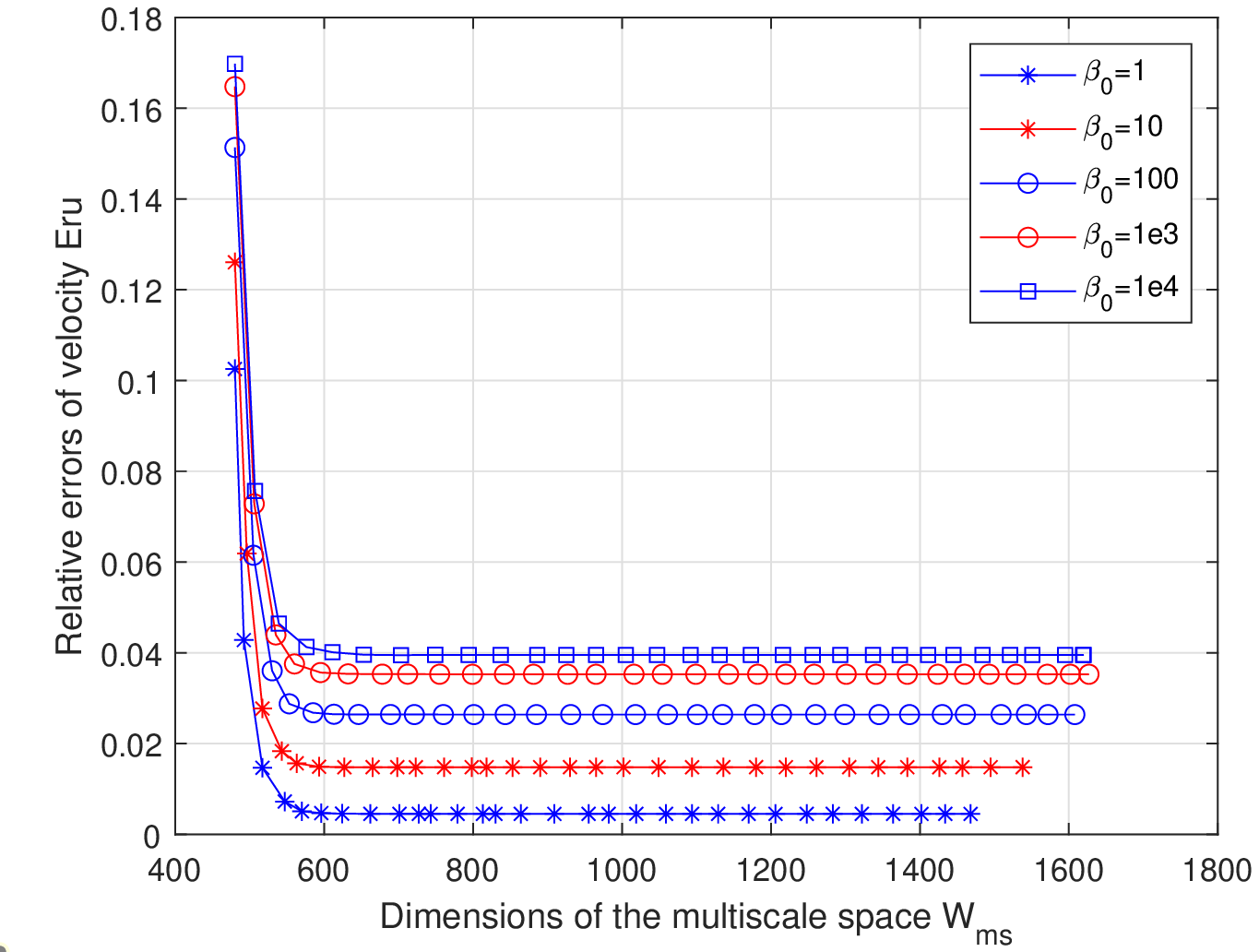}
	\end{minipage}
	\mbox{\hspace{0.00cm}}
	\begin{minipage}[b]{0.40\textwidth}
		\centering
		\includegraphics[scale=0.48]{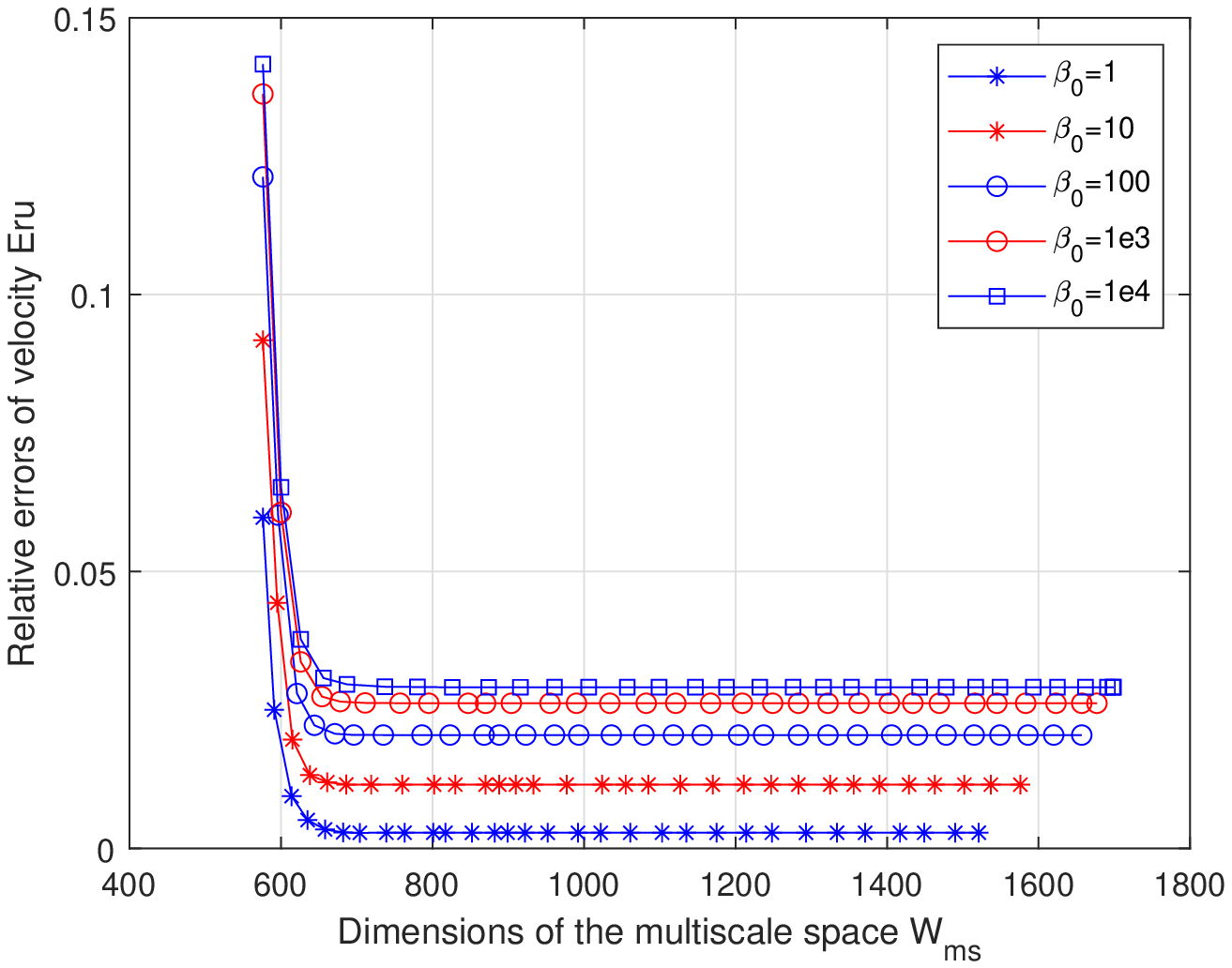}
	\end{minipage}
	\caption{Adaptive online enrichment (Remark 2): Convergence comparison for $\beta_0$ taking different values. Top left: initial dofs per T=3. Top right: initial dofs per T=4. Bottom left: initial dofs per T=5. Bottom right: initial dofs per T = 6.}\label{fig_velocityerr_3_rem2}
\end{figure}
\begin{figure}
	\centering
	\mbox{\hspace{0.00cm}}
	\begin{minipage}[b]{0.40\textwidth}
		\centering
		\includegraphics[scale=0.48]{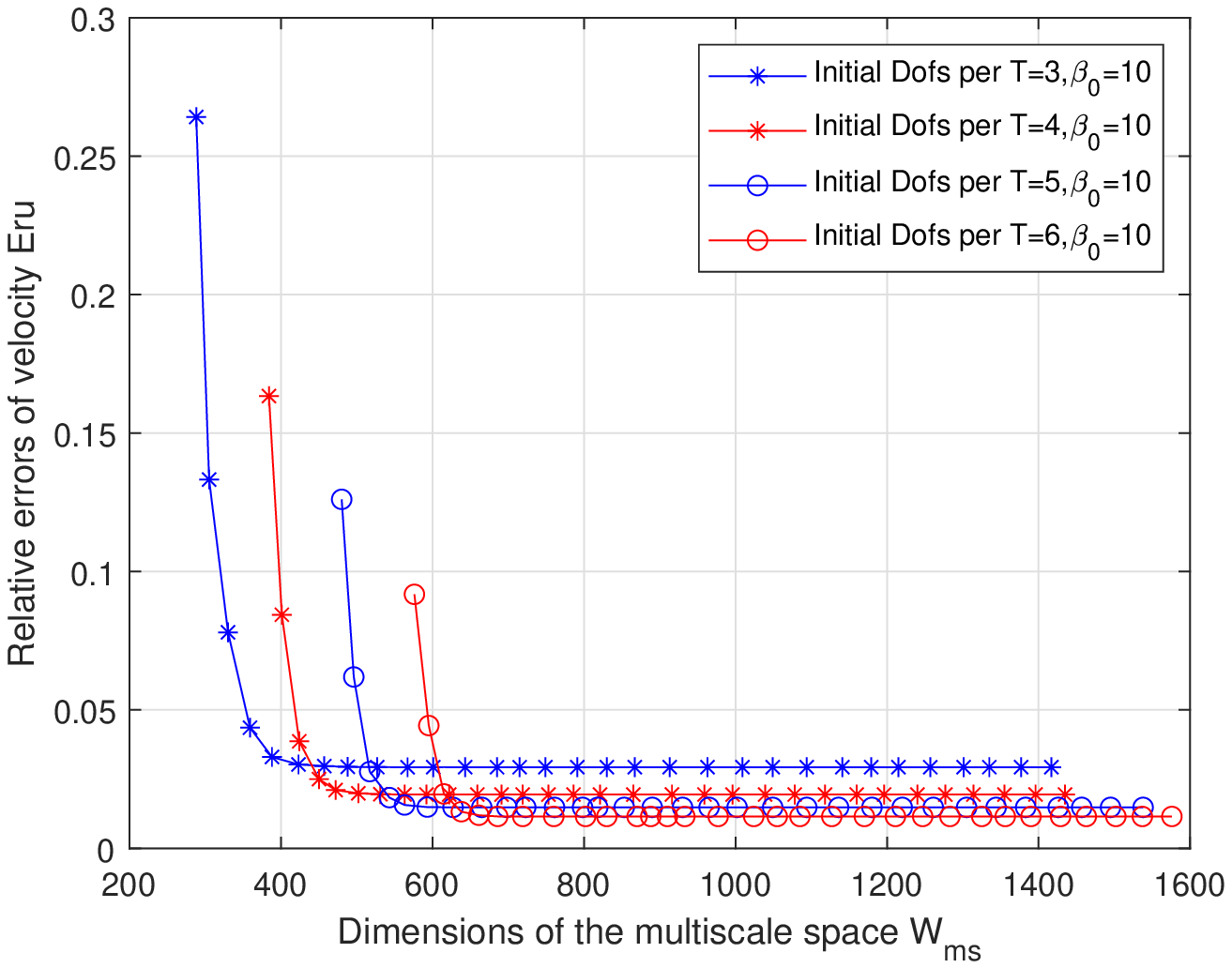}
	\end{minipage}
	\mbox{\hspace{0.00cm}}
	\begin{minipage}[b]{0.40\textwidth}
		\centering
		\includegraphics[scale=0.48]{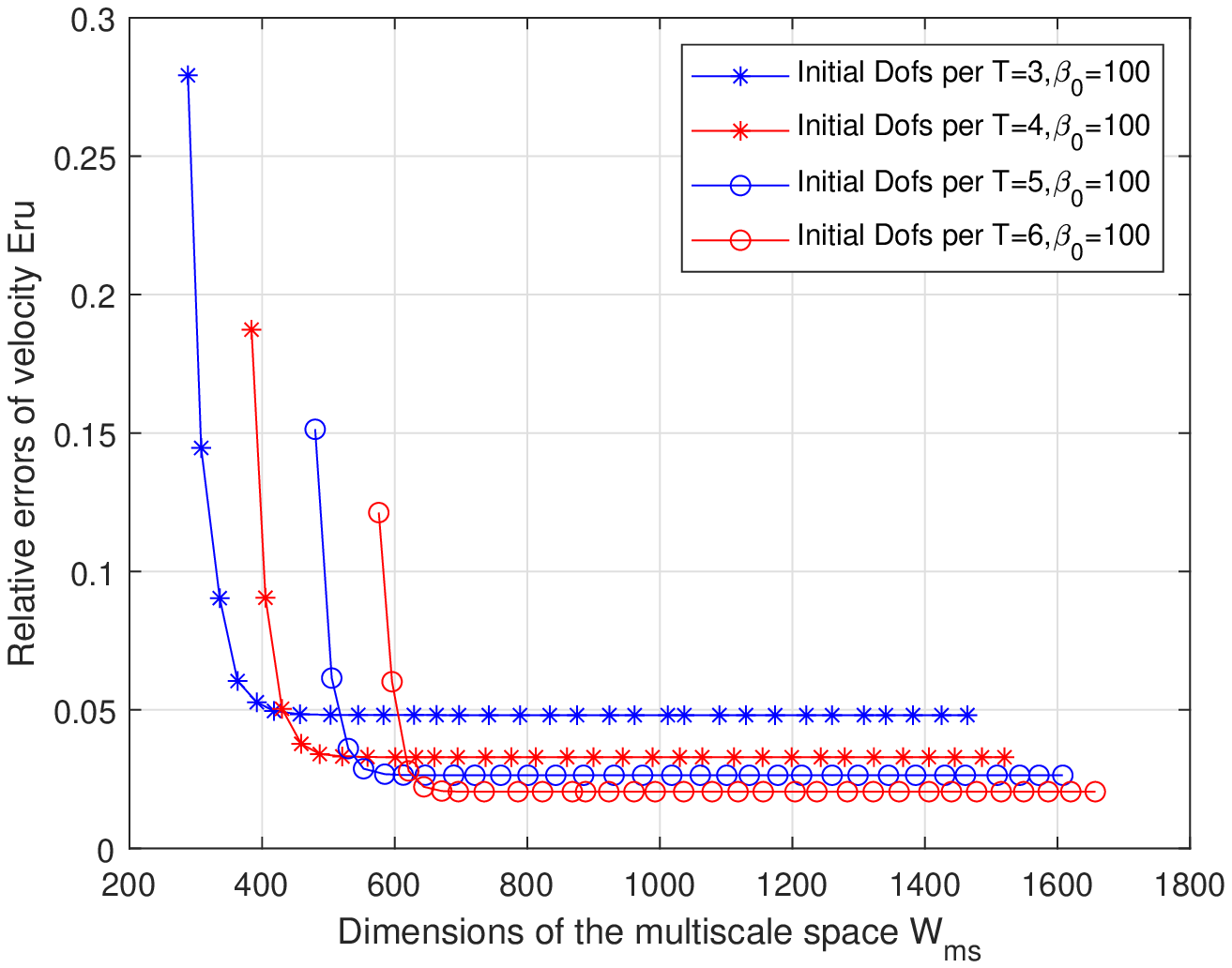}
	\end{minipage}
	\\
	\mbox{\hspace{0.00cm}}
	\begin{minipage}[b]{0.40\textwidth}
		\centering
		\includegraphics[scale=0.48]{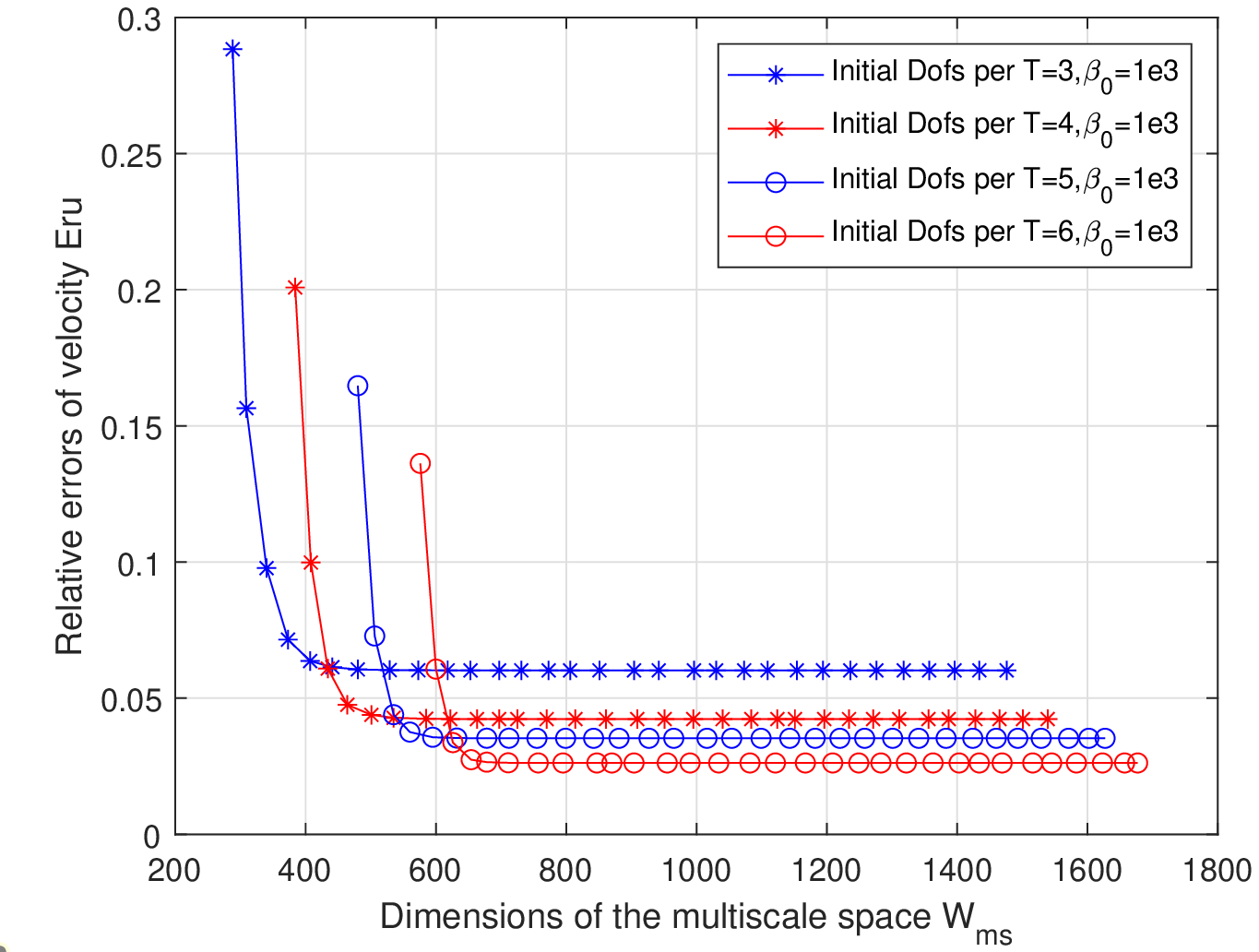}
	\end{minipage}
	\mbox{\hspace{0.00cm}}
	\begin{minipage}[b]{0.40\textwidth}
		\centering
		\includegraphics[scale=0.48]{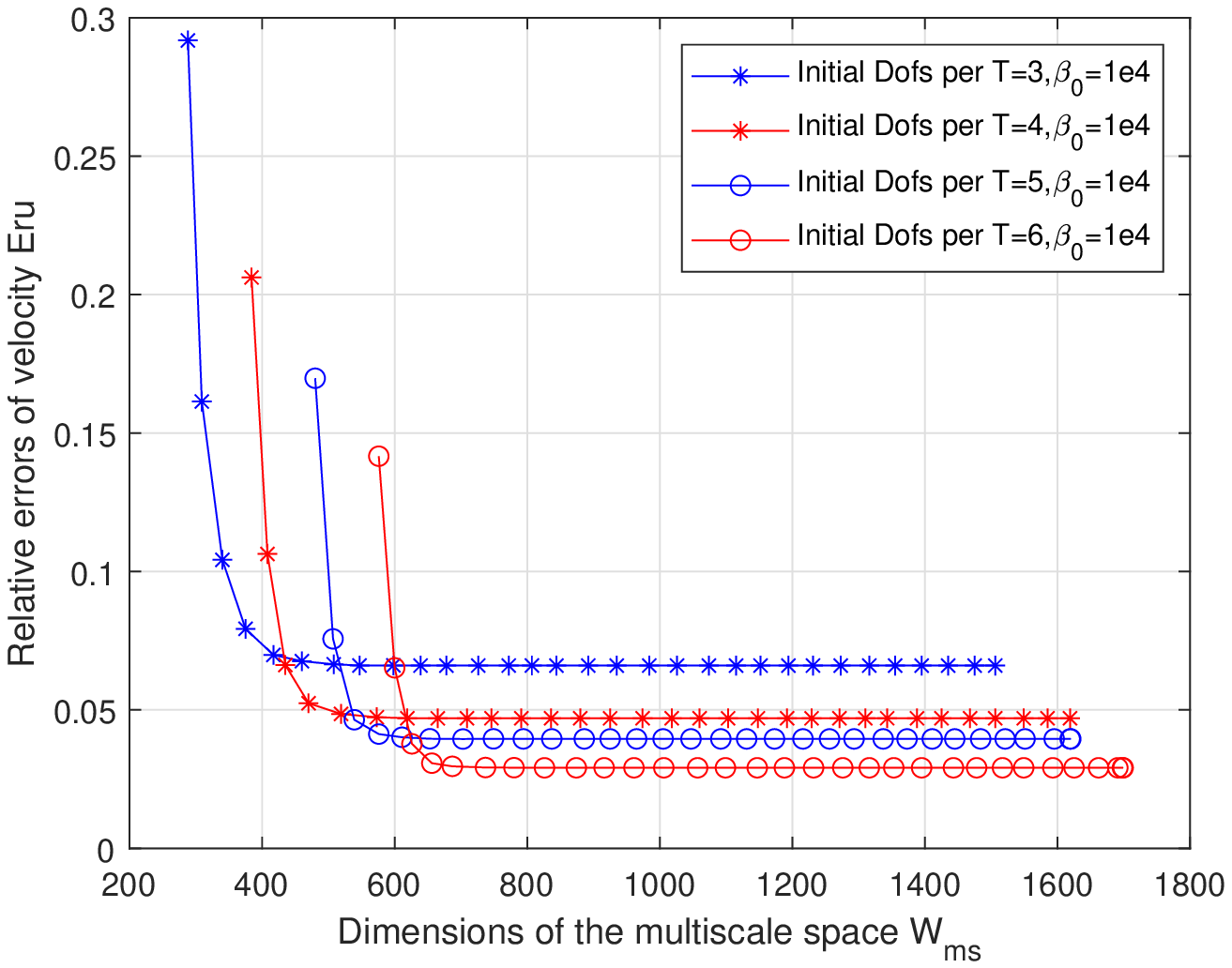}
	\end{minipage}
	\caption{Adaptive online enrichment (Remark 2) : Convergence comparison for different choices of the number of initial basis. Top left: $\beta_0=1$. Top right: $\beta_0=10$. Middle left: $\beta_0=100$. Middle right: $\beta_0=1e3$. Bottom: $\beta_0 = 1e4$.}\label{fig_velocityerr_4_rem2}
\end{figure}
\section{Conclusions}
In this paper, we employ the GMsFEM framework to solve the Darcy-Forchheimer model in highly heterogeneous porous media. An MFMFE method is applied for the discretization of the problem on the underlying fine grid. In the MFMFE method, $\textrm{BDM}_1$ mixed finite element spaces are used for approximating the velocity and pressure, and symmetric trapezoidal quadrature rule is employed for the integration of bilinear forms relating to the velocity variables, 
which allows for local velocity elimination and lead to a cell-centered system for the pressure. 
We construct the multiscale basis functions for approximating the pressure and solve the problem on the coarse grid following the GMsFEM framework. The computation of the local snapshot spaces and the smaller dimensional offline space by a series of local spectral decompositions are conducted in the offline stage. In the online stage, we use the Newton iterative algorithm to deal with the nonlinear term and obtain the offline solution, then based on the offline space and offline solution, we enrich the multiscale space by calculating online basis functions. 
In the end, some numerical examples are supplied to test the performance of the proposed multiscale method. The numerical results demonstrate that the number of Newton iterations is much less than the Picard iterations, the multiscale method provides a good approximation of the problem on the coarse grid even though the Darcy-Forchheimer parameter takes large values and the online basis functions are effective to improve the accuracy of the multiscale solution substantially.
\bibliographystyle{elsarticle-num}
\bibliography{Darcy-Forchheimer_REF}

\begin{thebibliography}{10}
\expandafter\ifx\csname url\endcsname\relax
  \def\url#1{\texttt{#1}}\fi
\expandafter\ifx\csname urlprefix\endcsname\relax\def\urlprefix{URL }\fi
\expandafter\ifx\csname href\endcsname\relax
  \def\href#1#2{#2} \def\path#1{#1}\fi

\bibitem{park2005mixed}
E.-J. Park, Mixed finite element methods for generalized {Forchheimer} flow in
  porous media, Numerical Methods for Partial Differential Equations: An
  International Journal 21~(2) (2005) 213--228.

\bibitem{girault2008numerical}
V.~Girault, M.~F. Wheeler, Numerical discretization of a {Darcy--Forchheimer}
  model, Numerische Mathematik 110~(2) (2008) 161--198.

\bibitem{pan2012mixed}
H.~Pan, H.~Rui, Mixed element method for two-dimensional {Darcy-Forchheimer}
  model, Journal of Scientific Computing 52~(3) (2012) 563--587.

\bibitem{rui2012block}
H.~Rui, H.~Pan, A block-centered finite difference method for the
  {Darcy--Forchheimer} model, SIAM Journal on Numerical Analysis 50~(5) (2012)
  2612--2631.

\bibitem{rui2015block}
H.~Rui, D.~Zhao, H.~Pan, A block-centered finite difference method for
  {Darcy--Forchheimer} model with variable forchheimer number, Numerical
  Methods for Partial Differential Equations 31~(5) (2015) 1603--1622.

\bibitem{rui2017block}
H.~Rui, H.~Pan, A block-centered finite difference method for slightly
  compressible {Darcy--Forchheimer} flow in porous media, Journal of Scientific
  Computing 73~(1) (2017) 70--92.

\bibitem{wang2015stabilized}
Y.~Wang, H.~Rui, Stabilized crouzeix--raviart element for {Darcy--Forchheimer}
  model, Numerical Methods for Partial Differential Equations 31~(5) (2015)
  1568--1588.

\bibitem{xu2017multipoint}
W.~Xu, D.~Liang, H.~Rui, A multipoint flux mixed finite element method for the
  compressible {Darcy--Forchheimer} models, Applied Mathematics and Computation
  315 (2017) 259--277.

\bibitem{huang2018multigrid}
J.~Huang, L.~Chen, H.~Rui, Multigrid methods for a mixed finite element method
  of the {Darcy--Forchheimer} model, Journal of scientific computing 74~(1)
  (2018) 396--411.

\bibitem{rui2015two}
H.~Rui, W.~Liu, A two-grid block-centered finite difference method for
  {Darcy--Forchheimer} flow in porous media, SIAM Journal on Numerical Analysis
  53~(4) (2015) 1941--1962.

\bibitem{fairag2020two}
F.~A. Fairag, J.~D. Audu, Two-level galerkin mixed finite element method for
  {Darcy--Forchheimer} model in porous media, SIAM Journal on Numerical
  Analysis 58~(1) (2020) 234--253.

\bibitem{zhang2020variational}
T.~Zhang, X.~Li, Variational multiscale interpolating element-free galerkin
  method for the nonlinear darcy--forchheimer model, Computers \& Mathematics
  with Applications 79~(2) (2020) 363--377.

\bibitem{spiridonov2019mixed}
D.~Spiridonov, J.~Huang, M.~Vasilyeva, Y.~Huang, E.~T. Chung, Mixed generalized
  multiscale finite element method for {Darcy-Forchheimer} model, Mathematics
  7~(12) (2019) 1212.

\bibitem{chung2015mixed}
E.~T. Chung, Y.~Efendiev, C.~S. Lee, Mixed generalized multiscale finite
  element methods and applications, Multiscale Modeling \& Simulation 13~(1)
  (2015) 338--366.

\bibitem{efendiev2013generalized}
Y.~Efendiev, J.~Galvis, T.~Y. Hou, Generalized multiscale finite element
  methods {(GMsFEM)}, Journal of Computational Physics 251 (2013) 116--135.

\bibitem{chung2014adaptive}
E.~T. Chung, Y.~Efendiev, G.~Li, An adaptive gmsfem for high-contrast flow
  problems, Journal of Computational Physics 273 (2014) 54--76.

\bibitem{efendiev2014generalized}
Y.~Efendiev, J.~Galvis, G.~Li, M.~Presho, Generalized multiscale finite element
  methods. nonlinear elliptic equations, Communications in Computational
  Physics 15~(3) (2014) 733--755.

\bibitem{chung2015residual}
E.~T. Chung, Y.~Efendiev, W.~T. Leung, Residual-driven online generalized
  multiscale finite element methods, Journal of Computational Physics 302
  (2015) 176--190.

\bibitem{chan2016adaptive}
H.~Y. Chan, E.~Chung, Y.~Efendiev, Adaptive mixed gmsfem for flows in
  heterogeneous media, Numerical Mathematics: Theory, Methods and Applications
  9~(4) (2016) 497--527.

\bibitem{chung2016adaptive}
E.~Chung, Y.~Efendiev, T.~Y. Hou, Adaptive multiscale model reduction with
  generalized multiscale finite element methods, Journal of Computational
  Physics 320 (2016) 69--95.

\bibitem{hou1997multiscale}
T.~Y. Hou, X.-H. Wu, A multiscale finite element method for elliptic problems
  in composite materials and porous media, Journal of computational physics
  134~(1) (1997) 169--189.

\bibitem{chen2020generalized}
J.~Chen, E.~T. Chung, Z.~He, S.~Sun, Generalized multiscale approximation of
  mixed finite elements with velocity elimination for subsurface flow, Journal
  of Computational Physics 404 (2020) 109133.

\bibitem{wheeler2006multipoint}
M.~F. Wheeler, I.~Yotov, A multipoint flux mixed finite element method, SIAM
  Journal on Numerical Analysis 44~(5) (2006) 2082--2106.

\bibitem{wheeler2012multipoint}
M.~Wheeler, G.~Xue, I.~Yotov, A multipoint flux mixed finite element method on
  distorted quadrilaterals and hexahedra, Numerische Mathematik 121~(1) (2012)
  165--204.

\bibitem{ingram2010multipoint}
R.~Ingram, M.~F. Wheeler, I.~Yotov, A multipoint flux mixed finite element
  method on hexahedra, SIAM Journal on Numerical Analysis 48~(4) (2010)
  1281--1312.

\bibitem{klausen2006robust}
R.~A. Klausen, R.~Winther, Robust convergence of multi point flux approximation
  on rough grids, Numerische Mathematik 104~(3) (2006) 317--337.

\bibitem{hou2014numerical}
J.~Hou, S.~Sun, Z.~Chen, Numerical comparison of robustness of some reduction
  methods in rough grids, Numerical Methods for Partial Differential Equations
  30~(5) (2014) 1484--1506.

\bibitem{chung2017online}
E.~T. Chung, Y.~Efendiev, W.~T. Leung, An online generalized multiscale
  discontinuous {Galerkin} method {(GMsDGM)} for flows in heterogeneous media,
  Communications in Computational Physics 21~(2) (2017) 401--422.

\bibitem{li2001literature}
D.~Li, T.~W. Engler, et~al., Literature review on correlations of the
  {non-Darcy} coefficient, in: SPE Permian Basin Oil and Gas Recovery
  Conference, Society of Petroleum Engineers, 2001.

\bibitem{muljadi2016impact}
B.~P. Muljadi, M.~J. Blunt, A.~Q. Raeini, B.~Bijeljic, The impact of porous
  media heterogeneity on {non-Darcy} flow behaviour from pore-scale simulation,
  Advances in Water Resources 95 (2016) 329--340.

\end{thebibliography}
\end{document}